\documentclass[12pt,letterpaper,reqno]{amsart}

\usepackage{amssymb}
\usepackage{amsmath}
\usepackage{amsthm}
\usepackage{xcolor}
\usepackage{graphicx}
\graphicspath{ {pics/} }

\usepackage{caption}
\usepackage{subcaption}

\usepackage{wrapfig}
\usepackage{float}
\usepackage{floatflt}

\usepackage{cancel}
\usepackage{comment}

\usepackage{mathtools} 

\usepackage{pgfplots}
\pgfplotsset{compat=1.18} 

\newcommand{\R}{\mathbb R}

\newcommand{\sech}{\textmd{sech}}

\newcommand{\ds}{\displaystyle}

\newtheorem{theorem}{Theorem}[section]
\newtheorem{proposition}[theorem]{Proposition}
\newtheorem{lemma}[theorem]{Lemma}

\theoremstyle{definition}

\newtheorem{conj}[theorem]{Conjecture}

\theoremstyle{remark}
\newtheorem{remark}[theorem]{Remark}

\title[NLS on a unit ball]{Nonlinear Schr\"odinger equation on a unit ball\\ in one and two dimensions}

\author[C. Klein]{Christian Klein}
\address{Université Bourgogne Europe, CNRS, IMB UMR 5584, 21000 Dijon, France, and 
Institut Universitaire de France}
\email{Christian.Klein@u-bourgogne.fr}

\author[S. Roudenko]{Svetlana Roudenko}
\address{Department of Mathematics \& Statistics\\
Florida International University,  Miami, FL 33199, USA}
\curraddr{}
\email{sroudenko@fiu.edu}

\author[N. Stoilov]{Nikola Stoilov}
\address{Université Bourgogne Europe, CNRS, IMB UMR 5584, 21000 Dijon, France} 
\email{Nikola.Stoilov@u-bourgogne.fr}



\begin{document}

\begin{abstract}

We investigate the nonlinear Schr\"odinger equation on a unit ball in one and two dimensions with Dirichlet boundary condition, focusing on ground states, their limiting behavior, stability and perturbation dynamics.
Our computations illustrate the stabilizing effect of the Dirichlet boundary compared with the whole-space setting. 
In the subcritical and critical cases, the mass curves have positive slope and the considered perturbations show recurrent dynamics near the ground state family. 
In the critical case, the mass of the ground states increases towards the mass of the whole-space ground state. Staying below the whole-space ground state mass makes perturbations of the bounded-domain ground states stable. 
In the supercritical case, the computations reveal a {\it new} branching phenomenon: the mass curve develops a turning point, separating stable and unstable branches. In 1D we prove existence of a global maximum in the mass curve and derive an integral equation for the turning points. 
Perturbations of the stable branch remain close to the corresponding ground state, while perturbations of the unstable branch lead to finite-time blow-up or recurrent oscillations between neighborhoods of at least two coherent profiles. 
We prove sufficient blow-up conditions in the critical and supercritical cases, in particular, for every fixed $A>1$, the data $A Q_b$ blows up for sufficiently large $b$ in these cases.  
Furthermore, the Dirichlet boundary reflects the solution, preventing mass from escaping via radiation, while coherent oscillations persist even for small initial data, thus, supporting a bounded-domain soliton-resolution conjecture.  
\end{abstract}

\maketitle


\section{Introduction}

We consider the nonlinear Schr\"odinger (NLS) equation  on a bounded domain: 
\begin{equation}\label{NLS}
i\,u_t + \Delta u + |u|^{\alpha} u = 0, \qquad x \in \Omega \subset \R^d,  
\quad t \in \R, \qquad
\end{equation}
where $u(t,x)$ is a complex-valued function, $\Delta$ is the standard Laplacian, 
the nonlinearity power is $\alpha>0$, and $\Omega$ is a bounded domain, 
which for this paper we take as a unit ball $B_1 = \{x \in \mathbb R^d : |x| < 1 \}$ in 1 or 2 dimensions.  
We are specifically interested in the Dirichlet boundary condition, namely,
\begin{equation}\label{BC-Omega}
u(t,x) = 0, \quad t \in \mathbb R, \quad x \in \partial \Omega,
\end{equation}
and its influence on the solution, including solitary waves and their stability.

The Cauchy problem for the equation \eqref{NLS} is locally well-posed 
in $H^1_0$: for $u_0 \in H_0^1(B_1)$ 
there exists $T= T (\|u_0\|_{H^1}) >0$ and 
a unique solution $u \in C([0,T), H_0^1(B_1))$ with the given initial condition $u_0$. 
The existence and continuous dependence (as well as mass and energy conservations) follow from \cite[Thm. 3.3.9]{Caz-book}. 
The uniqueness is addressed in dimension $d=1$ for any $\alpha$ in \cite[Thm. 3.5.1]{Caz-book},
in dimension $d = 2$ when $0 < \alpha \leq 2$ in \cite{V1984}, \cite{O1990}, 
and when $2<\alpha<\infty$, uniqueness holds among generalized
solutions in $L^\infty([0,T), H_0^1(B_1))$ satisfying the additional
condition ~${\rm ess\,sup}_{0<t<T}
\int_{B_1}\exp\!\left(\kappa|u(t,x)|^\alpha\right)\,dx<\infty$ ~
for some $\kappa>0$, as in \cite[Thm.~5]{V1986}\footnote{This condition holds on every compact time interval, on which the solution is bounded in $L^\infty(B_1)$. 
Uniqueness in the bounded solution class also follows directly
from \cite[Thm.~4]{V1986}.}.

As a consequence, during their lifespan, solutions $u(t)$ to \eqref{NLS} conserve mass and energy (or Hamiltonian):
\begin{equation}\label{MC}
M[u(t)]=\int_{B_1} |u(t)|^2 \, dx \equiv M[u(0)]
\end{equation}
and
\begin{equation}\label{EC}
E[u(t)]=\dfrac{1}{2}\int_{B_1} |\nabla u(t)|^2 \; dx - \dfrac1{\alpha+2} \int_{B_1} |u(t)|^{\alpha+2} \; dx \equiv E[u(0)].
\end{equation}

While the scaling invariance is absent in the equation \eqref{NLS} due to the boundary, one can still refer to the criticality of the equation according to the scaling of the NLS equation on the whole space:
$$
s_c = \frac{d}2 - \frac{2}{\alpha}.
$$
Thus, the NLS equation is ($L^2$ or mass)-critical if $\alpha=\frac4{d}$, subcritical when $\alpha<\frac4{d}$ and supercritical otherwise. 

From the conservation laws and the sharp Gagliardo--Nirenberg inequality on $H^1_0(B_1)$, solutions with $u_0\in H^1_0(B_1)$ are global in the subcritical case. In the critical case, a standard sufficient condition is that the mass lies below the sharp Gagliardo--Nirenberg threshold on $B_1$.
This threshold coincides with $M(\mathcal R)$, where $\mathcal R$ is the whole-space ground state, $M(u_0)<M(\mathcal R)$.
For the Dirichlet ground-state branch $Q_b$, this value is approached in the limit $b\to\infty$, namely, $M(Q_b)\nearrow M(\mathcal R)$, which we discuss in Lemma \ref{L:conv}(II).

\subsection{Motivation}
The cubic NLS on the unit ball in 2D with Dirichlet boundary condition
describes the propagation of the laser beam in a hollow-core fiber 
(as opposed to a solid-core, thus, reducing drastically light-matter interaction).
A hollow waveguide filled with a noble gas is used as an effective technique for extending the
interaction length between nonlinear optical materials and high-energy laser pulses.
The electric field outside the core is negligible and 
almost all of the laser beam reflects on the interface between the core and the cladding, implying the Dirichlet
boundary condition \eqref{BC-Omega}, \cite{GF2000}, \cite{FHK2012}, \cite{TB1998}. 
This setting influences the stability of ground states (and excited states) in the NLS on a unit ball 
as was indicated in \cite{FM2001} and \cite{FHK2012}, 
which differs from the NLS considered on the whole space $\R^d$. 
The purpose of this paper is exactly to investigate further 
solution behavior and dynamics of ground states in this model and the difference with the whole space.

\subsection{Solitary waves in a ball}\label{S:0.2}
Let $\alpha>0$ for $d=1,2$ and $0<\alpha<4/(d-2)$ for $d \geq 3$.

The NLS equation \eqref{NLS} on a bounded domain has a family of standing solitary waves, 
\begin{equation}\label{Eq:Qb}
u(t,x) = e^{ibt} \, Q_b(x), 
\end{equation}
with a parameter $b \in \mathbb R$ and 
$Q_b(x)$ a {\it ground state} solution  
in $H^1_0(B_1)$ with the  Dirichlet boundary condition on $B_1$, 
i.e., 
\begin{equation}\label{E:GS-NLS}
\left\{ 
\begin{array}{l}
-\Delta Q +  b\, Q  - |Q|^{\alpha} Q = 0, \quad x \in B_1,  \\
\quad Q(x) = 0, \qquad x \in \partial B_1.
\end{array} 
\right.
\end{equation}
(We note that on the whole space $\mathbb R^d$ ground states  $Q_b$ are defined only for positive values of $b$, $b>0$; also $Q_b(x) \to 0$, decaying exponentially, as $|x| \to \infty$.)

In this work we define the ground states variationally, following Fukuizumi et al. \cite{FHK2012}. 
For $b \in \mathbb R$ define an action (often called the Weinstein or Lyapunov functional in the dispersive literature) 
$$
S_b (u) = E(u) +\frac{b}2 M(u).
$$ 
Then $S_b'(u) = 0$ if and only if $u = u_b$ 
$\in H^1_0(B_1)$ is a weak solution of \eqref{E:GS-NLS} for a given $b$. 
The {\it ground state} solution $Q_b \in H^1_0(B_1)$ for a given $b$ is 
the {\it least-energy} solution of \eqref{E:GS-NLS} if $Q_b$ is a non-trivial weak solution of \eqref{E:GS-NLS} 
and it minimizes $S_b$ among all non-trivial solutions of \eqref{E:GS-NLS}. 
The least-energy solution exists for $b>-\lambda_1$, where 
$\lambda_1>0$ is the {\it first (lowest) eigenvalue} of the operator $(-\Delta)$ on $B_1$ with Dirichlet boundary condition; this solution is unique (up to a constant phase), radially symmetric, positive after a phase rotation, and strictly decreasing in the radial variable $r = |x|$, see \cite{FHK2012} and references therein. 
\smallskip

The standing wave $e^{i b t}Q_b$ is said to be {\it orbitally stable} in $H_0^1(B_1)$ if for any $\varepsilon >0$ there exists $\delta > 0$ such that whenever $u_0 \in H_0^1(B_1)$ with $\|u_0-Q_b\|_{H^1} < \delta$, the solution $u(t)$ of \eqref{NLS} with $u(0) = u_0$ exists globally and satisfies  
$$
\inf_{\theta \in \mathbb R} \|u(t) - e^{i \theta} Q_b\|_{H^1} <\varepsilon \qquad \mbox{for ~~all} \quad t \geq 0.
$$
Otherwise, the standing wave is {\it unstable}. We mention that when arbitrarily small perturbations produce finite-time blow-up, it is often referred to as strong instability. 
\smallskip

It was found in \cite{FM2001} ($d=2, \alpha=2$) and \cite{FHK2012} (general dimension and power) that the Dirichlet boundary has a {\it stabilizing} effect: 
standing solitary waves $e^{ibt}Q_b$ of the equation \eqref{NLS} are orbitally stable in $H_0^1(B_1)$ not only in the subcritical case $0 < \alpha < \frac4{d}$ (for the entire range of $b$), but also for some range of $b$ in the critical case $\alpha = \frac4{d}$. 
Furthermore, in the supercritical case, the ground state is also stable in a small neighborhood of $-\lambda_1$, which shows a significant difference with the ground state  behavior in the critical case 
on the whole space $\mathbb R^d$ as well as in the supercritical case (unstable for any $b$). 
We summarize the results on orbital stability and instability of ground states in this setting in the following statement (denoting by $\frac4{(d-2)^+}$ the value $\frac{4}{d-2}$ if $d \geq 3$ and $\infty$ for $d=1,2$).

\begin{theorem}\label{T:1}
Let $d \geq 1$ and $Q_b$ be the least-energy (ground state) solution. 
\begin{itemize}
\item[\bf(I)] 
\underline{\rm (Orbital stability)}

\noindent{\rm(a)} For $0 < \alpha \leq \frac4{d}$, the standing wave $e^{ib t } Q_b$ is orbitally stable in $H^1_0 (B_1)$ for any $b \in (-\lambda_1,\infty)$.

\noindent{\rm(b)} For $\frac4{d} < \alpha < \frac4{(d-2)^+}$, 
there exists $\epsilon > 0$ such that the standing wave $e^{ibt}Q_b$ is orbitally stable in 
$H_0^1(B_1)$ for any $b \in (-\lambda_1, -\lambda_1 + \epsilon) $.
\smallskip

\item[\bf(II)] 
\underline{\rm (Instability)} 

\noindent Let $\frac4{d} < \alpha < \frac4{(d-2)^+}$. There exists $b_2 > 0$ such that the standing wave $e^{i b t} Q_b$ is unstable for any $b \in (b_2,\infty)$.
\end{itemize}
\end{theorem}

The proof of Theorem \ref{T:1} is contained in several works. Initially, the stability problem was considered for the 2D critical NLS by Fibich \& Merle in \cite{FM2001}, where  orbital stability of standing waves was proved close to the endpoints, namely,  
for $b \in (-\lambda_1, -\lambda_1 + \epsilon) \cup (b_1, \infty)$, for some $b_1 > 0$ and $\epsilon > 0$. This was intriguing, since it exhibited stable regimes in a critical problem.
This result was extended to all dimensions by Fukuizumi et al. in \cite{FHK2012}; 
in particular, the 1D critical case was treated on the full interval $(-\lambda_1,\infty)$.
The extension to higher dimensions in the critical case was obtained by a different constrained-minimization approach with {\it two} constraints in \cite[Thm. 1.8]{NTV2014}. (A recent work \cite{FIK2025} gives a simplified one-dimensional approach via a standard one-constraint minimization obtaining a similar result in one dimension and discusses possible extensions, see \cite[Remark 1.4]{FIK2025}.)
The stability part in an intercritical regime slightly above $-\lambda_1$ and the instability in Theorem \ref{T:1} are also shown in \cite{FHK2012}, \cite{FM2001}, \cite{NTV2014}.

While Theorem \ref{T:1} gives important analytical stability and instability information, it does not compute the constants $\epsilon$ and $b_2$, nor does it  determine whether the stable and unstable intervals meet at a single transition point (or there may be several transition points), nor does it describe the {\it finite-amplitude} (vs. infinitely-small) perturbation dynamics. These are the questions that we address in the present work. 
\smallskip

In this paper, we investigate standing waves in dimensions $d=1, 2$:

$\bullet$ In the subcritical and critical cases, we examine stability regimes across the entire interval $(-\lambda_1, \infty)$ under various perturbations (including finite-amplitude). In particular, we confirm Theorem \ref{T:1}(I) numerically and complement it by showing the resulting perturbation dynamics.

$\bullet$ In the supercritical case, we identify complementary ranges of stability and instability and determine the threshold values $b^*$, where the stability transition occurs. For the one-dimensional problem, we prove Proposition \ref{P:1D-b*}, which provides an exact integral characterization of the mass turning points. Taking the ideas of Kfoury-LeCoz-Tsai in \cite{KLT}, we study the mass and the length integrals, and prove the existence of the maximum for the mass curve, giving the exact formulas for obtaining its location (namely, at the turning point $b^*$) and value $M(Q_{b^*})$. This characterizes the mass turning points (though, we do not address their uniqueness).
This implies that in the supercritical case, {\it branching} of ground states occurs, which we demonstrate in our computations, leading to the emergence of stable and unstable branches (similar to the branching observed for the bi-harmonic NLS in \cite{KPRS} and \cite{KR2026}; see also combined NLS in \cite{CKS} and \cite{HRR2026}).

$\bullet$ In Proposition \ref{P:blowup} we give sufficient blow-up conditions from the Dirichlet virial identity and either negative energy or the scale-invariant mass-energy dichotomy argument of Holmer-Roudenko \cite{HR2007, HR2008}. In particular, for every $A>1$, the data $A\, Q_b$ blow up for all sufficiently large $b$ in the mass-critical and intercritical cases (the statement is also extended for more general profiles $f$).
\smallskip

To investigate orbital stability, one way is to prove it analytically by verifying the slope condition, namely, the sign of $\partial_b M(Q_b)$ (the Vakhitov--Kolokolov condition in the Grillakis--Shatah--Strauss framework), see Proposition 5 in \cite{FHK2012} and the general theory \cite{GSS}. In the present paper we use this slope as a numerical diagnostic and combine it with direct time evolutions of perturbed ground states. Our contribution is the computation of the ground states in this setting, the mass-energy curves, the turning points, and the associated dynamics in the ranges considered below, including analytical proof of turning points in 1D (see comments about 2D in Remark \ref{R:2D-turning}) and blow-up criteria that rigorously establish a large-$b$ regime of increasing amplitude blow-up. 
 
We would like to point out the following difference in stability investigations: the analytical orbital-stability statements, quoted above, describe {\it infinitesimally} small perturbations: for each $\epsilon > 0$ there exists a $\delta > 0$ such that $\delta$-small $H^1$-perturbations stay in an $\epsilon$-neighborhood of the orbit. They do not describe the behavior of {\it finite-amplitude} perturbations,
nor the dynamics after a solution leaves a small neighborhood of the initial ground state. The simulations below are intended to probe these finite perturbations and their resulting dynamics. 
In fact, we provide numerical evidence as well as some analysis towards the following conjecture. (We note that the finite-amplitude stability assertions below concern perturbations within ranges that may depend on the ground state, and thus, determining these ranges is part of the conjecture.)

\begin{figure}[htb!]
\includegraphics[width=0.52\textwidth,height=0.38\textwidth]{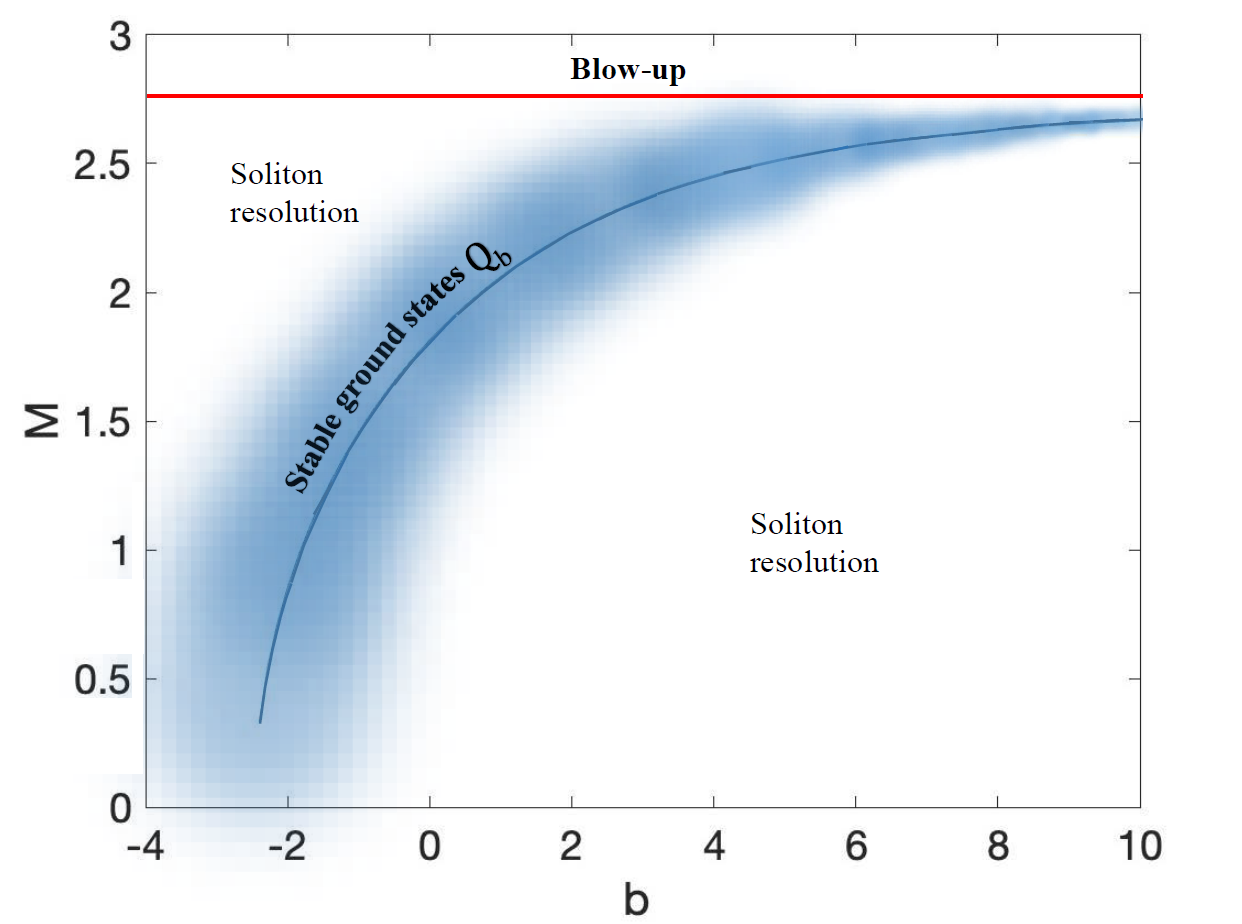}
   \includegraphics[width=0.47\textwidth,height=0.38\textwidth]{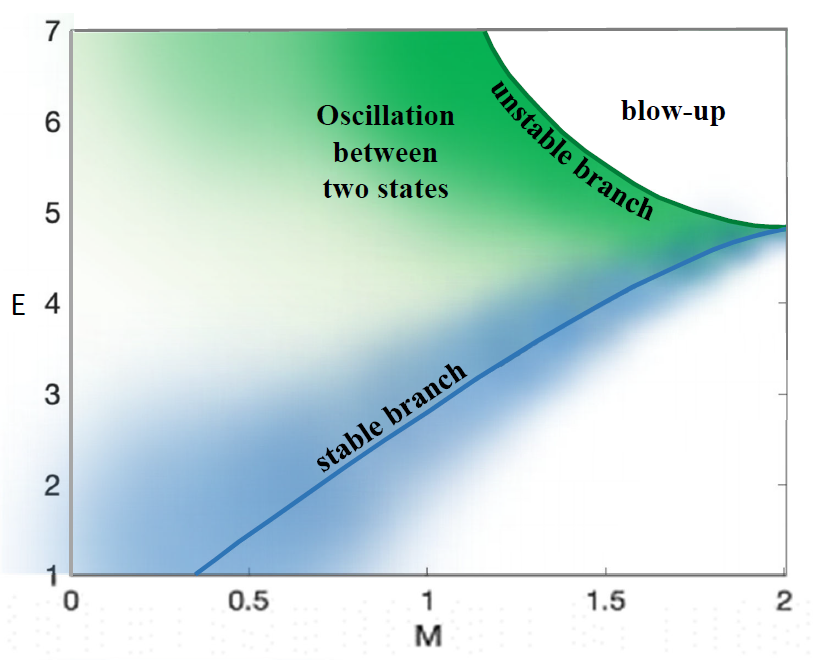}
\caption{\small Schematic stability of ground states $Q_b$ in critical (left) and supercritical (right) cases. 
The axes labels on the left are from the example of 1D quintic NLS, Figure \ref{F:ME-1D-p3-5}(D) 
and on the right are from the example of 2D quintic NLS, Figure \ref{F:ME-2D}(I).}
\label{F:cloud}
\end{figure}

\begin{conj}\label{C:1}
Let $d \geq 1$ and let $Q_b$ be the least-energy ground state solution of \eqref{E:GS-NLS} for $b \in (-\lambda_1, \infty)$.  

\smallskip

\noindent\underline{\textbf{Part I: Subcritical case.}} Let $0<\alpha<4/d$. 
For any $b \in (-\lambda_1,\infty)$ the slope condition has a positive sign, ${\partial_b} \,  M(Q_b) >0$. This implies that the standing wave $e^{ibt} Q_b$ is orbitally stable in $H^1_0(B_1)$. 

Furthermore, even sufficiently large initial perturbations remain close to the ground-state family and exhibit only small recurrent oscillations (and converge to oscillations around the final state with sufficiently small amplitude difference).

\smallskip

\noindent\underline{\textbf{Part II: Critical case.}} Let $\alpha=4/d$.
Then 
\begin{itemize}
\item[(i)]
for any $b \in (-\lambda_1,\infty)$,  $M(Q_b) <  M(\mathcal R)$ and the slope condition has a positive sign, ${\partial_b} \,  M(Q_b) >0$. This implies that the standing wave $e^{ibt} Q_b$ is orbitally stable in $H^1_0(B_1)$, including sufficiently large perturbations as long as $M(u_0) <M(\mathcal R)$. 

\noindent Furthermore, $ M(Q_b)$ increases monotonically with  $\ds \lim_{b \to \infty} M(Q_b) =  M(\mathcal R)$.

\noindent{\it (see Figure \ref{F:cloud} left: the blue line represents $M(Q_b)$, the area around it shows stable perturbations)}.  

\item[(ii)]
If $M(u_0) > M(\mathcal R)$, then the solution $u(t)$ with the initial condition $u_0$ and negative energy blows up in finite time. 
\noindent{\it (see Figure \ref{F:cloud} left: the red line represents the mass of the ground state on the whole space, $M(\mathcal{R})$, above which initial data with negative energy blow up in finite time)}. 
\end{itemize}

\smallskip

\noindent\underline{\textbf{Part III: Supercritical case.}} Let $\alpha > \frac4{d}$ (with $\alpha < \frac4{d-2}$ if $d \geq 3$). 
Then there exists a single turning point $b^\ast \in (-\lambda_1,\infty)$ such that 
\begin{itemize}
\item[(i)] ({\it existence of $b^*$})
the graph of $M(Q_b)$ as a function of $b$ is a single-peaked curve, which increases monotonically up to its maximum value at $b^\ast$ and then decreases to zero \\
{\it (for example, as in Figure \ref{F:ME-2D} (G))};

\item[(ii)] ({\it stable branch})
for any $b \in (-\lambda_1, b^\ast)~~$ the slope $~~{\partial_b} \,  M(Q_b) >0$. This implies that the standing wave $e^{ibt} Q_b$ is orbitally stable in $H^1_0(B_1)$, including  sufficiently large perturbations 
{\it (see Figure \ref{F:cloud} right: the blue curve and the shaded area around it, and also Remark \ref{R:1})}. 

\item[(iii)] ({\it unstable branch})
for any $b \in (b^\ast, \infty)~$ the slope $~~{\partial_b} \,  M(Q_b) <0$, and thus, the standing wave $e^{ibt} Q_b$ is unstable.  The mechanism of the instability is as follows:
\begin{itemize}
\item
perturbations of larger mass than the ground state ($A>1$ in $u_0 = AQ_{b}$) produce blow-up in finite time;  
\item
perturbations of smaller mass than the ground state ($A<1$ in $u_0=AQ_b$) approach oscillation  between two states 
\end{itemize} 
{\it (see Figure \ref{F:cloud} right: the green curve and shaded area {\rm below} it, and also  Remark \ref{R:1})}. 
\end{itemize}
\end{conj}

\begin{remark}\label{R:1}
In Figure \ref{F:cloud} we schematically depict the stability of ground states in 
the critical and supercritical cases as stated in Conjecture \ref{C:1} Parts II and III. 

On the left, using the 1D quintic NLS example from Figure \ref{F:ME-1D-p3-5}(D), we show the stability of $Q_b$ for different $b$ in the {\it critical} case: the blue curve indicates the mass of $Q_b$. Perturbations of $Q_b$ are stable, which is indicated by a blue cloud area around the blue curve. Note that the shaded perturbation zone shrinks as $b \to \infty$ (thus, allowing smaller perturbations), since $M(Q_b) \to M(\mathcal R)$ (the red line), above which solutions with negative energy (and possibly some with positive energy) blow up in finite time. In the rest of the plot, solutions are conjectured to satisfy the {\it soliton resolution}, splitting into a sum of coherent structures (such as solitons), even for the small initial data, as the mass of $Q_b$ covers the entire interval $(0,M(\mathcal{R}))$ (with possible small oscillations around the final states). There seems to be no radiation forming in this problem: the Dirichlet boundary condition prevents mass from escaping, instead it has a stabilizing effect and reflects solutions completely, and thus, the reflected dynamics exhibits oscillations around the coherent states in our simulations. 

On the right, using the 2D quintic NLS example from Figure \ref{F:ME-2D} (I), we show the {\bf branching} of ground states in the {\it supercritical} case: the blue curve indicates the stable branch of $Q_b$ and the green curve shows the unstable branch of $Q_b$. Perturbations of the stable branch of $Q_b$ (blue curve) are stable, as indicated by the shaded blue area around the blue curve. Perturbations of the unstable branch either produce a blow-up in finite time (if $A>1$ in the perturbations $u_0 = AQ_{b}$) or oscillations between two states (if $A<1$), shaded in green below the green curve. 
\end{remark}

\begin{remark}\label{R:2} 
Here, we point out which parts of the conjecture have been proved analytically and to which extent, in connection to the above Theorem \ref{T:1}. 
In dimension 1, Part I of Conjecture \ref{C:1} is proved (for $\delta$-small perturbations) in Proposition 20 of \cite{FHK2012}, which is also a part of Theorem \ref{T:1}(I.a). In dimensions two and higher Theorems 17, 18 of \cite{FHK2012} address the orbital stability (again for small $\delta$-perturbations) either in the neighborhood of $-\lambda_1$ or for sufficiently large $b$; the entire range is treated in Theorem 1.8 in \cite{NTV2014}.
The slope condition in Part II, in dimension 2, was partially proved in \cite{FM2001}, and in all dimensions is partially addressed in Theorem 17(i) of \cite{FHK2012}, which is a part of Theorem \ref{T:1}(I.a); on the entire interval the orbital stability (for $\delta$-small perturbations) is in \cite{NTV2014}.  
The intervals of stability and instability near the respective endpoints in Part III are also obtained in the same works \cite{FHK2012}, \cite{FM2001}, \cite{NTV2014}.
\end{remark}

In this work we provide analytical and numerical evidence for all three parts of Conjecture \ref{C:1} in dimensions 1 and 2. In particular, we show numerically that the standing waves are {\it stable} (including sufficiently large perturbations) in the subcritical and critical cases considered, showing the resulting time evolutions. In the supercritical case computations, the stable and unstable regimes meet at the turning point $b^*$. 
In Proposition \ref{P:1D-b*} we prove existence of a mass maximum and characterize its location in 1D. Furthermore, at $b^\ast$ we observe that branching of ground states occurs, which not only determines different branches of stability, but also different types of stable behavior, depending on the branch as discussed in Remark \ref{R:1} (see further details in Section \ref{S:branching}). We also prove blow-up conditions for Part II\,(ii) and part of the blow-up statement in Part III\,(iii) of the conjecture in Proposition \ref{P:blowup}.  
\smallskip

This paper is organized as follows: in Section \ref{S:GS} we review useful identities and convergence properties of the ground states. In Section \ref{S:numerics} we describe our numerical methods and their accuracy. In Section \ref{S:results} we present our analytical and numerical results, first discussing the convergence in two limiting regimes, mass-energy curves, and the 1D turning-point characterization in Proposition \ref{P:1D-b*}. We then examine branching and perturbation dynamics, prove the blow-up criteria in Proposition \ref{P:blowup}, and conclude with small-amplitude dynamics, thus, providing evidence for Conjecture \ref{C:1}. Section \ref{S:conclusions} summarizes the results. 
\smallskip

{\bf Acknowledgments.} The work of C. K. and N.S. was partially supported by the ANR project 
ANR-17-EURE-0002 EIPHI and by the ANR project 
ISAAC-ANR-23-CE40-0015-01. \\
S.R. would like to thank IMB for its hospitality and support, where most of this work was done. S.R. was partially supported by the NSF grant DMS-2452782.

\section{Properties of ground states}\label{S:GS}
In this section we discuss several useful properties of ground states $Q_b$ from \eqref{E:GS-NLS}.

\subsection{Pokhozhaev identities on the unit ball}\label{S:Pokh}
The following Pokhozhaev identities are useful for numerical verification of the ground state $Q_b$. 
The first identity is obtained by multiplying \eqref{E:GS-NLS} by $Q_b$ and integrating over $B_1$, 
\begin{equation}\label{E:P1}
\|\nabla Q_b\|^2_{L^2(B_1)} + b\|Q_b\|^2_{L^2(B_1)} - \|Q_b\|^{\alpha+2}_{L^{\alpha+2}(B_1)} = 0.
\end{equation}
The second identity is obtained by multiplying \eqref{E:GS-NLS} by $x \cdot \nabla Q_b$  
and integrating over $B_1$ (or in spherical coordinates multiplying by  $\sigma(B_1) r^d \partial_r Q_b$, where $\sigma(B_1)$ is the surface area of the unit ball in $\mathbb R^{d}$, and then integrating in $r$), 
\begin{equation}\label{E:P2}
-\frac{1}{2} \sigma(B_1) |\partial_r Q_b(1)|^2 + \frac{2-d}{2} \|\nabla Q\|^2_{L^2(B_1)}  - \frac{b\,d}{2} \|Q\|^2_{L^2(B_1)} + \frac{d}{\alpha+2} \|Q\|^{\alpha+2}_{L^{\alpha+2}(B_1)} = 0.
\end{equation}

Solving for  $\|\nabla Q\|_{L^2(B_1)}^2$ and $\| Q\|^{\alpha+2}_{L^{\alpha+2}(B_1)}$ from \eqref{E:P1} and \eqref{E:P2}, 
we obtain
\begin{align}
\|\nabla Q_b\|^2_{L^2(B_1)} 
&= \frac{1}{\alpha(2-d)+4} \big(d\alpha \, b \, M(Q_b) + (\alpha+2)\sigma(B_1)|\partial_r Q_b(1)|^2 \big), \label{E:gradient} \\ 
\|Q_b\|_{L^{\alpha+2}(B_1)}^{\alpha+2} & = 
\frac{(\alpha+2)}{\alpha(2-d) +4} \big( 2 \, b \, M(Q_b) + 
\sigma(B_1)|\partial_r Q_b(1)|^2 \big). \label{E:potential}
\end{align}

Recalling the energy from \eqref{EC}, we obtain the following relation between the mass, energy and the boundary value term:
\begin{equation}\label{E:EMP}
E[Q_b] =  \frac{d\alpha-4}{2(\alpha(2-d) +4)}\,b\, M(Q_b) + \frac{\alpha}{2(\alpha(2-d) +4)}\sigma(B_1)|\partial_r Q_b(1)|^2.
\end{equation}
Note that in the critical case, $\alpha = \frac4{d}$, the energy is 
\begin{equation}\label{E:EMP2}
E[Q_b] = \frac{\sigma(B_1)}{4} \,|\partial_r Q_b(1)|^2 >0,
\end{equation}
which is different from the NLS on the whole space, where the ground state in the critical case has zero energy. 

\subsection{Changes of mass and energy with respect to $b$}
We next record the dependence on $b$  of energy and mass (of the ground state). Recalling the energy from \eqref{EC} and the ground state equation, we have
$$
\partial_b E(Q_b) = - \int_{B_1} (\Delta Q_b +Q_b^{\alpha+1}) \partial_b (Q_b) = -b \int_{B_1} Q_b \partial_b(Q_b) = -\frac{b}2 {\partial_b} \, M(Q_b).
$$ 
As we track the energy and the mass of $Q_b$ in our simulations, we also verify the above dependence, namely, 
\begin{equation}\label{E:EM-check} 
{\partial_b} E(Q_b) + \frac{b}2 \,  {\partial_b} \,  M(Q_b) = 0, 
\end{equation}
e.g., see right column in Figure \ref{F:ME-1D-p3-5}.

\subsection{Two regimes of convergence of ground states}

As $b$ varies in the ground states $Q_b$, there are two regimes that are typically identified, as $b \to -\lambda_1$ and as $b \to \infty$ (called small and large amplitude regimes in \cite{FM2001}). 
More precisely, the following convergences hold (see Propositions 11, 12 in \cite{FHK2012}, additionally for 2D cubic NLS see \cite{FM2001}). For these ground state convergence statements, we assume the energy-subcritical regime: for $d \geq 3$ $0<\alpha<4/(d-2)$ and any $\alpha>0$ for $d=1,2$.

\begin{lemma}\label{L:conv} 
Let $Q_b$ be the least-energy (ground state) solution of \eqref{E:GS-NLS} for $b \in (-\lambda_1, \infty)$, where $\lambda_1$ is the first (lowest) eigenvalue of the Laplacian $-\Delta$ with Dirichlet boundary conditions in the unit ball $B_1$. For $d \geq 3$ assume $0< \alpha < 4/(d-2)$ and for $d = 1,2$, any $\alpha>0$.

I. {\rm (Small amplitude regime.)} The normalized ground states converge strongly in $H^1_0(B_1)$ to the first (normalized) eigenfunction of the Laplacian operator $(-\Delta)$ on the ball $B_1$ as $b$ decreases down to $-\lambda_1$: 
$$
\frac{Q_b}{\|Q_b\|_{L^{2}}} \to \chi_1  \quad \mbox{as} \quad b \to -\lambda_1,
$$
provided  $\|\chi_1\|_{L^2} = 1$. 

II. {\rm (Large amplitude regime.)} The rescaled ground states $\tilde Q_b(x) = b^{-1/\alpha} Q_b (\frac{x}{\sqrt{b}})$ for $x \in B_{\sqrt b}$ and zero otherwise, converge strongly in $H^1(\mathbb R^d)$ to $\mathcal R (x)$  as $b \to \infty$, where $\mathcal R$ is the ground state for the NLS equation on the whole space, i.e., 
\begin{equation}\label{E:R}
-\Delta \mathcal  R + \mathcal R - |\mathcal R|^{\alpha} \mathcal R = 0, \quad  ~~ x \in \mathbb R^d. 
\end{equation}
\end{lemma}
We illustrate this behavior in one and two dimensions in Sections \ref{S:convergence} and \ref{S:convergence2}. We note that from rescaling \eqref{E:Q-rescale} and Lemma \ref{L:conv}, we have 
$$
M(Q_b) = b^{2/\alpha - d/2} \big(M(\mathcal R) +o(1) \big), \quad \mbox{as} ~~b \to \infty.
$$
Thus, the (unscaled) mass of $Q_b$ tends to infinity in subcritical case, is increasing to the value $M(\mathcal R)$ as $b \to \infty$ in the critical case (e.g., see  $M = M(b)$ curve dependence in Figures \ref{F:ME-1D-p3-5} and \ref{F:ME-2D} plots (A) \& (D) ), and  
goes to zero in the supercritical case (which is essential later for the existence of a turning point in \S \ref{S:stab}).

\section{Numerical approach}\label{S:numerics}
In this section we describe the numerical methods used to construct ground states and to evolve perturbations. Since the paper is a numerical study, we include the accuracy diagnostics used in the computations. The spatial discretization is based on Chebyshev spectral collocation, while the time evolution is computed with a Crank--Nicolson method combined with a fixed-point/Newton iteration.

\subsection{Numerical construction of ground states in 1D}
In the 1D case, we work with $x\in [-1,1]$ and apply standard 
Chebyshev collocation  points $x_{n}=\cos(n\pi/N)$, $n=0,1,\ldots,N$, 
$N\in\mathbb{N}$, see \cite{trefethen} for a 
comprehensive review of such methods. Using standard Lagrange 
interpolants on these collocation points and approximating the 
derivative of a function by the derivative of the interpolant, we get an approximation of 
derivatives via well known Chebyshev differentiation matrices $D$, see 
\cite{trefethen, weireddy}. 

With this discretization, equation (\ref{E:GS-NLS}) can be 
approximated by a system of algebraic equations
\begin{equation}
	(D^{2}-b) \mathbf{Q}+\mathbf{Q}^{\alpha+1}=0,
	\label{GSdisc}
\end{equation}
where $\mathbf{Q}=(Q_{0},\ldots,Q_{N})$ with $Q_{n}:=Q(x_{n})$, 
$n=0,\ldots,N$. The vanishing conditions for $x=\pm1$ are implemented 
as in \cite{trefethen} by omitting the columns and lines in 
\eqref{GSdisc} with the indices 0, $N$. Denoting 
$\tilde{\mathbf{Q}}:=(Q_{1},\ldots,Q_{N-1})$ and 
$\tilde{D}^{2}_{nm}=D^{2}_{nm}$, $n,m=1,\ldots,N-1$, we get for 
\eqref{GSdisc} the $N-1$ dimensional system
\begin{equation}
	F:= (\tilde{D}^{2}-b) \mathbf{\tilde{Q}}+\mathbf{\tilde{Q}}^{\alpha+1}=0
	\label{GSdisct},
\end{equation}
which is solved by a standard Newton iteration. Typically, we take for 
values of $|b|\sim 1$ initial iterates of the form 
$Q= \lambda \sin (\pi(x+1)/2)$, where $\lambda>0$ is chosen in 
dependence of $\alpha$ and $b$ in order 
to get rapid convergence. The iteration is stopped when the residual 
of $F$ is smaller than some prescribed threshold, generally 
$10^{-10}$. The resulting solution $Q$ is then used in 
a tracing approach as an initial iterate for a nearby $b$. This 
allows us to reach values of $b$ on the entire considered interval $(-\lambda_1, \infty)$, i.e., very large values of $b$ 
and values of $b$ approaching $-\lambda_{1}$. 

The choice of the number $N$ of collocation points is controlled by 
the \emph{Chebyshev coefficients}: an analytic function can be expanded in 
terms of  Chebyshev polynomials $T_{n}(x)= \cos(n\arccos(x))$, 
$n=0,1,\ldots$, and be approximated by a sum
\begin{equation}
	Q(x)\approx\sum_{n=0}^{N}a_{n}T_{n}(x)
	\label{cheb}.
\end{equation}
The \emph{Chebyshev coefficients} can be determined by a collocation 
method, i.e., the condition
\begin{equation}
	Q(x_{k})= \sum_{n=0}^{N}a_{n}T_{n}(x_{k}),\quad k=0,\ldots,N
	\label{chebcol}.
\end{equation}
The coefficients $a_{n}$ can be computed with a \emph{Fast Cosine 
transform}, see \cite{trefethen}. It is known that they decrease 
exponentially with the index for analytic functions. The numerical 
error is controlled by the coefficients with the largest values of 
$n$. Since the ground states $Q$ are expected to be analytic, we 
always choose $N$ large enough that $|a_{N}|$ is on the order of 
machine precision, here, roughly $10^{-16}$. The numerical error in 
approximating analytic functions thus decreases exponentially with $N$, which is called \emph{spectral convergence}. 

An advantage of the use of Chebyshev collocation points is that 
the integrals can be approximated on these points by the Clenshaw-Curtis 
algorithm \cite{clencurt}:
\begin{equation}
	\int_{-1}^{1}f(x)dx\approx \sum_{n=0}^{N}f(x_{n})w_{n}
	\label{CC},
\end{equation}
where the $w_{n}$, $n=0,\ldots,N$, are known weights, see 
\cite{trefethen}. 
The Clenshaw-Curtis algorithm is once more a spectral method. 

\subsection{Numerical construction of ground states in 2D}
In the 2D case, a similar approach is applied for the radial 
coordinate. We write equation \eqref{E:GS-NLS} in the form
\begin{equation}
	Q_{rr}+\frac{1}{r}Q_{r}-bQ+Q^{\alpha+1}=0,
	\label{GSrr}
\end{equation}
and introduce as in \cite{CKS} the coordinate $s=r^{2}$ to obtain a 
less singular form of the equation (\ref{GSrr}),
\begin{equation}
	4sQ_{ss} + 4 Q_{s}-bQ + Q^{\alpha+1}=0.
	\label{GSs}
\end{equation}
As in the 1D case we apply a Chebyshev collocation method, i.e., we 
discretize $s\in[0,1]$ via $s_{n}=(1+\cos(n\pi/N))/2$, $n=0,\ldots,N$. 
The rest of the approach is as in the 1D case. 

We verify our computation of the 2D ground states via the virial and energy identities, as discussed in Section \ref{S:Pokh} in \eqref{E:gradient}-\eqref{E:EMP}. 
In particular, we compute the following errors for the ground states in the cubic $\alpha=2$ and quintic $\alpha=4$ cases: 
$$
\mathcal{E}_1(b): = \|\partial_r Q_b\|^2_{L^2(B_1)} - \frac{\alpha b}{2}  M(Q_b) - \frac{\pi(\alpha+2)}{2} |\partial_r Q_b(1)|^2,
$$
$$
\mathcal{E}_2(b) := E(Q_b)- \frac{(\alpha-2) b}{4}  M(Q_b) - \frac{\pi \,\alpha}{4}|\partial_r Q_b(1)|^2.
$$
For example, for $b=1$ when the ground state $Q_b$ is computed up to 
the resolution of $10^{-9}$, we obtain that both errors $\mathcal 
E_1$ and $\mathcal E_2$ produce the values of the order $10^{-11}$ 
for the $\alpha=2$ case and $10^{-12}$ for the $\alpha=4$ case. For 
other values of $b$ we find similar precision.

\subsection{Time integration}
Applying the spatial discretization above in $x$, we obtain for the 
NLS equation \eqref{NLS} in one dimension an $N-1$ dimensional system of the form 
\begin{equation}
\tilde{u}_{t} = i\tilde{D}^{2}\tilde{u}+i|\tilde{u}|^{\alpha}\tilde{u}.
\label{ut}
\end{equation}
In the radial two-dimensional case, we replace $\tilde D^2$ by the discretization for $4s \partial_s^2+4\partial_s$, with the origin treated by regularity and the outer boundary has the Dirichlet condition. 

We discretise also in time, $t=0,\ldots, t_{n},\ldots,T$, where $T$ is 
the final time, and where $N_{t}>0$ time steps of length $h = T/N_t$ are used. 
The system \eqref{ut} is integrated by the 
Crank-Nicolson method: integrating the  equation \eqref{ut} from $t_{n}$ 
to $t_{n+1}$ and applying on the right-hand side the trapezoidal rule, 
leads to 
\begin{equation}
(1-\frac{h}2 \, i\tilde{D}^{2})	\tilde{u}(t_{n+1})=(1+\frac{h}2\, i\tilde{D}^{2})	\tilde{u}(t_{n})
+i\frac{h}{2}\,|\tilde{u}(t_{n+1})|^{\alpha}\tilde{u}(t_{n+1})
+i\frac{h}{2}\,|\tilde{u}(t_{n})|^{\alpha}\tilde{u}(t_{n})
	\label{CN}.
\end{equation}
The nonlinear system \eqref{CN} is solved for $\tilde{u}(t_{n+1})$ by a 
simplified Newton/fixed-point iteration: in each inner iteration we invert the linear 
operator $(1-\frac{h}2 \,i\tilde{D}^{2})$ (so not the full 
Jacobian is inverted), while the nonlinear term is updated from the previous inner iterate. 
The iteration is stopped when the difference between two consecutive inner iterates is below the prescribed tolerance.

The accuracy of the time integration is monitored through the relative
variation of the numerically computed energy. This diagnostic measures the
change of the discrete energy between $t=0$ and later times, normalized by
the size of the initial discrete energy. It is not an a priori error estimate, but
a standard practical check for Hamiltonian dispersive equations. As discussed
in \cite{etna}, relative conservation of the energy typically overestimates the
actual accuracy by two to three orders of magnitude. Except in blow-up
simulations, where the computation is stopped once the relative energy
variation exceeds $10^{-3}$, we choose the time step so that the monitored
quantity remains below $10^{-6}$.


\section{Analytical and Numerical results}\label{S:results}

We start with illustrating the two convergence regimes of Lemma \ref{L:conv}, namely $b\to\infty$ and $b\to-\lambda_1$. We then study the slope of the mass curve, including its integral characterization in 1D in Proposition \ref{P:1D-b*}, and time evolutions of perturbed ground states, which provide numerical evidence for Conjecture \ref{C:1} in dimensions 1 and 2. We conclude with the blow-up criteria in Proposition \ref{P:blowup} and small-amplitude dynamics. 

\subsection{Convergence as $b \to \infty$}\label{S:convergence}
\smallskip

Recall that $Q_b$ satisfies 
\begin{equation}\label{E:Q_b}
-\Delta Q +b Q -|Q|^\alpha Q = 0, \quad \mbox{for} ~~ x \in B_1, 
\quad \mbox{and} ~~ Q_{|_{\{|x|=1 \} }} = 0. 
\end{equation}

To investigate the large amplitude regime convergence, we define a rescaled $\tilde Q_b$ by 
\begin{equation}\label{E:Q-resc}
Q_{b}(x) \stackrel{def}{=} b^{1/\alpha} \tilde Q_b(\sqrt b \,x).
\end{equation} 
Then $\tilde Q_b$ satisfies
\begin{equation}\label{E:tildeQ_b}
\qquad -\Delta \tilde Q + \tilde Q - |\tilde Q|^\alpha \tilde Q = 0, 
\quad \mbox{for} ~~ y \in B_{\sqrt{b}} \quad \mbox{and} ~~ \tilde{Q}_{ |_{\{|y|=\sqrt b \} } }=0.
\end{equation}
Notice that there is no $b$ in front of the second term in \eqref{E:tildeQ_b}, instead the domain of $\tilde Q_b$ depends on $\sqrt b$ and expands to the whole space as $b$ grows to infinity, thus, giving justification to the convergence of the rescaled ground states $\tilde Q_b$ (which are extended to the whole space by zero) to $\mathcal R$ 
as $b$ grows: $\tilde Q_b(x)  \to \mathcal R (x)$  as $b \to \infty$, where $\mathcal R$ is the ground state on the whole 
space $\mathbb R^d$, solving \eqref{E:R}.  

\begin{figure}[h!]
\centering
\begin{minipage}[b]{0.43\linewidth} 
\centering
\begin{tikzpicture}
\begin{axis}[width=\textwidth, height=5.5cm,     
                xlabel style={font=\tiny},
                ylabel style={font=\tiny},
                tick label style={font=\tiny},
			    xlabel={$x$},                
                ylabel={$\mathcal{R}_\alpha$},
                xmin=-6, xmax=6,
                ymin=0, ymax=1.6,
                grid=both,
                legend style={at={(1, 1)}, anchor=north east, draw=black, fill=white, font=\tiny},]
\node at (current axis.north west)
[anchor=south west, yshift=0.3cm, align=left]
{\large $\mathcal{R}_\alpha$};
  \addplot[domain=-6:6, samples=100, smooth, blue, thin] { 1.5 * (1/cosh(0.5*x))^2 };
  \addlegendentry{$\alpha=1$}
\addplot[domain=-6:6, samples=100, smooth, green!70!black, thin] { 1.41421356 * (1/cosh(x)) };
 \addlegendentry{$\alpha=2$}
\addplot[domain=-6:6, samples=100, smooth, magenta, thin] { 1.316074 * (1/cosh(2*x))^0.5 };
\addlegendentry{$\alpha=4$}
\addplot[domain=-6:6, samples=100, smooth, orange, thin] { 1.196 * (1/cosh(5*x))^0.2 };
\addlegendentry{$\alpha=10$}
\addplot[domain=-6:6, samples=100, smooth, red, thin] { 1.066 * (1/cosh(25*x))^0.04 };
 \addlegendentry{$\alpha=50$}
\end{axis}
\end{tikzpicture}
\end{minipage}\hfill
\begin{minipage}[b]{0.54\linewidth} 
\includegraphics[width=1\textwidth,height=.57\textwidth]{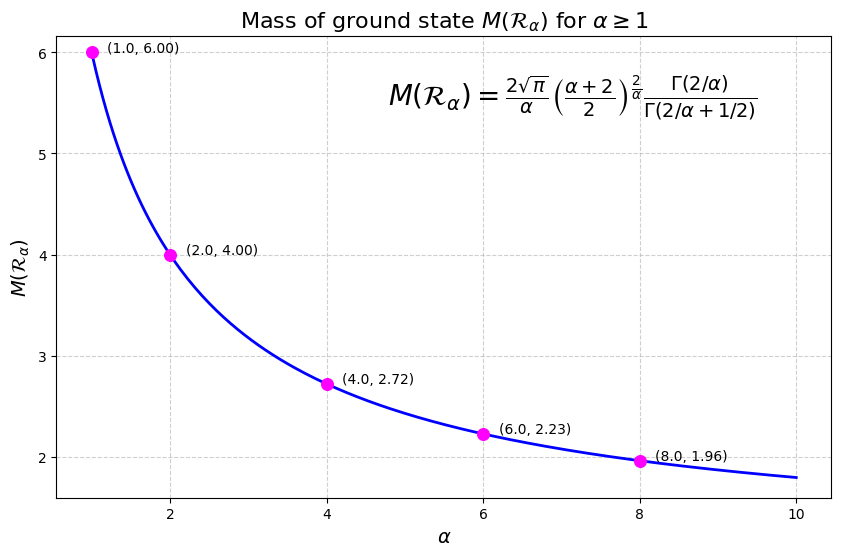}   
\end{minipage}
\caption{\small The profiles of $\mathcal{R}_{\alpha}$ for various values of $\alpha$ (left) and mass   $M(\mathcal{R}_{\alpha})$ (right).}
\label{F:R-MR}
\end{figure}
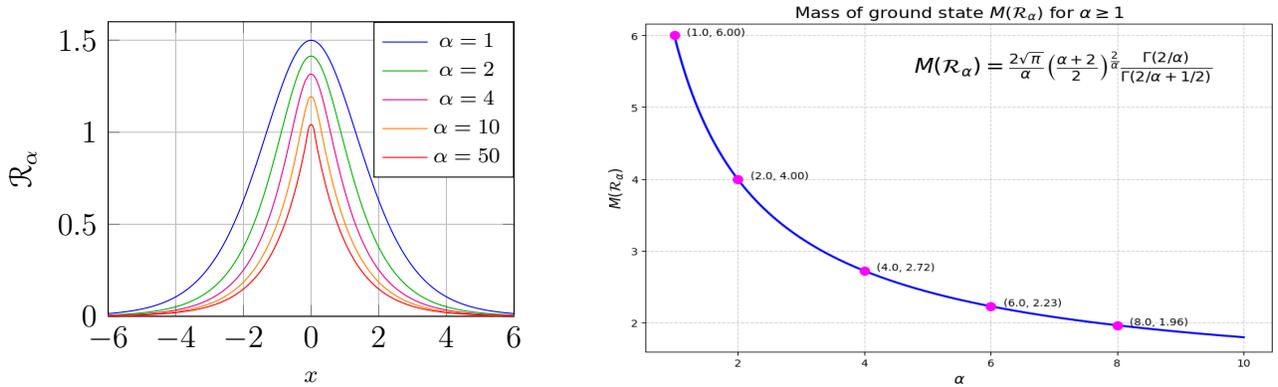

We note that, in particular, in 1D the ground state $\mathcal R = \mathcal{R}_\alpha $ is explicit 
\begin{equation}\label{E:R-alpha}
\mathcal{R}_\alpha (x) = \Big( \frac{\alpha+2}2 \Big)^\frac1{\alpha} \sech^{\frac2{\alpha}} \Big(\frac{\alpha}2 x \Big),
\end{equation} 
see profiles in the left plot of Figure \ref{F:R-MR}. The mass can also be computed explicitly, 
\begin{equation}\label{E:MR}
M(\mathcal{R}_\alpha)=\frac{2}{\alpha} \Big(\frac{\alpha+2}{2} \Big)^{\frac{2}{\alpha}} \int_{-\infty}^{\infty} \sech^{\frac{4}{\alpha}} x\, dx = \frac{2 \sqrt{\pi}}{\alpha} \Big(\frac{\alpha+2}{2}\Big)^{\frac2{\alpha}} \frac{\Gamma(\frac2{\alpha})}{\Gamma(\frac{2}{\alpha} + \frac12)},
\end{equation}
values of which are decreasing, see the right plot in Figure \ref{F:R-MR}. Note that both $M(\mathcal{R}_{\alpha}) \to 1$ and $\|\mathcal{R}_\alpha\|_{L^\infty} \to 1$ as $\alpha \to \infty$.

To avoid confusion between the expanded-domain variable and the unit-ball variable, denote $X=\sqrt b\,x$ and define
\begin{equation}\label{E:Q-rescale}
\widetilde Q_b(X)=b^{-1/\alpha}Q_b(X/\sqrt b),
\qquad X\in B_{\sqrt b}.
\end{equation}
Then $\widetilde Q_b$ solves
$$
-\Delta \widetilde Q_b+\widetilde Q_b-|\widetilde Q_b|^\alpha\widetilde Q_b=0
\quad \text{on } B_{\sqrt b},
$$
and $\widetilde Q_b\to \mathcal R$ as $b\to\infty$. Equivalently, on the
unit ball,
$$
Q_b(x)=b^{1/\alpha}\widetilde Q_b(\sqrt b\,x)
\approx b^{1/\alpha}\mathcal R(\sqrt b\,x).
$$
Thus, for the numerical comparison on $B_1$, we use the whole-space ground state rescaled back to the unit-ball variables, namely the profile
\begin{equation}\label{E:R-rescaled}
{\mathcal R}_b(x)=b^{1/\alpha}\mathcal R(\sqrt b\,x).
\end{equation}
\begin{figure}[htb!]
 \includegraphics[width=0.33\textwidth,height=0.25\textwidth]{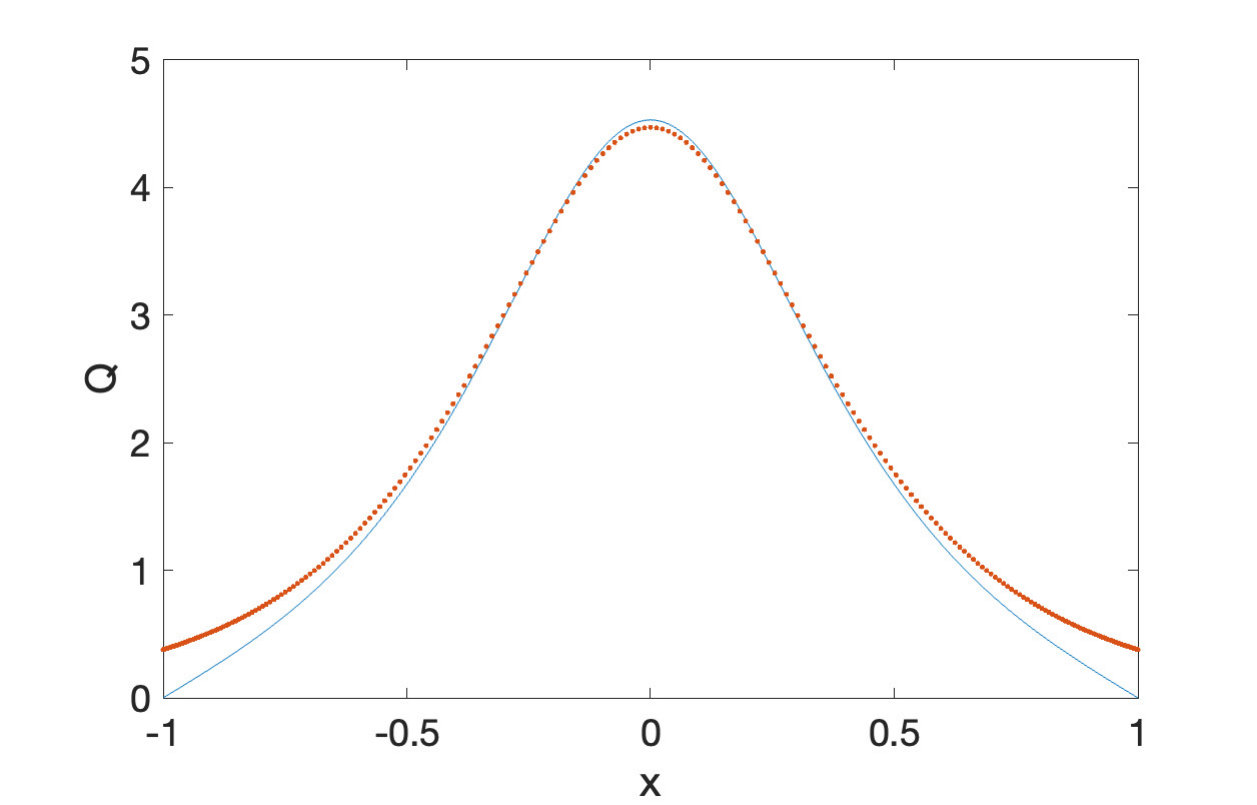}
 \includegraphics[width=0.33\textwidth,height=0.25\textwidth]{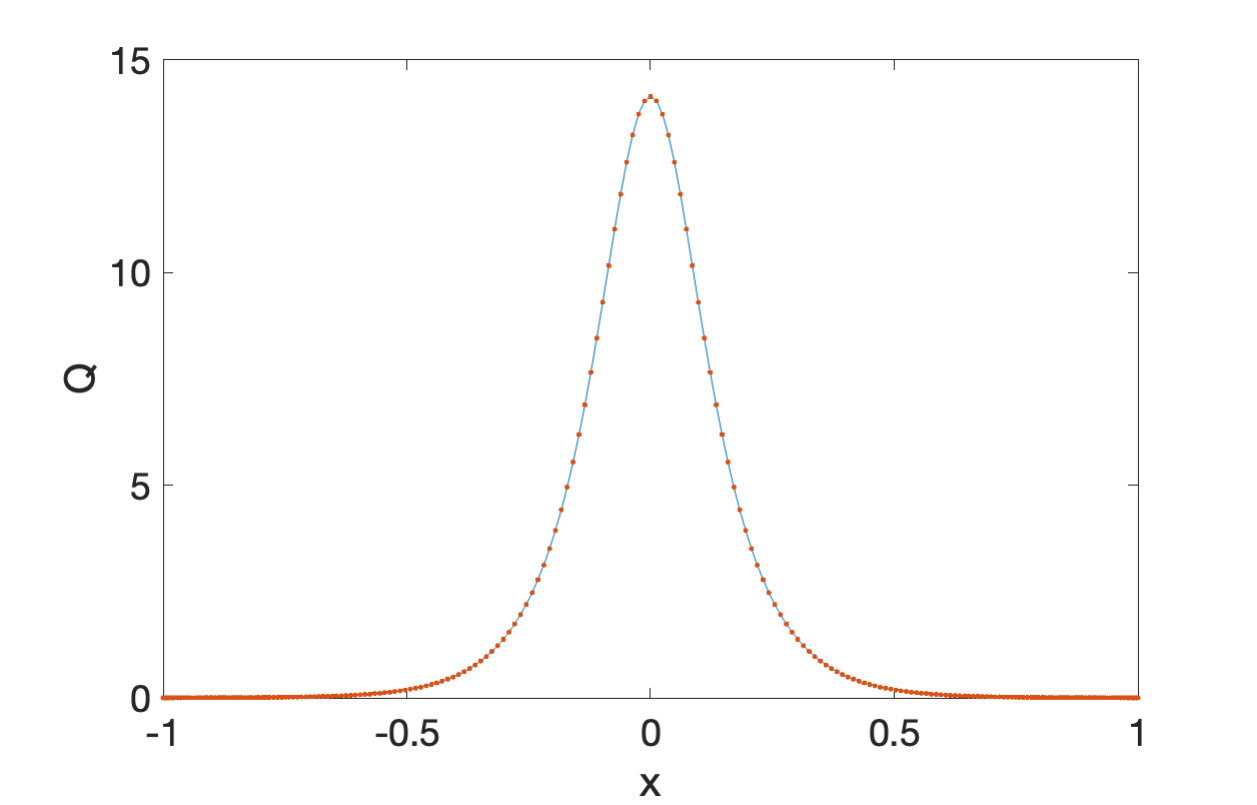}  
 \includegraphics[width=0.32\textwidth,height=0.25\textwidth]{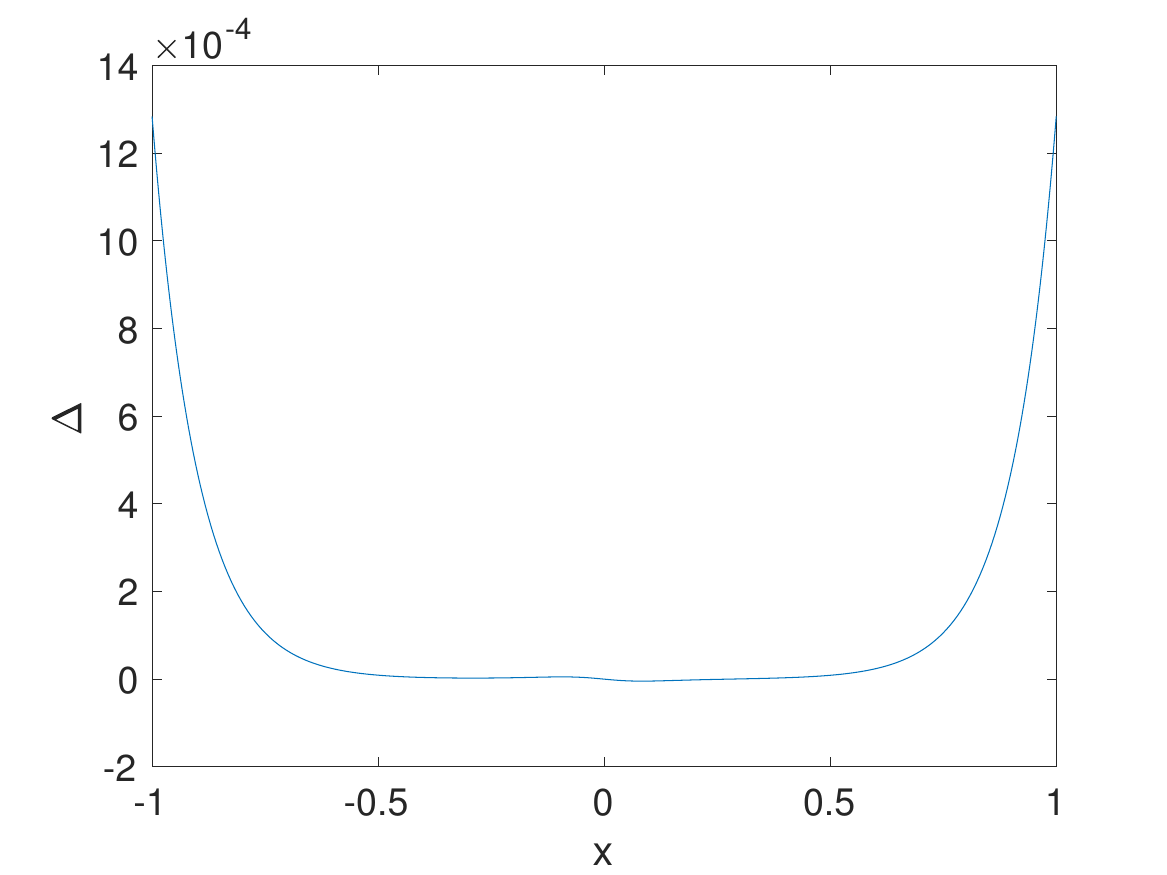}
\caption{\footnotesize 1D cubic NLS. Convergence of the ground states $Q_b$ (solid blue) to the rescaled 1D ground state $\mathcal R = \sqrt 2 \, \sech \,  x$ (dotted red), rescaled to the interval $[-1,1]$ 
as in \eqref{E:R-rescaled}, as $b \to \infty$: $b=10$ (left), $b=100$ (middle), difference for $b=100$ (right).}
 \label{F:1Dcon-p=3}
\end{figure}
In Figure \ref{F:1Dcon-p=3} we show convergence in 1D of $Q_b$ as $b$ becomes large to the rescaled explicit ground state ${\mathcal{R}_b}$ for $b=10$ and $b=100$. While there is some difference for $b=10$, especially on the edges, for $b=100$ the difference is not visible and the error is plotted on the right of Figure \ref{F:1Dcon-p=3}, which is on the order of 10$^{-3}$ at the edges or smaller in the interior.

\begin{figure}[h!]
 \includegraphics[width=0.33\textwidth,height=0.27\textwidth]{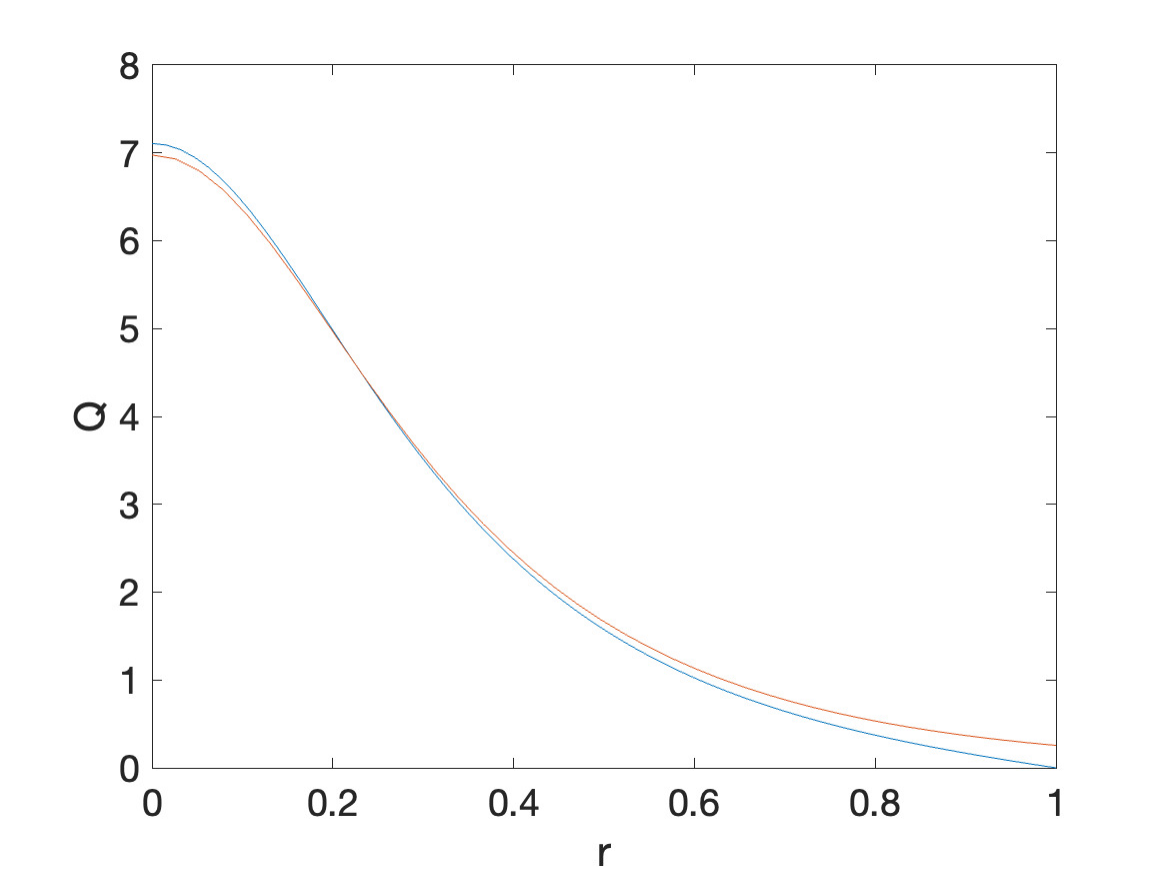}
 \includegraphics[width=0.33\textwidth,height=0.27\textwidth]{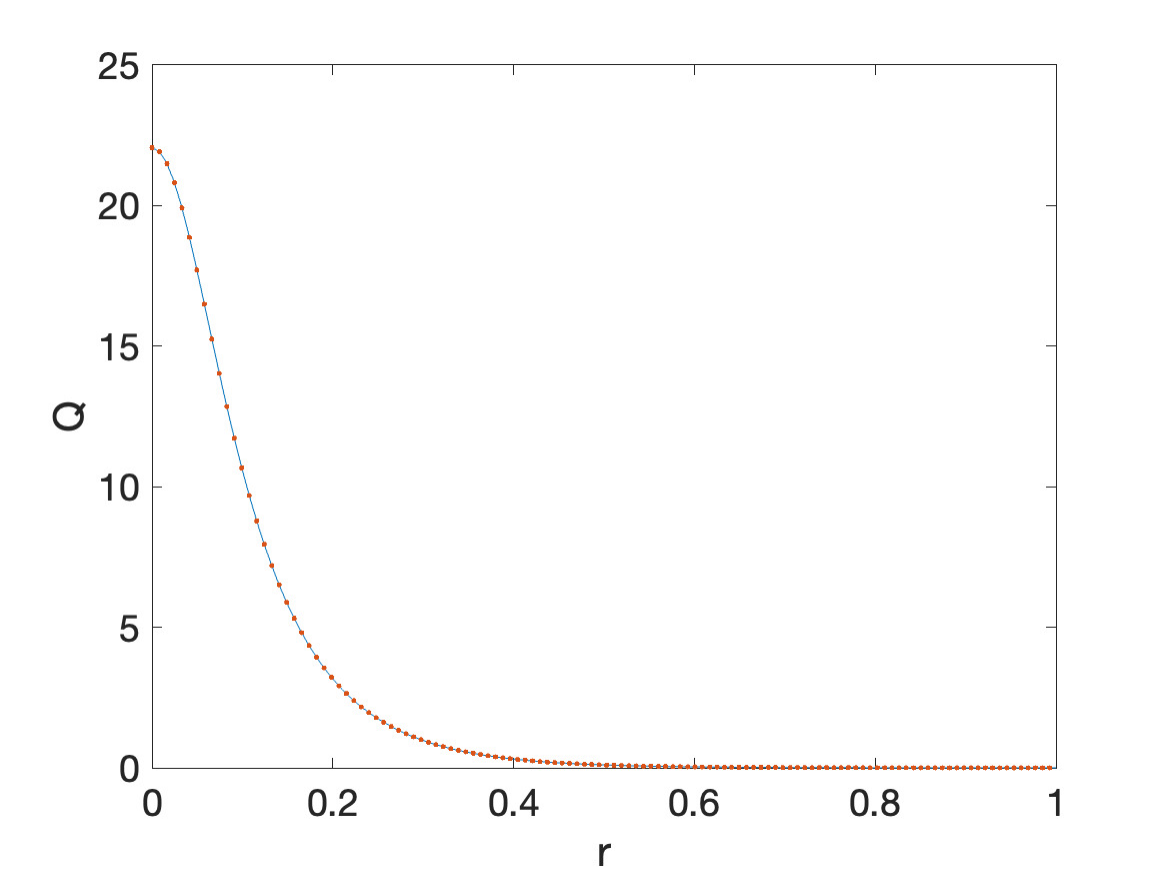}  
 \includegraphics[width=0.32\textwidth,height=0.27\textwidth]{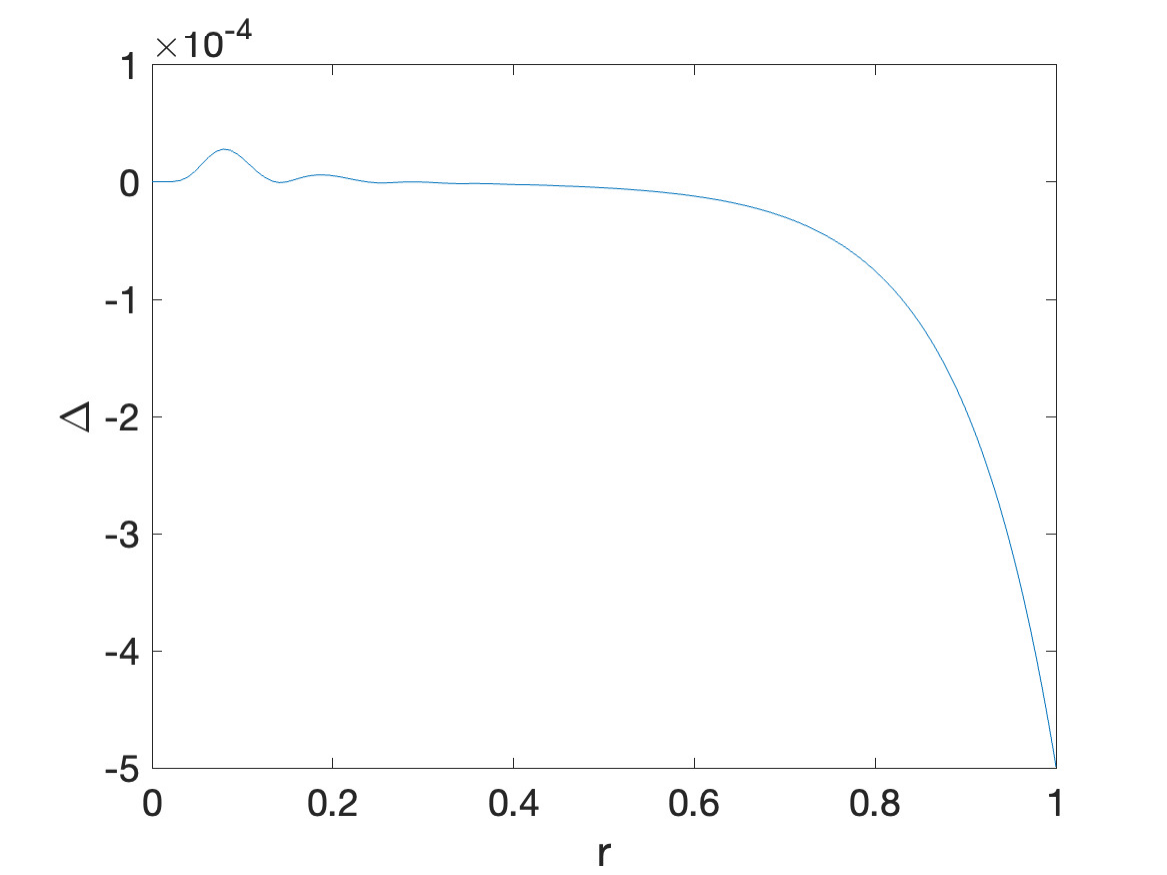}
\caption{\footnotesize 2D cubic NLS. Convergence of the ground states $Q_b$ (solid blue) to the 2D ground state $\mathcal R$ (computed numerically and rescaled to $|r| \leq 1$ as in \eqref{E:R-rescaled}, dotted red) as $b \to \infty$: $b=10$ (left), $b=100$ (middle), difference for $b=100$ (right). }
 \label{F:2Dcon-p=3}
\end{figure}
Similarly, Figure \ref{F:2Dcon-p=3} shows the convergence in 2D of $Q_b$ as $b$ becomes large to the rescaled (numerically computed) ground state $\mathcal{R}_b$ for $b=10$ and $b=100$. On the right plot notice the error (the difference) on the order $10^{-4}$ for $b=100$.

\subsection{Convergence as $b \to -\lambda_1$}\label{S:convergence2}
We next investigate the convergence of ground states in the small amplitude regime, in other words, as $b$ decreases. 

\subsubsection{1D case}
Recall that in 1D the Laplace operator  is $-\Delta = -\frac{d^2}{dx^2}$, and for $k \in \mathbb N$ the eigenvalues and the corresponding eigenfunctions, $-\Delta \chi_k(x) = \lambda_k \chi_k(x)$, $x \in [-1,1]$, are 
$$
\lambda_k = \Big(\frac{\pi k}{2} \Big)^2, \qquad  
\chi_k(x) = \sin \Big(\frac{\pi k}2 (x + 1) \Big).
$$ 
\begin{figure}[htb!]
 \includegraphics[width=0.49\textwidth,height=0.33\textwidth]{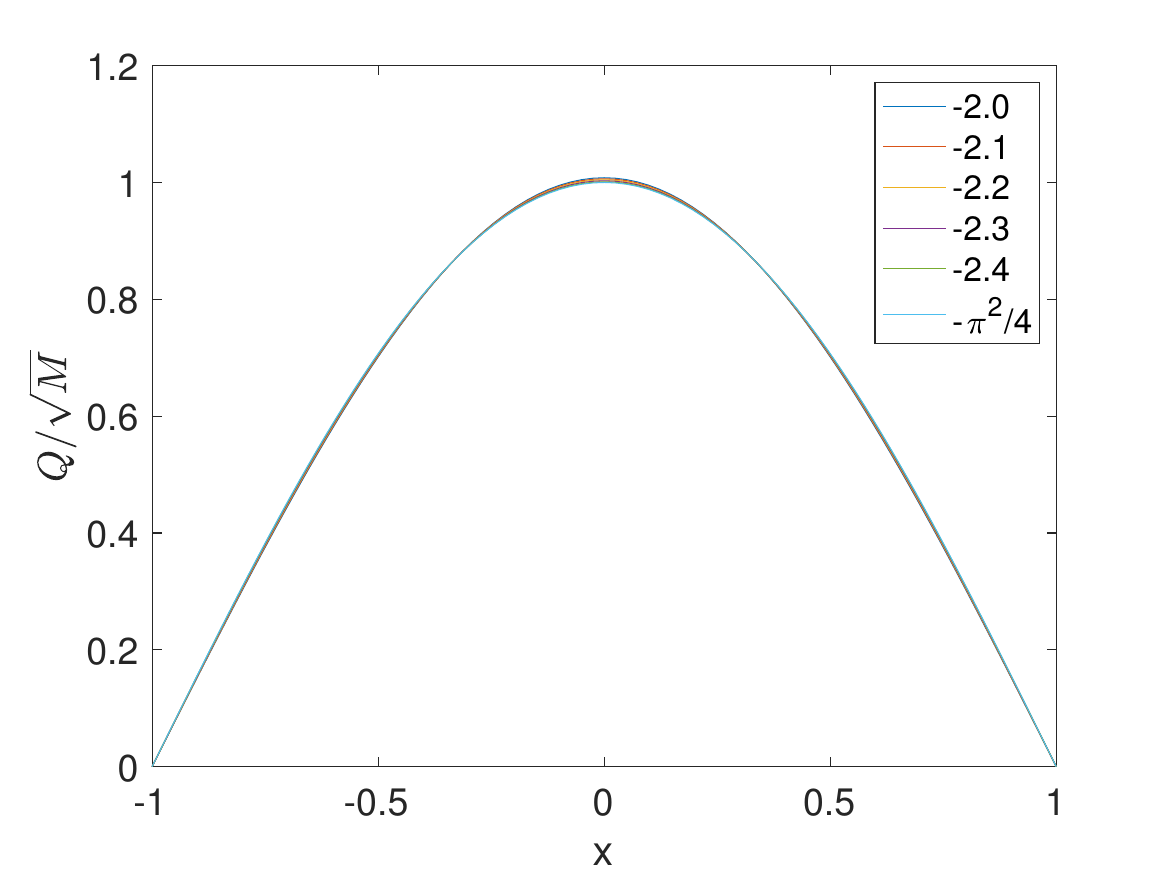}  
\includegraphics[width=0.5\textwidth,height=0.33\textwidth]{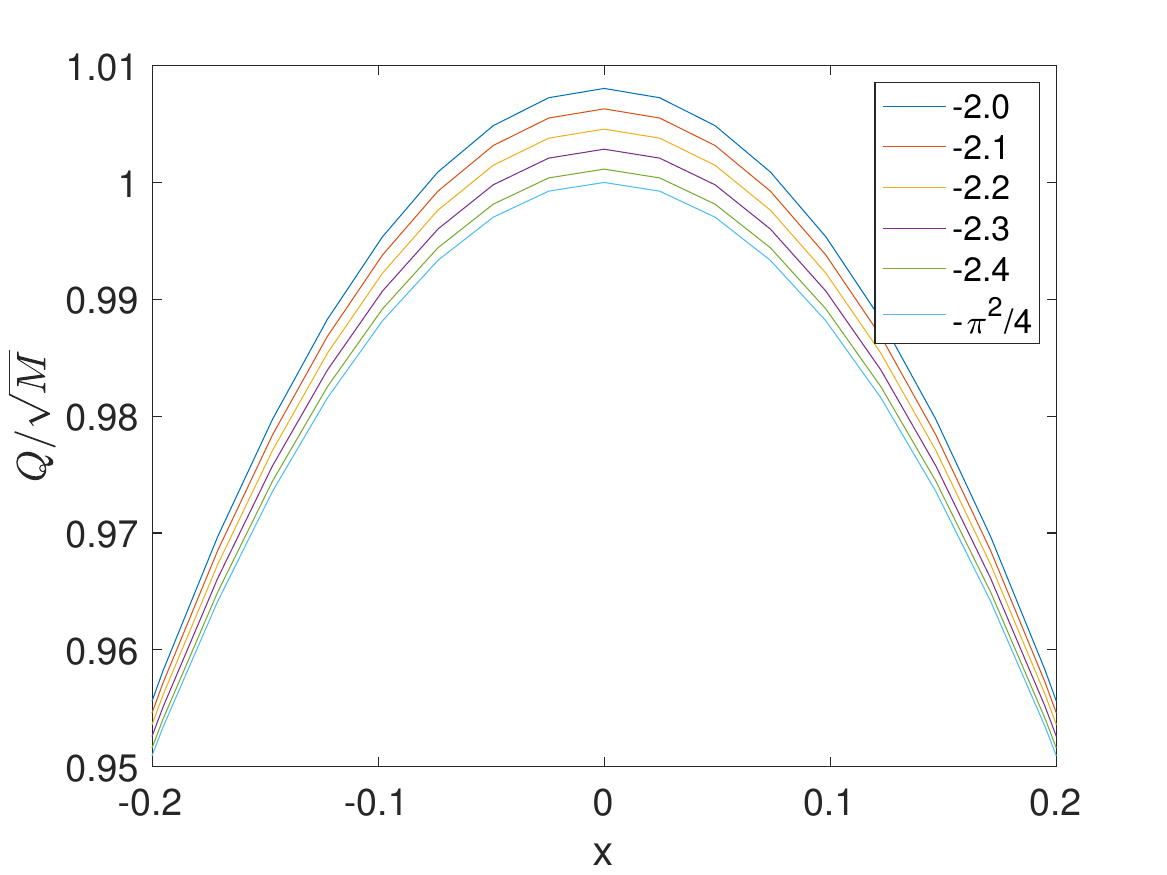}
\caption{\footnotesize 1D cubic NLS (subcritical): convergence of normalized ground states to $\chi_1$ as $b \to -\lambda_1$, which is indicated by the value $-\pi^2/4$ (light blue); zoom-in of the peak on the right. }
 \label{F:1-p=3-conv1}
\end{figure}
We consider subcritical and critical cases, $\alpha = 2$ and $\alpha = 4$, respectively, and show the convergence of normalized $Q_b$ in each case. 
\begin{figure}[htb!]
 \includegraphics[width=0.49\textwidth,height=0.31\textwidth]{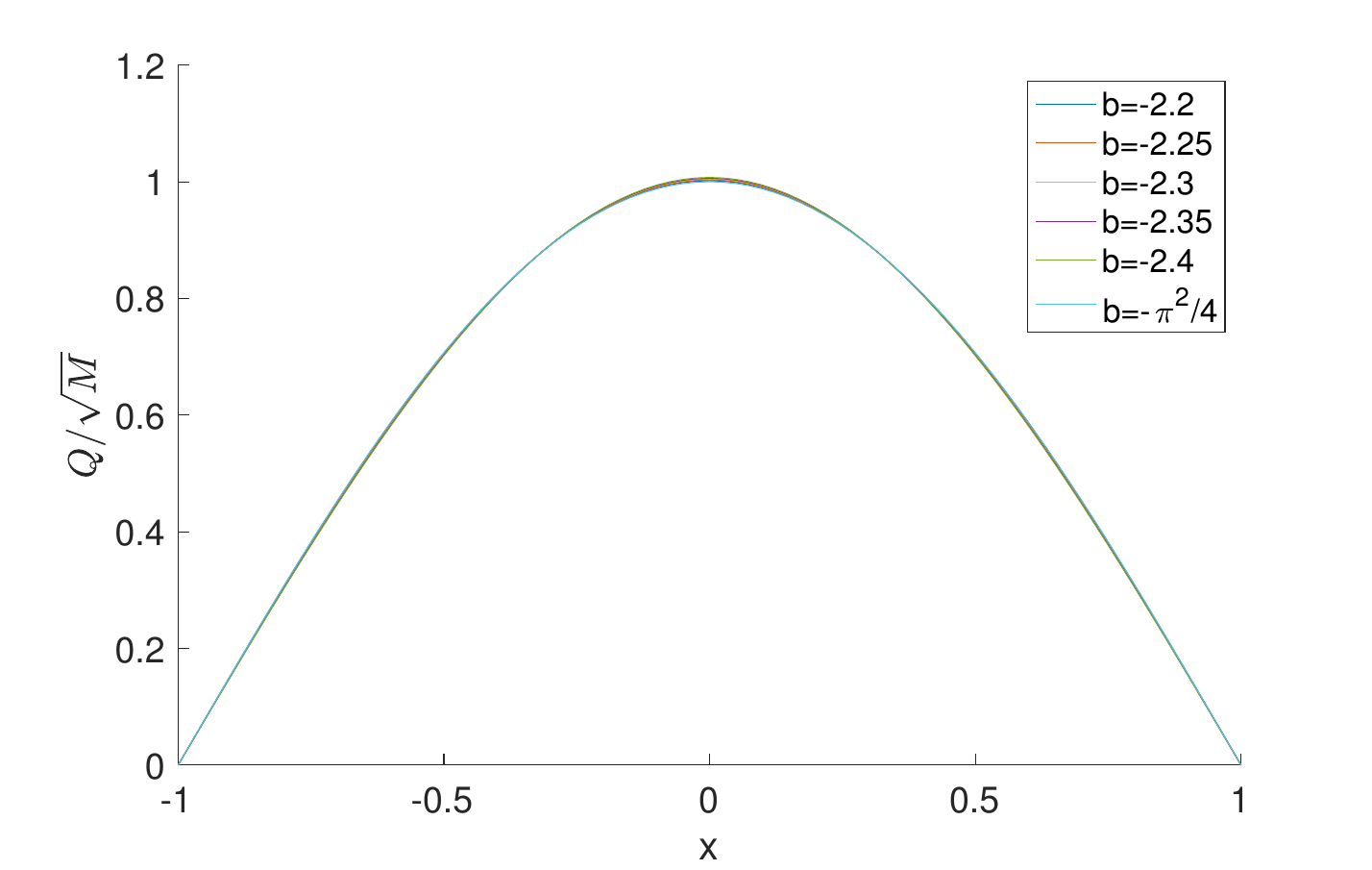}
 \includegraphics[width=0.5\textwidth,height=0.31\textwidth]{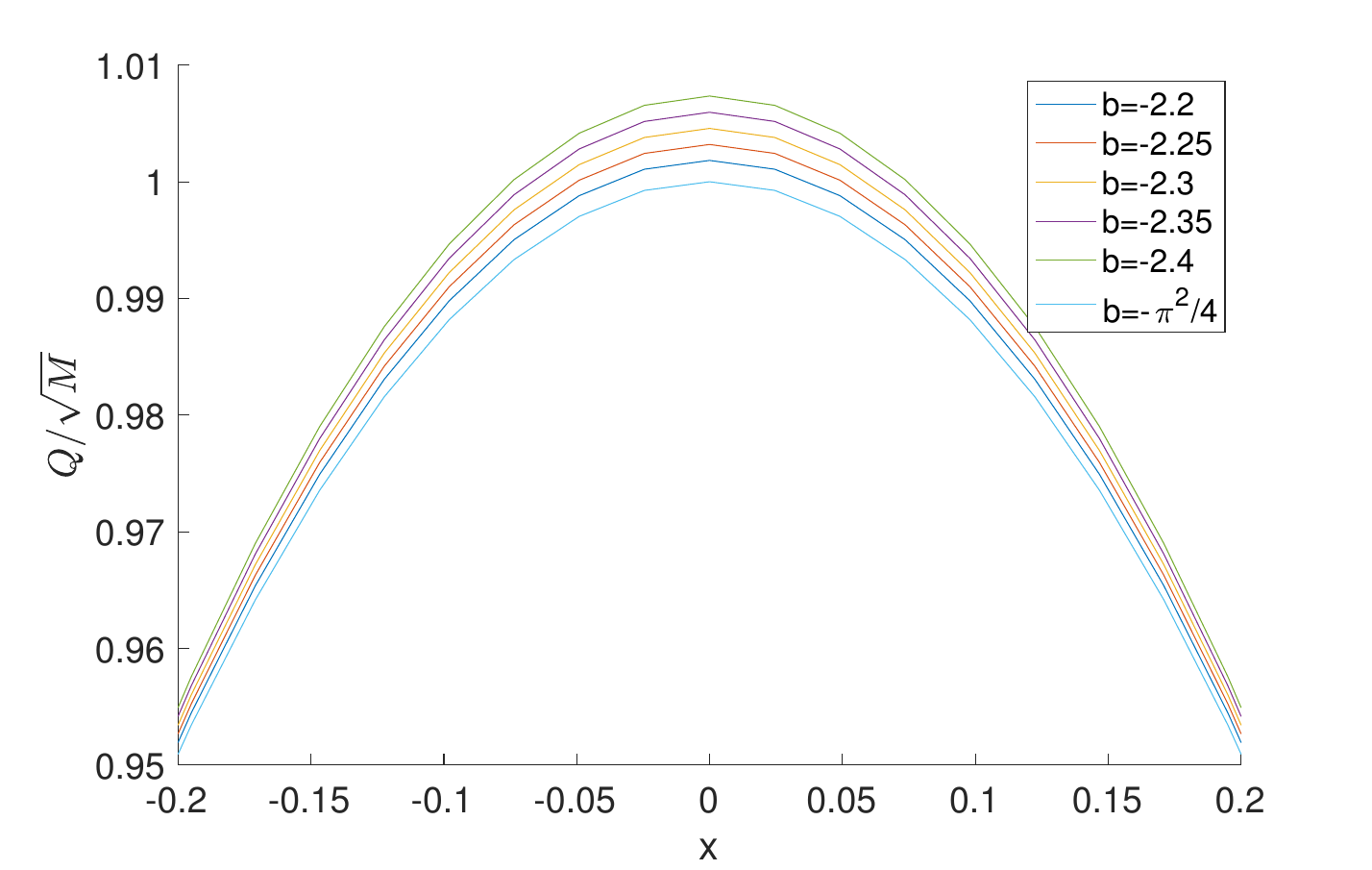} 
  \caption{\footnotesize 1D quintic NLS (critical): convergence of normalized ground states to $\chi_1$ as $b \to -\lambda_1$, which is indicated by the value $-\pi^2/4$ (light blue); zoom-in of the peak on the right. }
 \label{F:1-p=5-conv1}
\end{figure}
In Figures \ref{F:1-p=3-conv1} and \ref{F:1-p=5-conv1} we show that $\frac{Q_{b}\quad }{\|Q_b\|_{L^2([-1,1])}}$ numerically converges to the eigenfunction $\chi_1(x)= \sin \big(\frac{\pi}2 (x+1) \big)$ as $b \to -\lambda_1 = -\big(\frac{\pi}2 \big)^2 $ (note that $\|\chi_1\|_{L^2([-1,1])} = 1$). On the left of Figure \ref{F:1-p=3-conv1} one can see the convergence of the subcritical (cubic NLS) ground states on the whole interval $[-1,1]$, and on the right, the enlargement of the peak, where the convergence is shown for several values. The limiting profile of $\chi_1$ is the lowest curve indicated by the blue line with the eigenvalue $b = -\pi^2/4$.  

A similar convergence holds in the critical, quintic case, see Figure \ref{F:1-p=5-conv1}. 
{We have also checked the supercritical cases, and observed a similar convergence. }

\subsubsection{2D case}
In 2D the Laplacian operator is $-\Delta = -(\partial_{xx} + \partial_{yy}) \equiv -(\partial_{rr} +\frac1{r}\partial_r +\frac1{r^2}\partial_{\theta\theta})$. The eigenvalues and eigenfunctions of the Laplacian operator on a unit ball in 2D, $ - \Delta u = \lambda u$, can be written in polar coordinates as 
$u (r,\theta) = R(r) Y(\theta)$, where the radial function $R$ is expressed via the Bessel functions $J_m$ and their set of zeros $\{k_{m,n}\}$, and $Y(\theta) = \cos(m \theta)$ or $\sin (m \theta)$, namely,
$$
u_{m,n}(r, \theta) = J_m(k_{m,n}r) 
\begin{cases}
\cos(m \theta) \\
\sin(m \theta)
\end{cases}
\qquad \mbox{and} \qquad \lambda_{m,n} = (k_{m,n})^2.
$$
Thus, a full solution can be written as a general eigenfunction expansion, 
$$
u(r, \theta) = \sum_{m=0}^{\infty} \sum_{n=1}^{\infty} J_m (k_{m,n}r) \big( A_{m,n} \cos (m \theta) + B_{m,n} \sin (m \theta) \big).
$$
For this work, we only need the smallest eigenvalue, which is known to be $\lambda_{01} =(k_{01})^2 \approx (2.4048)^2 \approx 5.783$.
By Lemma \ref{L:conv}(I), in 2D we have  $\frac{Q_{b}}{\|Q_b\|_{L^2(B_1)}}$ converge to the normalized first eigenfunction $\frac{\chi_1}{\|\chi_1\|_{L^2}}$, where $\chi_1= J_0(k_{01}r)$, as $b \to -\lambda_1 = -(k_{01})^2 \approx -(2.404825557696086)^2 $ (note that for computational purposes, a very high precision of the constant $k_{01}$ is needed).

\begin{figure}[htb!]
 \includegraphics[width=0.49\textwidth,height=0.35\textwidth]{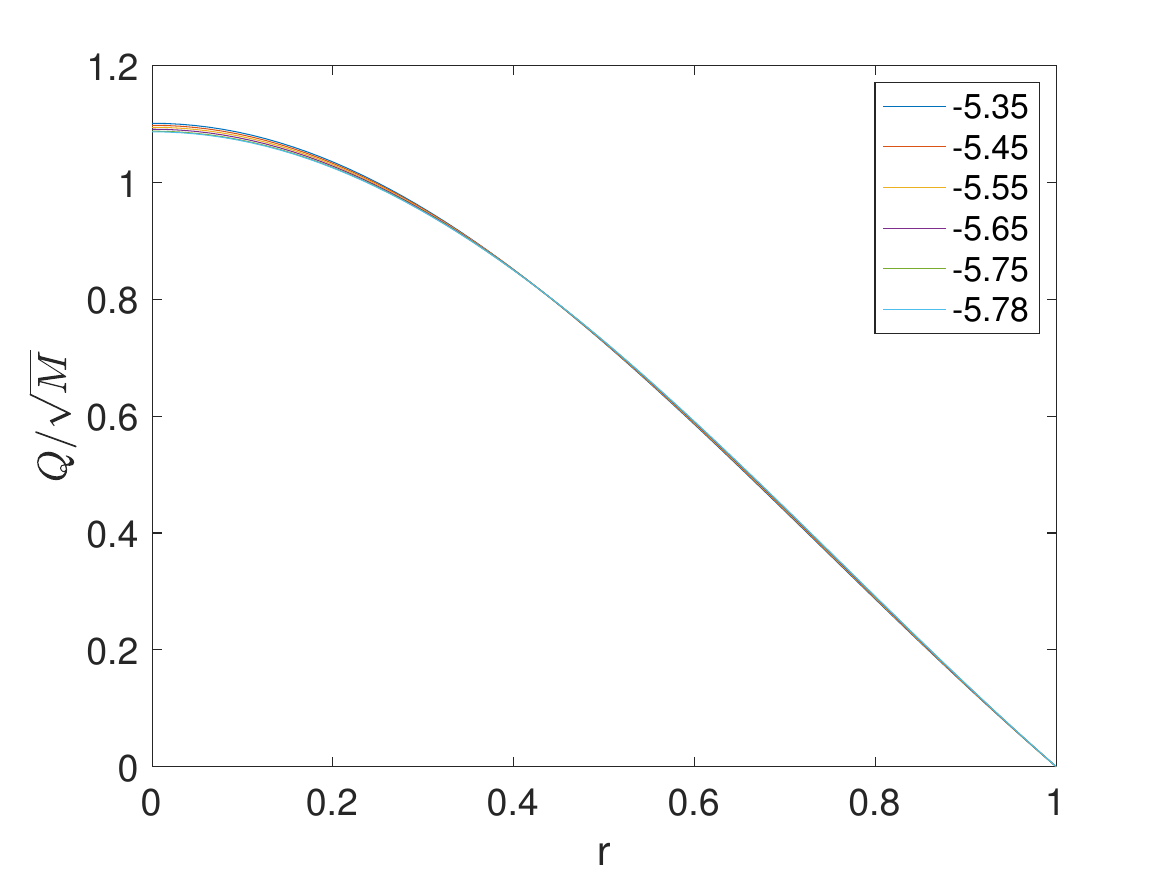}
 \includegraphics[width=0.5\textwidth,height=0.35\textwidth]{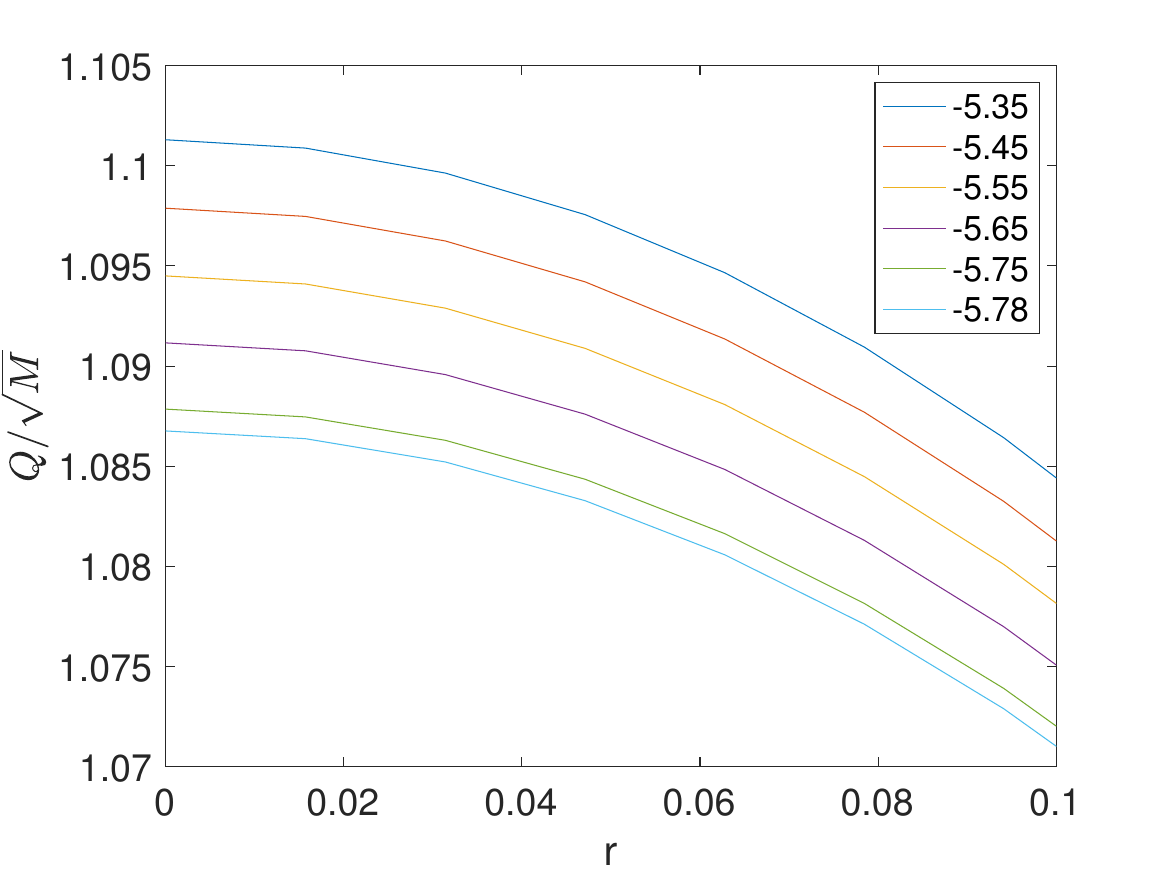}
\caption{ {\small 2D cubic NLS: convergence of the 2D normalized ground states $Q_b$ (with $b$ as labeled) to the normalized $J_0(k_{0,1}r)$,  labeled by $-(k_{0,1})^2 \approx -5.78$ in light blue. } }
 \label{F:2D-cubic}
\end{figure}
In the 2D cubic case ($\alpha=2$) in Figure \ref{F:2D-cubic} we show the convergence of normalized $Q_b$ to normalized $\chi_1$ as $b$ goes to $-(k_{01})^2$. We take a decreasing sequence by $0.1$ starting from $b = -5.35$ and show its convergence to the normalized first eigenfunction, indicated in Figure \ref{F:2D-cubic} in light blue via the value $-5.78$.

\subsection{Stability of solitary waves}\label{S:stab} 
In this section we compare the numerical results with the analytical stability summarized in Theorem \ref{T:1} and investigate stability further, providing numerical evidence towards Conjecture \ref{C:1}. 
\smallskip

{\it Mass and energy as functions of $b$.}
We start by checking how the mass $M(Q_b$) and energy $E(Q_b$) depend on $b$. We discuss the 1D case first, providing the numerical findings and then Prop. \ref{P:1D-b*}, where we show the existence and location of a turning point (in branching) $b^*$.

\begin{figure}[!htb]
\begin{subfigure}{.32\textwidth}
\includegraphics[width=1\linewidth,height=0.75\linewidth]{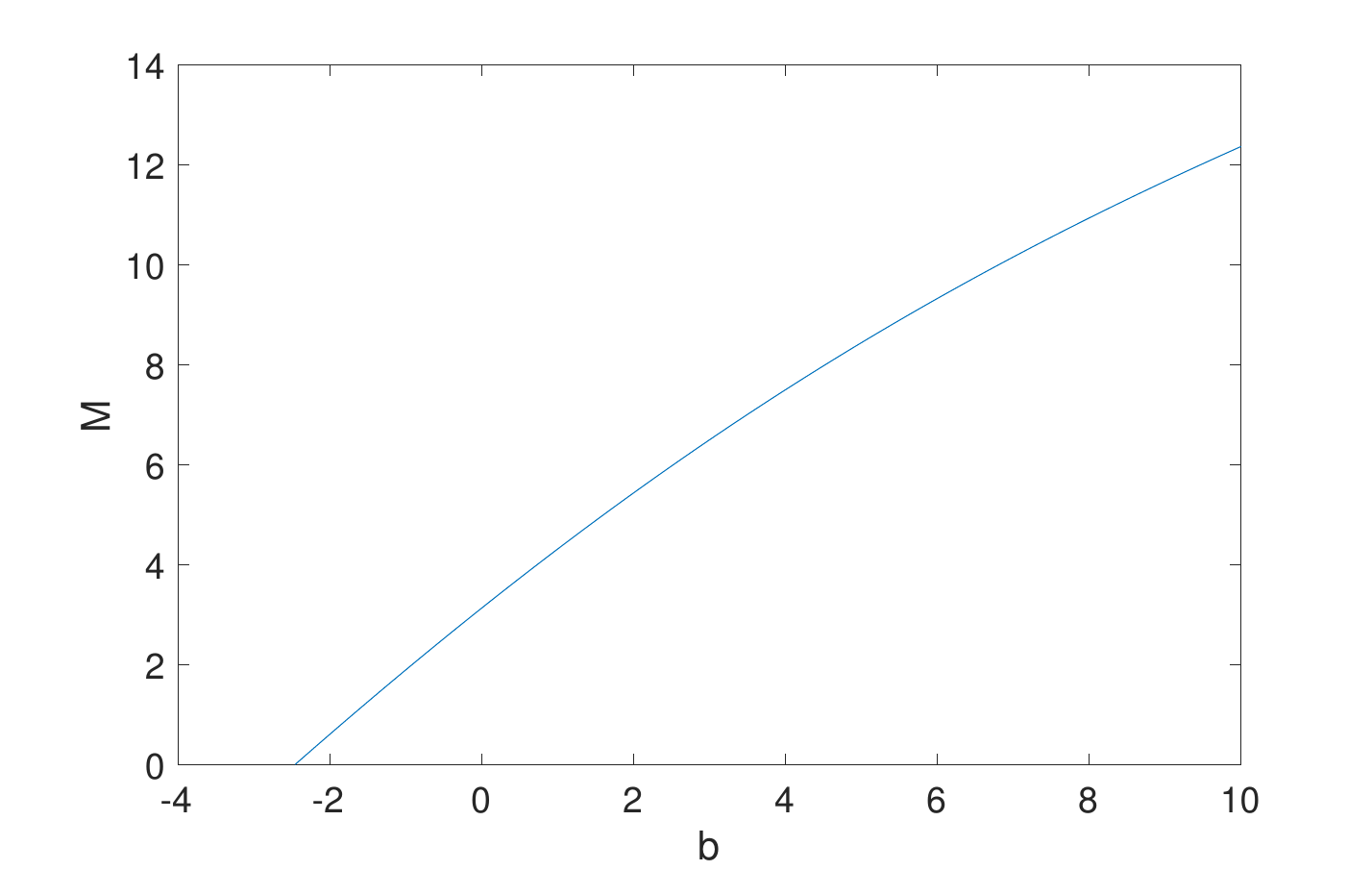}
\subcaption[]{{\footnotesize {$\alpha=2$, $M(Q_b) \to \infty$}}}
\end{subfigure}
\begin{subfigure}{.32\textwidth}
\includegraphics[width=1\linewidth,height=0.75\linewidth]{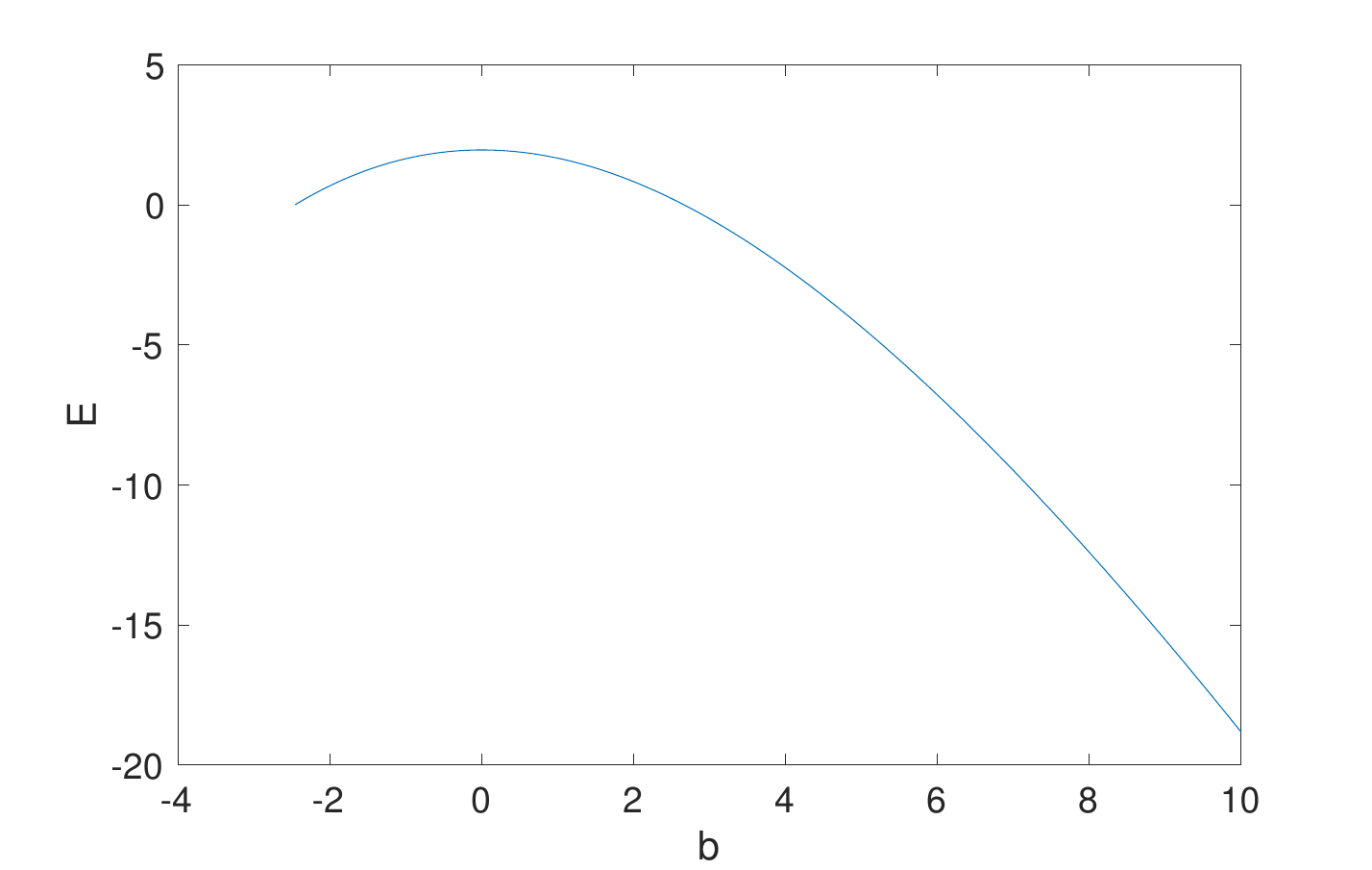}
\subcaption[]{{\footnotesize $E(Q_b)$ as function of $b$}}
\end{subfigure} 
\begin{subfigure}{.32\textwidth}
  \includegraphics[width=1\linewidth,height=0.75\linewidth]{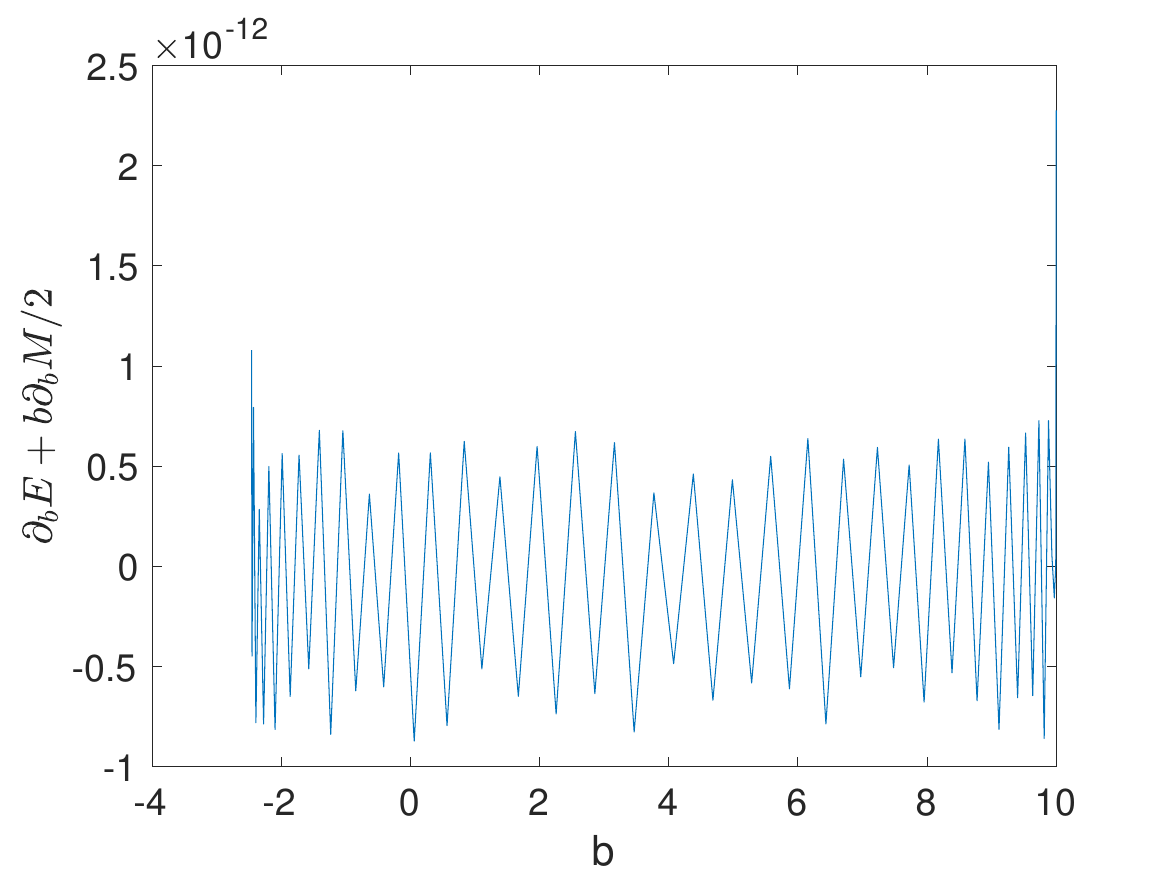}
\subcaption[]{{\footnotesize checking \eqref{E:EM-check}}}
\end{subfigure}\\
\begin{subfigure}{.32\textwidth}
\includegraphics[width=1\linewidth,height=0.75\linewidth]{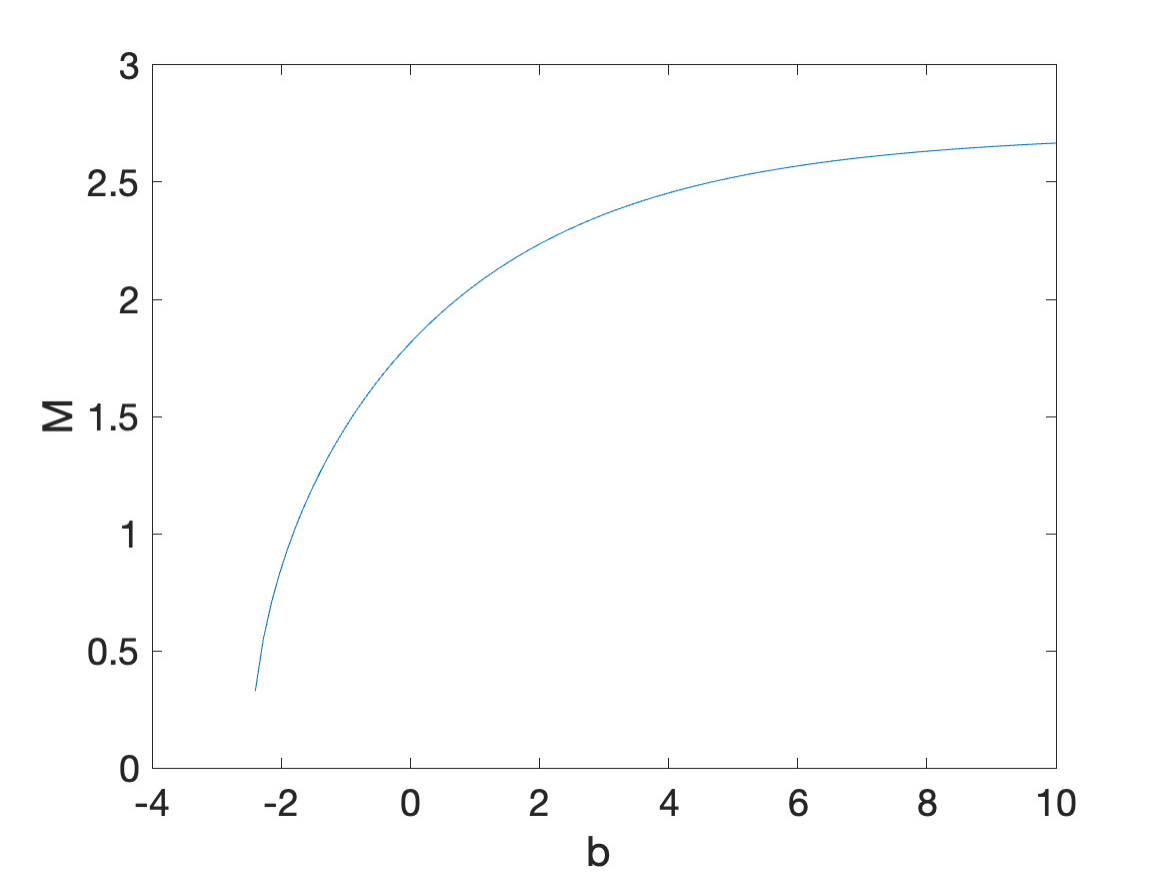}
\subcaption[]{{\footnotesize $\alpha=4$, $M(Q_b)\to 2.72$}}
\end{subfigure}
\begin{subfigure}{.32\textwidth}
\includegraphics[width=1\linewidth,height=0.75\linewidth]{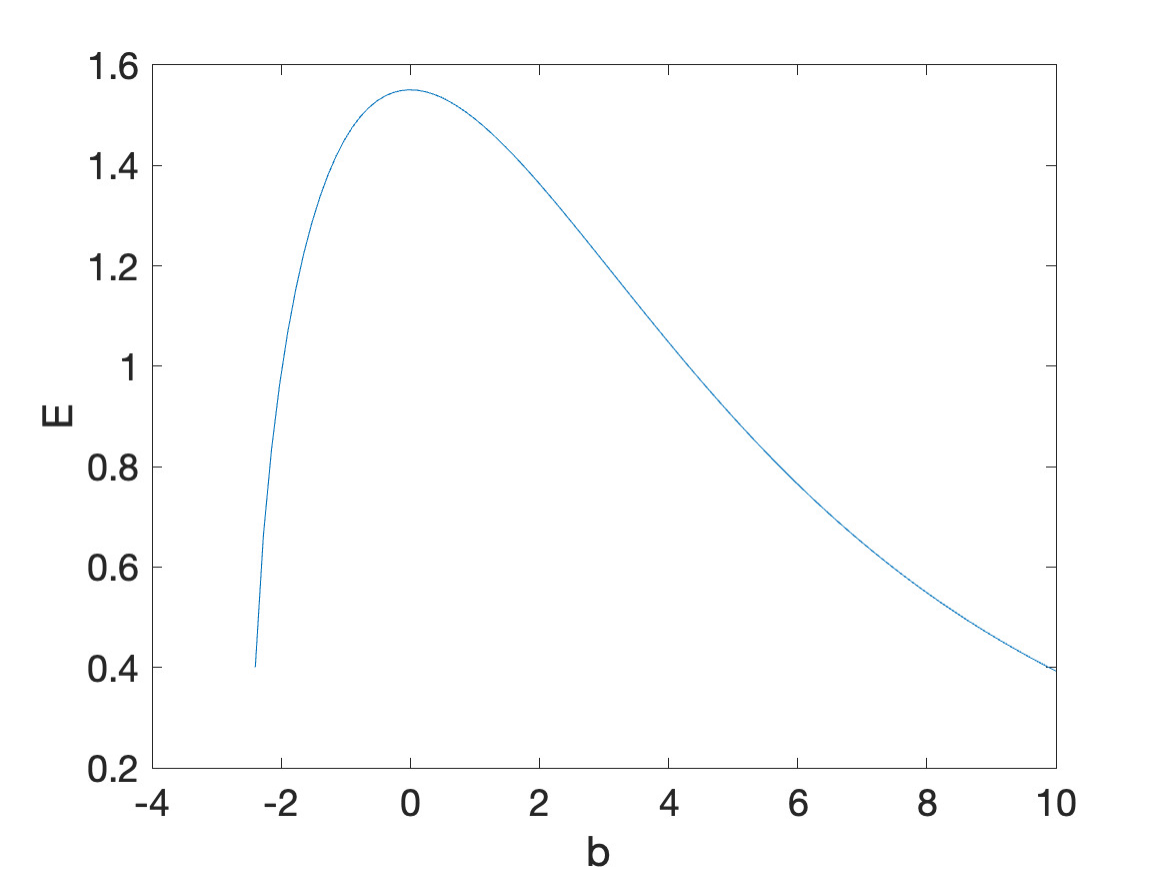}
\subcaption[]{{\footnotesize $E(Q_b)$ as function of $b$}}
\end{subfigure}
\begin{subfigure}{.32\textwidth}
\includegraphics[width=1\linewidth,height=0.75\linewidth]{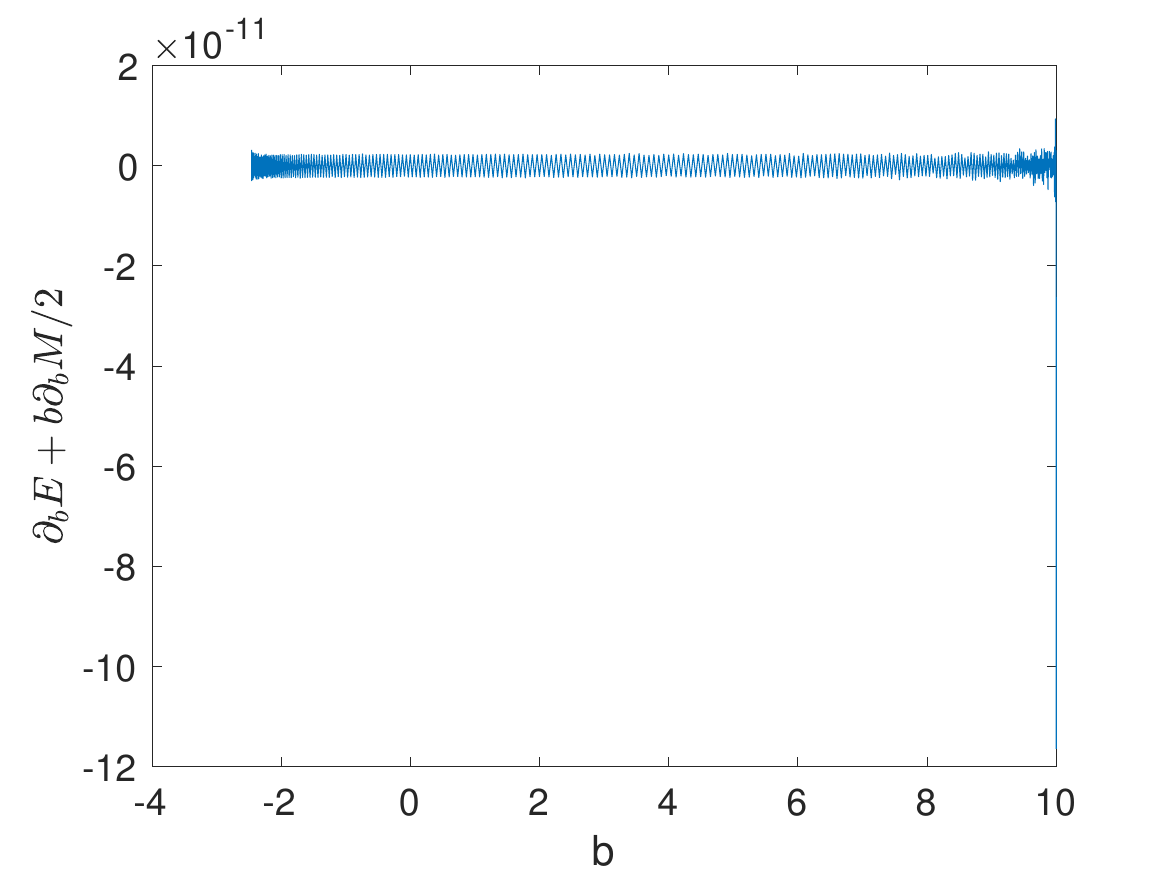} 
\subcaption[]{{\footnotesize checking \eqref{E:EM-check}}}
\end{subfigure}
\caption{\small 1D NLS: subcritical, $\alpha=2$ (top) and critical, $\alpha=4$ (bottom). Mass and energy dependence on $b$
(left and middle columns). Verification of \eqref{E:EM-check}, which is at least on the order of $10^{-11}$(right column).}
\label{F:ME-1D-p3-5}
\end{figure}

\subsubsection{1D case}

\begin{figure}[!htb]
\begin{subfigure}{.32\textwidth}
\includegraphics[width=1\linewidth,height=0.85\linewidth]{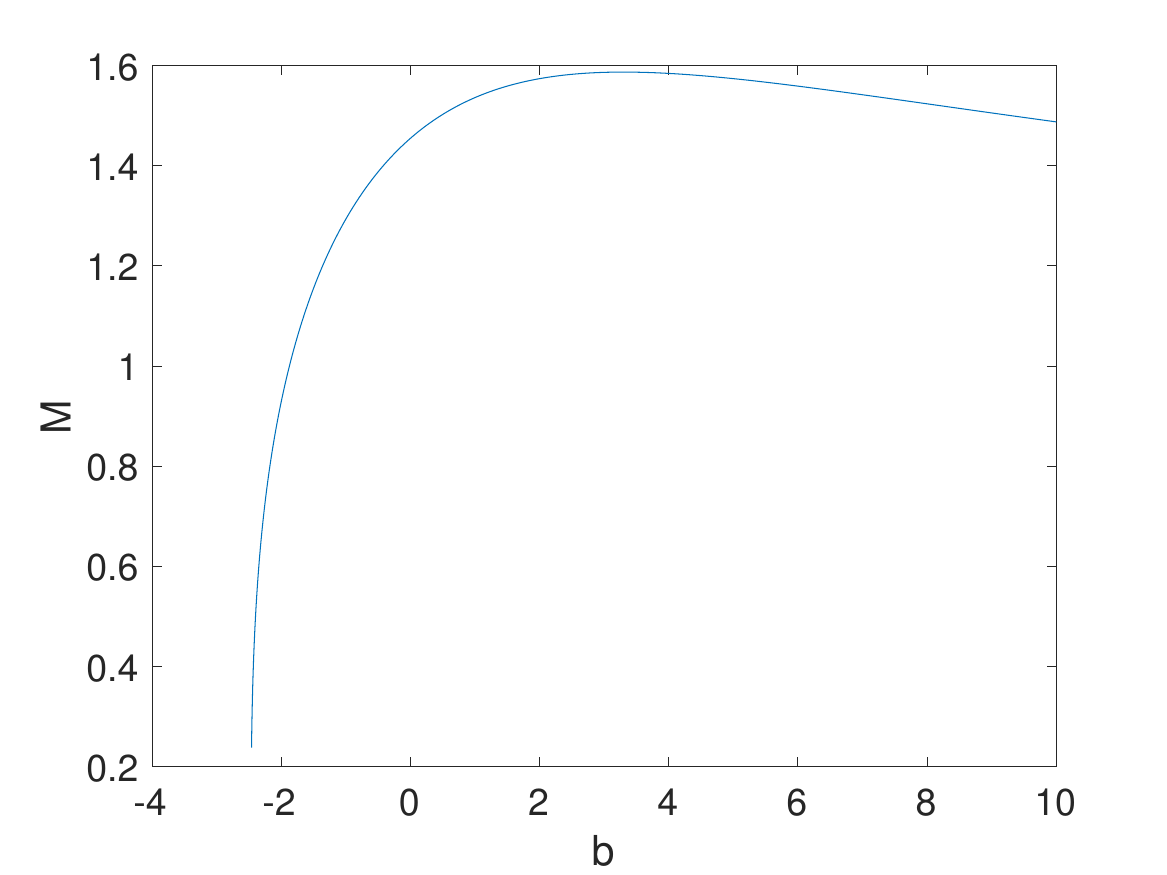}
\subcaption[]{{\footnotesize $\alpha=6$, $M(Q_b) \searrow 0$}}
\end{subfigure}
\begin{subfigure}{.32\textwidth}
\includegraphics[width=1\linewidth,height=0.85\linewidth]{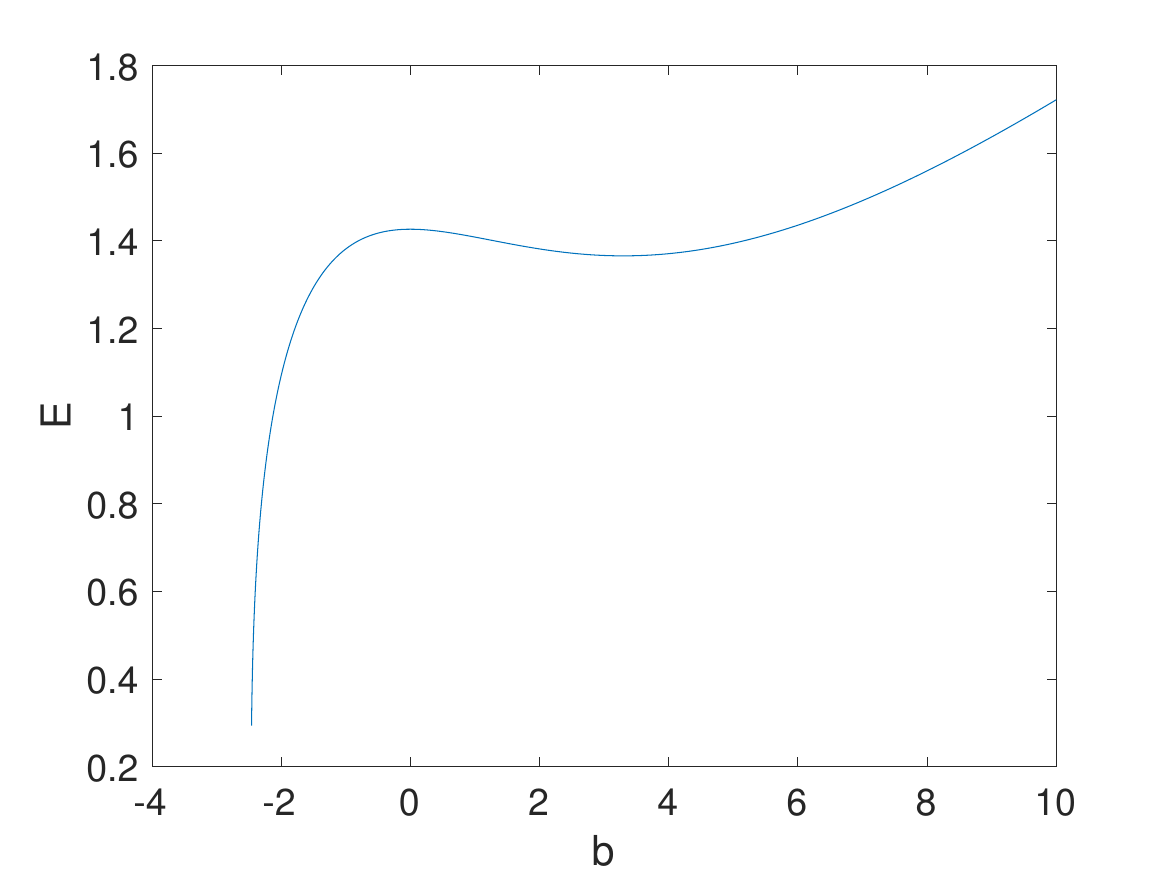}
\subcaption[]{{\footnotesize $E(Q_b)$ as function of $b$}}
\end{subfigure}
\begin{subfigure}{.32\textwidth}
\includegraphics[width=1\linewidth,height=0.85\linewidth]{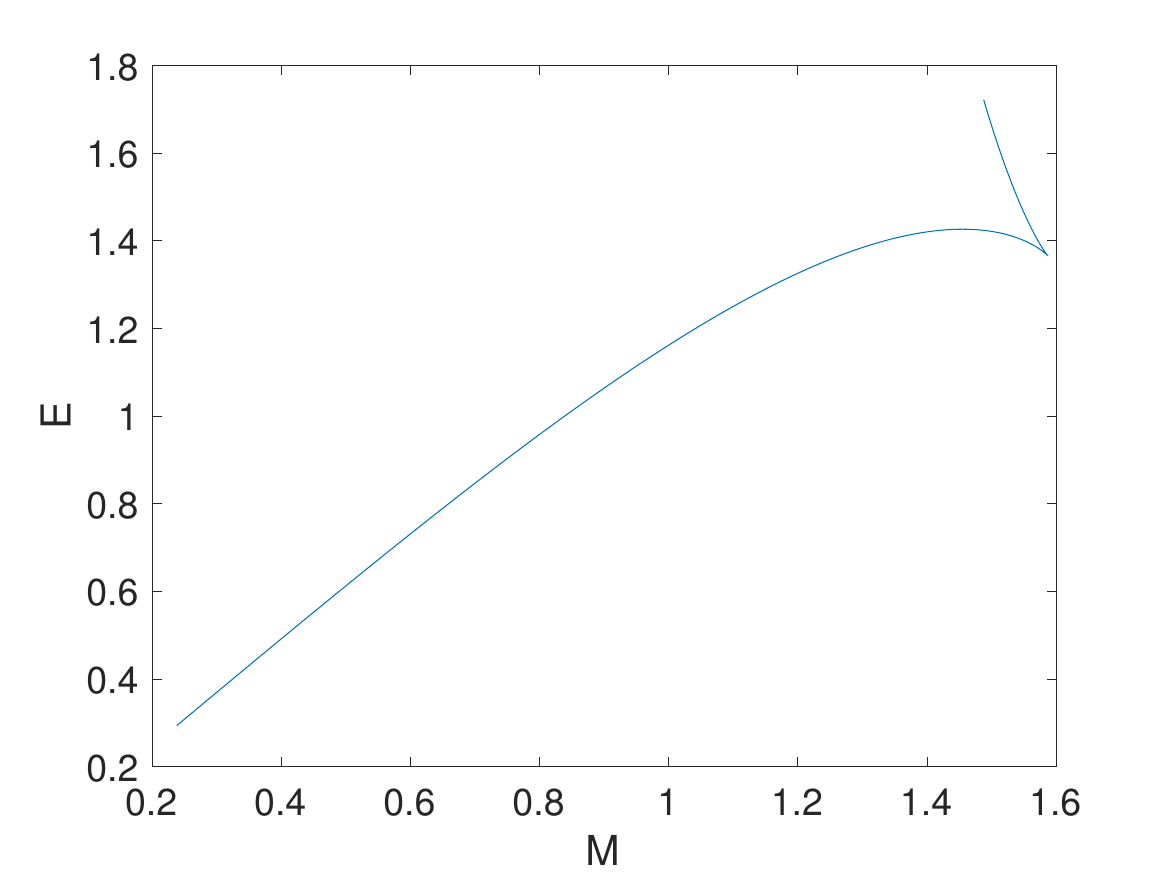}
\subcaption[]{{\footnotesize $E=E(M)$}}
\end{subfigure}\\
\begin{subfigure}{.32\textwidth}
\includegraphics[width=1\linewidth,height=0.85\linewidth]{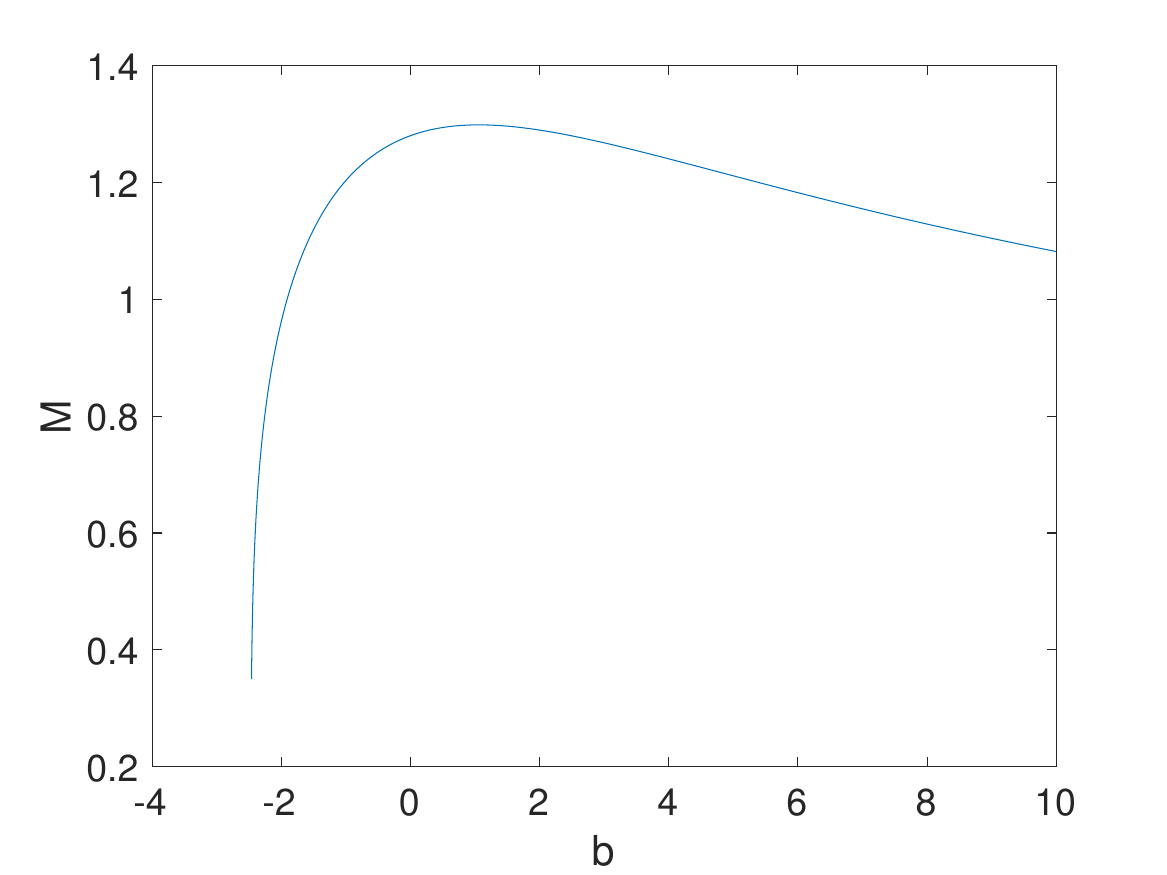}
\subcaption[]{{\footnotesize $\alpha=8$, $M(Q_b) \searrow 0$}}
\end{subfigure}
\begin{subfigure}{.32\textwidth}
\includegraphics[width=1\linewidth,height=0.85\linewidth]{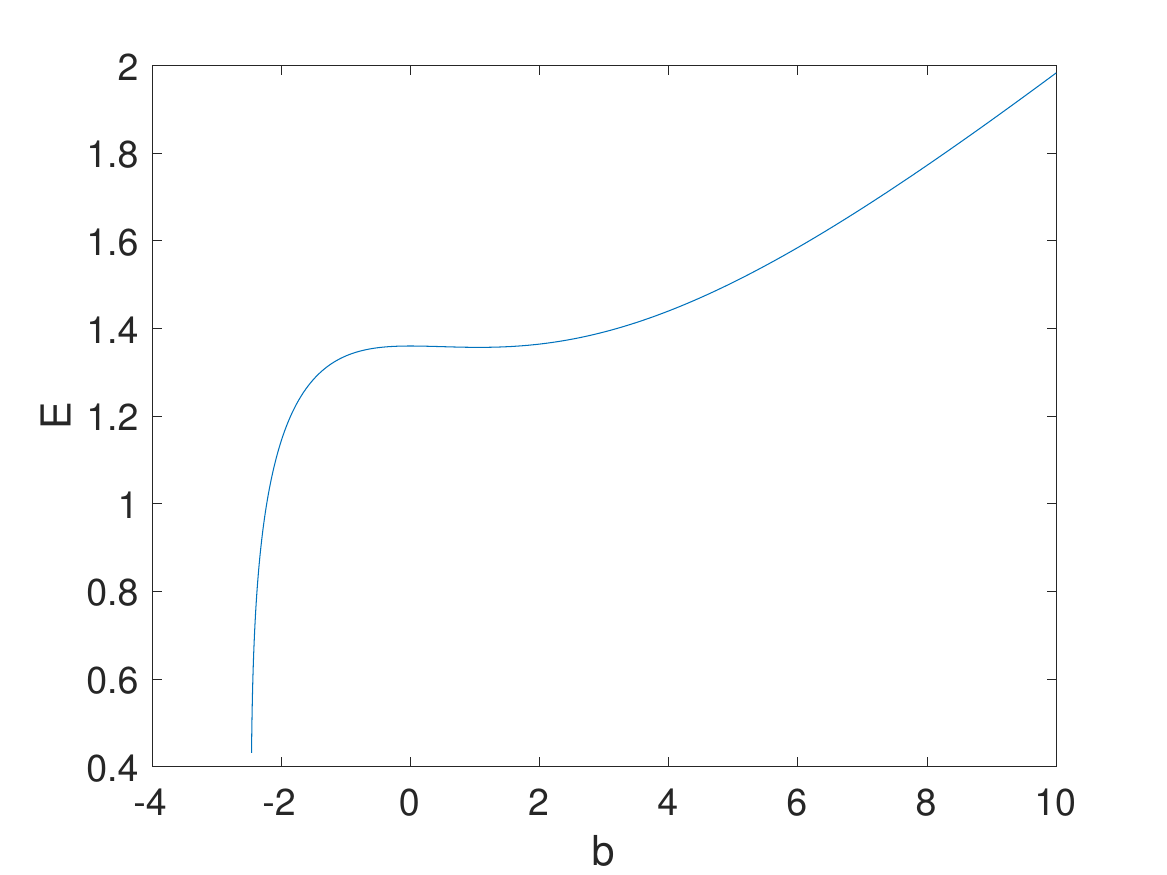}
\subcaption[]{{\footnotesize $E(Q_b)$ as function of $b$}}
\end{subfigure}
\begin{subfigure}{.32\textwidth}
\includegraphics[width=1\linewidth,height=0.85\linewidth]{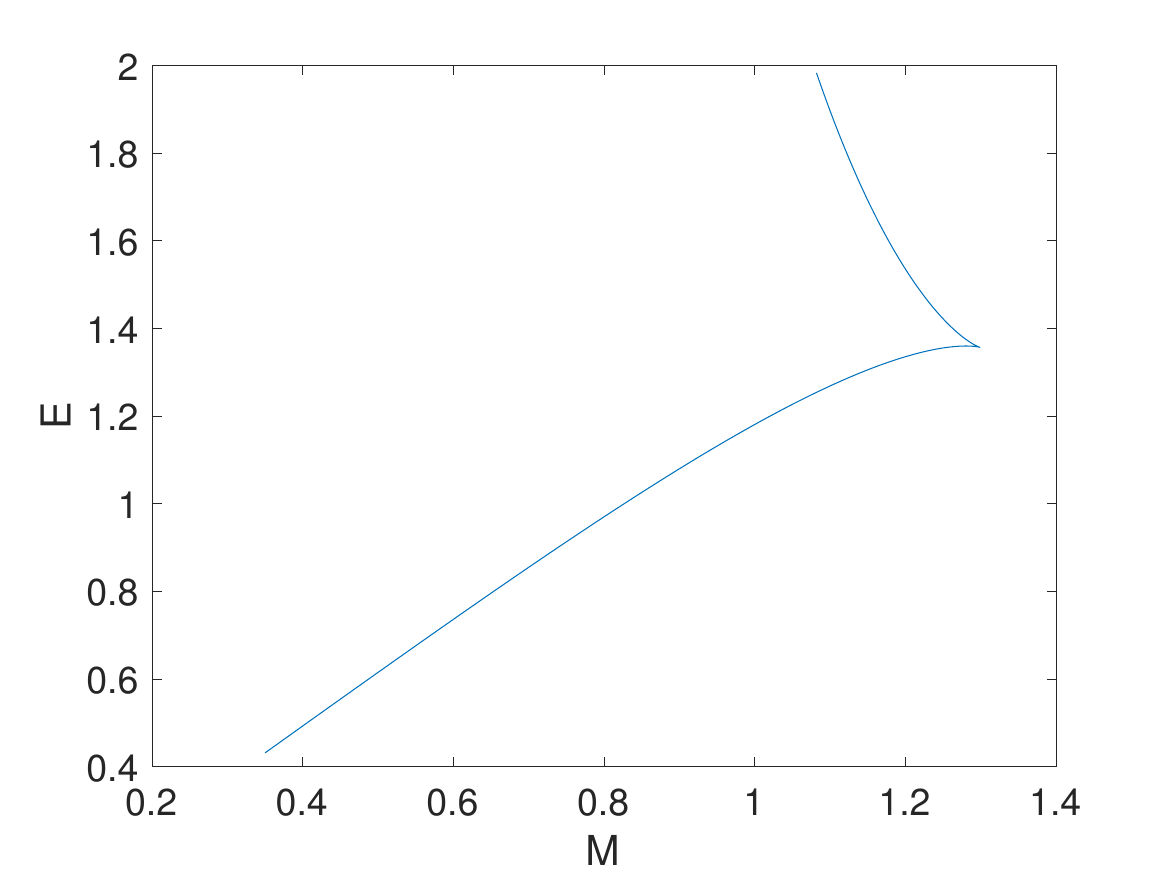}
\subcaption[]{{\footnotesize $E=E(M)$}}
\end{subfigure}
\caption{\small 1D NLS, supercritical cases: $\alpha=6$ (top) and $\alpha=8$ (bottom). Mass and energy dependence on $b$  (left and middle columns). Dependence of energy as a function of mass, $E = E(M)$ (right column). The maximum values of the mass as a function of $b$ are reached at  $b^\ast \approx 3.3117$ (for $\alpha=6$) and $b^\ast \approx 1.0718$ (for $\alpha=8$). This is the point where the split into two branches occurs in the right column.} 
\label{F:ME-1D-p7-9}
\end{figure}

By Lemma \ref{L:conv}(II), $\widetilde Q_b\to\mathcal R_\alpha$, and therefore, $M(Q_b)\sim b^{2/\alpha-d/2}M(\mathcal R_\alpha)$ as $b\to\infty$. In 1D this exponent is $(4-\alpha)/(2\alpha)$, thus, in the critical case $\alpha=4$, the mass of $Q_b$ approaches $M(\mathcal R_4)\approx2.72$, see \eqref{E:MR}, Figure \ref{F:R-MR}, and the left plot of Figure \ref{F:cloud}. In the subcritical case it grows to infinity, as in Figure \ref{F:ME-1D-p3-5}(A), while in the supercritical case it tends to zero, as in Figure \ref{F:ME-1D-p7-9} (left column).

In the middle column of Figure \ref{F:ME-1D-p3-5} we show how energy changes with respect to $b$ and in the right column we provide verification of how the slope of the mass is connected with the slope of energy in \eqref{E:EM-check}, which are on the order of $10^{-11}$ or better. In fact, we checked the identity \eqref{E:EM-check} in all other computed cases and it holds at least to the same accuracy. 

In Figure \ref{F:ME-1D-p7-9} we consider 1D supercritical cases and show the dependence of the mass and energy on $b$ in the left and middle columns. The mass is no longer monotone: it initially grows, reaches a maximum, and then decays. We therefore also plot the dependence of energy on mass in the right column and observe two branches. Thus, for the same mass there are {\bf two} computed ground-state profiles with different values of $b$; the branch with positive mass slope (and the lower energy) is the one selected as {\it stable} by the slope criterion and by the time evolutions below, whereas the branch with negative mass slope is {\it unstable}. 
\smallskip

We now use the approach of Kfoury-LeCoz-Tsai \cite{KLT} to find turning points in the supercritical case. Since we have boundary terms, it is not a straightforward adaptation, as we need to deal with them, so in our functions $I_j$ we have extra terms. The nonzero boundary contribution is retained explicitly through the parameter $\nu$ and the length condition below. This result shows how to locate the turning point $b^*$ analytically via a scalar integral equation. 

\begin{proposition}[Turning point $b^*$ in 1D]\label{P:1D-b*}
Let $\alpha>0$, and for $b>-\frac{\pi^2}4$ let $Q_b$ be the positive, even
ground state in 1D, satisfying
$$
-Q_b''+bQ_b-Q_b^{\alpha+1}=0 \quad\text{on }(-1,1),
\quad \mbox{with} \quad
Q_b(-1)=Q_b(1)=0.
$$
Denote
$$
A=Q_b(0) \quad \mbox{and} \quad \nu=\tfrac{(\alpha+2)}{2} \,b \,A^{-\alpha}.
$$
Then $\nu<1$.
(The parameter $\nu$ has the same sign as $b$. In particular,
$0<\nu<1$ only when $b>0$.)

Fix $\alpha>0$ and define $I_0(\nu)$ and $I_2(\nu)$ by
\begin{equation}\label{E:Ij}
I_j(\nu) = \int_0^1 \frac{s^j\,ds}
{\sqrt{1-s^{\alpha+2}-\nu(1-s^2)}},
\quad j=0,2,\quad \nu<1.
\end{equation}
For a given $b$, $\nu$ is defined implicitly via $b = \nu \big(I_0(\nu)\big)^2$, which has a unique solution $v<1$, and then the peak amplitude is given by
\begin{equation}\label{E:1D-A}
A \equiv Q_b(0) =\big(\tfrac{\alpha+2}{2}I_0(\nu)^2 \big)^{\frac1{\alpha}}.
\end{equation}

For $\alpha>4$ (supercritical case) the mass curve $b\mapsto M(Q_b)$
attains its global maximum at some frequency ~$b^*\in(-\frac{\pi^2}4,\infty)$, which is given by
\begin{equation}\label{E:1D-bstar}
b^*=\nu^* (I_0(\nu^*))^2,
\end{equation}
where $\nu^*<1$ solves
\begin{equation}\label{E:1D-b*-eq}
\frac{I_2'(\nu)}{I_2(\nu)}
= \Big(1-\frac4\alpha \Big) \frac{I_0'(\nu)}{I_0(\nu)}.
\end{equation}
Conversely, every solution of \eqref{E:1D-b*-eq} determines
a stationary point of the mass curve through \eqref{E:1D-bstar}.

At $b^\ast$, we have 
\begin{equation}\label{E:M*}
M(Q_{b^*}) = 2 \Big(\frac{\alpha+2}{2} \Big)^{\frac2{\alpha}}
\big(I_0(\nu^*) \big)^{\frac{4}{\alpha}-1} I_2(\nu^*).
\end{equation}

\end{proposition}

\begin{proof}
Recall the ground state equation $-Q_b'' + b Q_b - Q_b^{\alpha+1} = 0$ with $Q_b' (0) = 0$ and $Q_b(1)=0$. Multiplying it by $Q_b'$ and integrating from $0$ to $x \in [0,1]$, we obtain
\begin{equation}\label{E:1D-integralQ}
\left[ -\frac12 (Q_b'(y))^2 + \frac{b}2 (Q_b(y))^2 - \frac1{\alpha+2} (Q_b(y))^{\alpha+2} \right]_{y=0}^{y=x} = 0.
\end{equation}
Using $Q'_b(0) = 0$ and $Q_b(0)=A$ gives 
\begin{equation}\label{E:Qb}
(Q_b'(x))^2  = {b}\, (Q_b(x))^2 - \frac2{\alpha+2} (Q_b(x))^{\alpha+2} + \Big(\frac2{\alpha+2} A^{\alpha+2} - b \, A^2 \Big).
\end{equation}
If we substitute $x=1$ into the above and recall the zero boundary condition, then we get 
$$
(Q_b'(1))^2  = \frac2{\alpha+2} A^{\alpha+2} - b \, A^2.
$$
Since $(Q_b'(1))^2 >0$, we deduce that 
$\frac2{\alpha+2} A^{\alpha+2} - b \, A^2>0$, and hence, $\nu: = \frac{\alpha+2}2 \, b  A^{-\alpha} <1$. 

We proceed as in \cite[Lemma 4]{KLT}, we rewrite \eqref{E:Qb} (choosing the negative sign, since $Q_b$ is a decreasing function on $(0,1)$) as
$$
Q_b'(x)  = - \sqrt {{b}\, (Q_b(x))^2 - \tfrac2{\alpha+2} (Q_b(x))^{\alpha+2} +  \tfrac2{\alpha+2} A^{\alpha+2} \big(1-\nu) }.
$$
Let $t = Q_b(x)$, then 
\begin{align}
dx = \frac{dt}{Q_b'(x)} & = -\frac{dt}{\sqrt{ {b}\, (Q_b(x))^2 - \tfrac2{\alpha+2} (Q_b(x))^{\alpha+2} +  \tfrac2{\alpha+2} A^{\alpha+2} \big(1-\nu) }} \notag \\
& = -\frac{dt}{\sqrt{ {b}\, t^2 - \tfrac2{\alpha+2} t^{\alpha+2} +  \tfrac2{\alpha+2} A^{\alpha+2} \big(1-\nu) }} \label{E:dt}.
\end{align}
Noting that $x=0$ changes to $t=Q_b(0)=A$ and $x=1$ to $t=Q_b(1)=0$, integrating the left-hand side of \eqref{E:dt} from $0$ to $1$ simply gives length of the interval $[0,1]$, and thus, we obtain
$$
1 = \int_0^1 dx = \int_0^A \frac{dt}{\sqrt{ {b}\, t^2 - \tfrac2{\alpha+2} t^{\alpha+2} +  \tfrac2{\alpha+2} A^{\alpha+2} \big(1-\nu) }}.
$$
Similarly, we can express mass in a new variable 
$$
M(Q_b) = 2 \int_0^1 |Q_b(x)|^2 \, dx = 2 \int_0^A t^2 \, \frac{dt}{\sqrt{ {b}\, t^2 - \tfrac2{\alpha+2} t^{\alpha+2} +  \tfrac2{\alpha+2} A^{\alpha+2} \big(1-\nu) }}.
$$
Renormalizing $t = As$ (and thus, $0 \leq s \leq 1$), the last two integrals become
\begin{align}
1 &= \frac1{\sqrt{\tfrac2{\alpha+2} A^{\alpha} }} \int_0^1 \frac{ds}{\sqrt{ (1-s^{\alpha+2}) -   \nu(1-s^2)}} = \sqrt{\tfrac{\alpha+2}2} A^{-\frac{\alpha}2} \, {I_0}, \label{E:length} \\
M(Q_b) &= \frac{2 A^2}{\sqrt{\tfrac2{\alpha+2} A^{\alpha} }} \int_0^1 \frac{s^2 \,ds}{\sqrt{ (1-s^{\alpha+2}) -   \nu(1-s^2)}} = 2 \sqrt{\tfrac{\alpha+2}2} A^{2-\frac{\alpha}2} I_2,
\end{align}
where $I_0$ and $I_2$ were defined in \eqref{E:Ij}.

In particular, solving for $A$ in the first integral, gives \eqref{E:1D-A}, and substituting $A$ into the second, yields 
\begin{equation}\label{E:Mb}
M(Q_b) = 2 \Big( \frac{\alpha+2}{2} \Big)^{\frac2{\alpha}}I_0^{\frac4{\alpha}-1} I_2,
\end{equation}
which gives \eqref{E:M*} once we substitute $\nu^*$, which we find next. 

From the length equation \eqref{E:length} and the definition of $\nu$, we have  $b = \nu \big(I_0(\nu) \big)^2$.  We investigate the dependence $b(\nu)$ in more detail. 
Note that $I_0$ and $I_2$ are finite and positive for $\nu <1$ (and a fixed $\alpha>0$). 
Furthermore, factoring $\sqrt{-\nu}$ out (for negative $\nu$), we get
$$
I_0(\nu) = \frac1{\sqrt{-\nu}} \int_0^1 \frac{ds}{\sqrt{ 1-s^2+\frac{1-s^{\alpha+2}}{-\nu} } },
$$
and thus, as $\nu \to -\infty$, we obtain
$$
\sqrt{-\nu} I_0(\nu) \to \int_0^1 \frac{ds}{\sqrt{1-s^2}} = \frac{\pi}{2}, 
$$
implying that $b(\nu) \equiv - \big(\sqrt{-\nu} I_0(\nu)\big)^2 \to -\frac{\pi^2}{4}$.

For the other endpoint, as $\nu \to 1$, the denominator in the integrand of $I_0$ goes to 
$s\sqrt{1-s^\alpha}$, and thus, 
$$
I_0(\nu) \to \int_0^1 \frac{ds}{s\sqrt{1-s^{\alpha}}} = \infty \quad \mbox{as} ~~\nu \to 1.
$$
Thus, the function $b(\nu) \to -{\pi^2}/4$ as $\nu \to -\infty$ and $b(\nu) \to \infty$ as $\nu \to 1$, so $b: (-\infty,1) \to (-\frac{\pi^2}{4},\infty)$.
 
Now, differentiate $M(Q_b)$ with respect to $b$ to find the turning point, which amounts to 
\begin{equation}\label{E:derM}
\partial_b M(Q_b) = \frac{{\partial M(Q_b)}/{\partial \nu}}{{\partial b}/{\partial \nu}}.
\end{equation}
Differentiating the right-hand side in \eqref{E:Mb} with respect to $\nu$, we get 
\begin{equation}\label{E:dnu-M}
\partial_{\nu} M(Q_b) = 2 \Big( \frac{\alpha+2}{2} \Big)^{\frac2{\alpha}} \,  I_0^{\frac4{\alpha}-2} \left[ \Big({\frac4{\alpha}-1} \Big)  I_0' I_2+ I_0 I_2' \right].
\end{equation}
From the length equation \eqref{E:length} and definition of $b(\nu)$, we get 
$$
\frac{\partial b}{\partial \nu} = I_0^2 + 2\nu I_0 I_0' 
=  I_0 \int_0^1 \frac{(1-s^{\alpha+2}) \,ds}{ [ (1-s^{\alpha+2}) -   \nu(1-s^2) ]^{3/2} } > 0, \quad \nu < 1.
$$ 
(The essential fact is that it is not equal to zero for the purpose of \eqref{E:derM}). Together with the endpoint limits, positivity shows that $\nu \mapsto b(\nu)$ is a strictly increasing bijection from $(-\infty,1)$ onto $(-\pi^2/4,\infty)$. Putting these derivatives into \eqref{E:derM}, we obtain $\partial_b M(Q_b)$.

Then $\partial_b M(Q_b) = 0$ is equivalent to the expression in the square brackets in \eqref{E:dnu-M} being zero, or $({1-\frac4{\alpha}})I_0' I_2 =  I_0 I_2'$, which determines $\nu^*$ via \eqref{E:1D-b*-eq}.

Finally, note that for $\alpha>4$, from the definition of $I_0$ and $I_2$, as well as \eqref{E:Mb}, it follows that $M(Q_b) \to 0$ as $b \to -\pi^2/4$ and $b \to \infty$. From the continuity (and positivity) of mass, we can conclude that it attains a global maximum at an interior point $b^*$ (where its derivative vanishes). 
Conversely, every value $\nu^*<1$ satisfying \eqref{E:1D-b*-eq} defines $b^* = \nu^* (I_0(\nu^*))^2 > -\pi^2/4$, and thus, corresponds to the ground state $Q_{b^*}$, with \eqref{E:derM} making its mass derivative vanish. 
\end{proof}

\smallskip

\begin{remark}[An explicit equation for $\nu^*$]
\label{R:1D-bstar-integrals}
The substitution $s=\sin\theta$ gives
$$
I_0(\nu) = \int_0^{\pi/2}
\left[\frac{1-\sin^{\alpha+2}\theta}{\cos^2\theta}-\nu \right]^{-1/2} \,d\theta,
$$
and
$$
I_2(\nu) = \int_0^{\pi/2} \sin^2\theta
\left[\frac{1-\sin^{\alpha+2}\theta}{\cos^2\theta}-\nu \right]^{-1/2} \,d\theta,
$$
with the integrals defined for any $\alpha>0$ and $\nu<1$, and the integrands continuous on $(0,\frac{\pi}2)$. 

When $\alpha$ is an even positive integer, the quotient in square brackets simplifies to a sum of even powers of sines up to $\alpha$, and denoting by  
$S(\theta):=1+\sin^2\theta+\sin^4\theta + \cdots+\sin^\alpha\theta$, we get
\begin{align}\label{E:I0-I2}
I_0(\nu)
&=\int_0^{\pi/2}\frac{d\theta}{\bigl[S(\theta)-\nu\bigr]^{1/2}},
\qquad 
I_2(\nu) =\int_0^{\pi/2}
\frac{\sin^2\theta\,d\theta}{\bigl[S(\theta)-\nu\bigr]^{1/2}}, \\[.2mm]
I_0'(\nu) 
&=\frac12\int_0^{\pi/2}
\frac{d\theta}{\bigl[S(\theta)-\nu\bigr]^{3/2}},
\qquad
I_2'(\nu)=\frac12\int_0^{\pi/2}
\frac{\sin^2\theta\,d\theta}{\bigl[S(\theta)-\nu\bigr]^{3/2}}. \notag
\end{align}
Thus, for even $\alpha>4$ and $\nu<1$, \eqref{E:1D-b*-eq} becomes
\begin{equation*} 
\begin{aligned}
&\frac{\ds \int_0^{\pi/2}\frac{\sin^2 \theta \, d\theta}
{\bigl[S(\theta)-\nu \bigr]^{3/2}}}{\ds \int_0^{\pi/2} \frac{\sin^2 \theta \, d\theta}
{\bigl[S(\theta)-\nu \bigr]^{1/2}}}
= \Big(1-\frac4{\alpha} \Big)
\frac{\ds\int_0^{\pi/2}\frac{d\theta}{\bigl[S(\theta) - \nu \bigr]^{3/2}}
}{\ds \int_0^{\pi/2} \frac{d\theta}
{\bigl[S(\theta)-\nu \bigr]^{1/2}}}.
\end{aligned}
\end{equation*}
Cross-multiplying and keeping the derivative integrals as outer integration, we obtain
\begin{equation}\label{E:1D-b*-3}
\int_0^{\pi/2} \frac{I_0(\nu)\sin^2\theta
-\left(1-\frac4\alpha\right)I_2(\nu)}{\bigl[S(\theta)-\nu \bigr]^{3/2}}
\,d\theta=0,
\end{equation}
with the left-hand side being an error, or a residual, when searching for  $\nu^*$.
 
Thus, the process of finding $\nu^*$ requires computing three integrals: $I_0$ and $I_2$ from \eqref{E:I0-I2}, and the integral in \eqref{E:1D-b*-3}.
\end{remark}

\begin{remark}[Numerical computation of $b^*$ for even $\alpha>4$]
\label{R:1D-bstar-computation}
Here we outline the steps to find $b^*$ for even integer cases $\alpha= 6,8, \dots$.
\begin{itemize}
\item Fix even $\alpha>4$ and take a trial value $\nu<1$. 

\item Compute two integrals $I_0(\nu)$ and $I_2(\nu)$ at the trial
value $\nu$. 

\item Substitute the computed values of $I_0$ and $I_2$ into \eqref{E:1D-b*-3}, compute that integral and check its sign and how close it is to zero. 

\item Proceed (e.g., apply bisection) to obtain $\nu^*$. 

\item Then $b^*=\nu^* \big(I_0(\nu^*)\big)^2$.
\end{itemize}
This procedure locates a stationary point. If several roots are found, compare their masses to identify a global maximum (note that we don't prove uniqueness).  
\end{remark}

\underline{Examples.}  For $\alpha=6$, 
Prop. \ref{P:1D-b*} (or following the outline in Remark \ref{R:1D-bstar-computation}) gives $\nu^* \approx 0.8538$, and thus, $b^* \approx 3.2945$ and $M(Q_{b^*}) \approx 1.5871$.
For $\alpha=8$, 
we obtain $\nu^* \approx   0.5274$,
and thus, $b^* \approx 1.0602$ and $M(Q_{b^*}) \approx 1.2986$.

The $b^*$ values can be compared with our numerical simulations shown in Figure \ref{F:ME-1D-p7-9}, and in particular, for $\alpha=6$ we found the maximum of mass occurs approximately at $b^* \approx 3.3117$, and for $\alpha=8$, we got $b^* \approx 1.0718$, which indicates consistency with the results we obtained from Prop. \ref{P:1D-b*}.

\subsubsection{2D case}

\begin{figure}[!htb]
\begin{subfigure}{.32\textwidth}
\includegraphics[width=1\linewidth,height=0.85\linewidth]{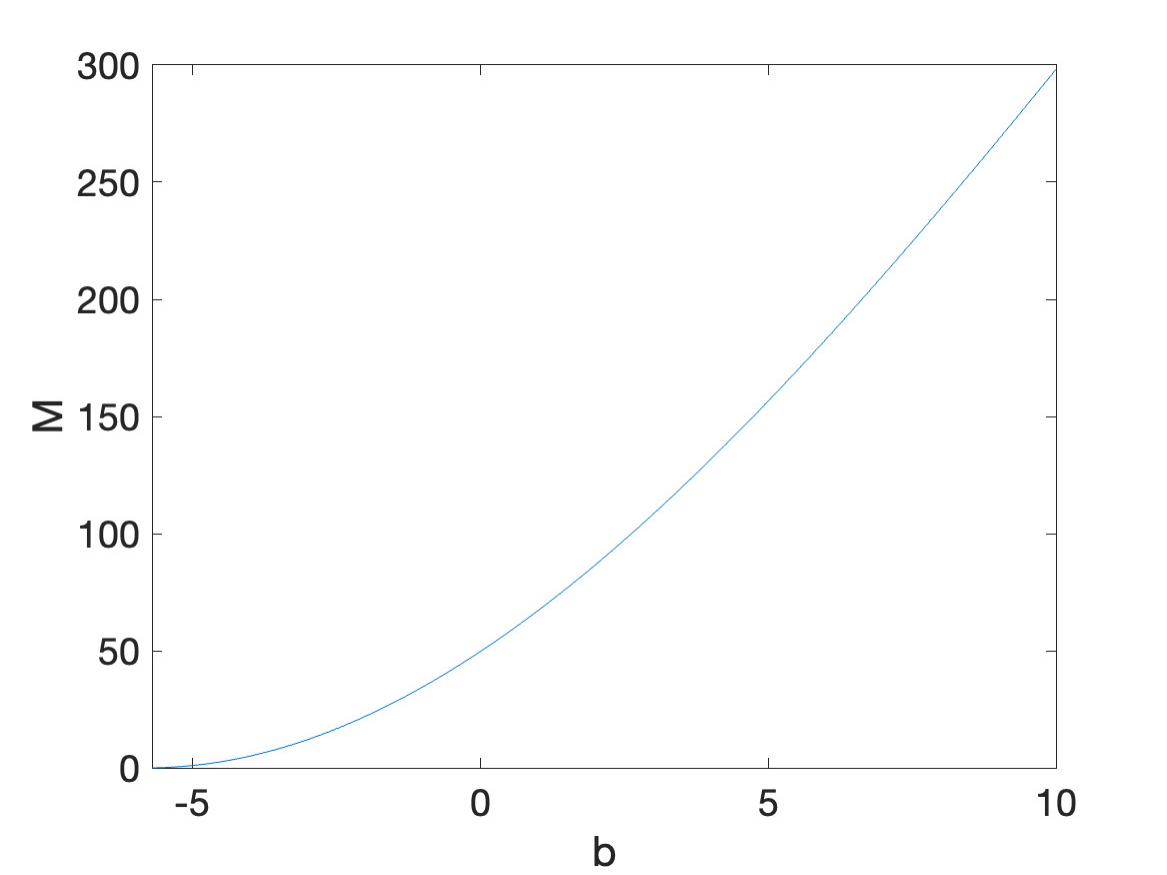}
\subcaption[]{{\footnotesize {$\alpha=1$, $M(Q_b) \nearrow \infty$}}}
\end{subfigure}
\begin{subfigure}{.32\textwidth}
\includegraphics[width=1\linewidth,height=0.85\linewidth]{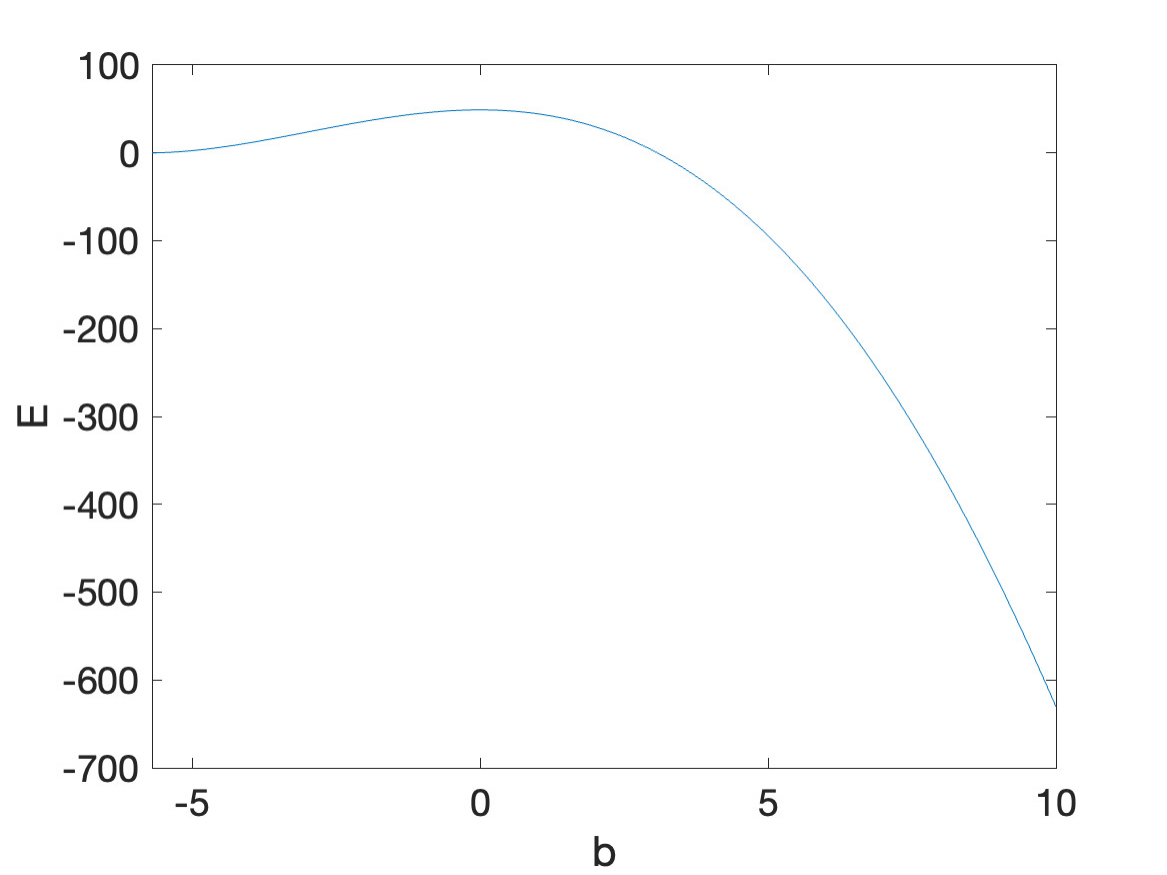}
\subcaption[]{{\footnotesize $E(Q_b)$ as function of $b$}}
\end{subfigure} 
\begin{subfigure}{.32\textwidth}
\includegraphics[width=1\linewidth,height=0.85\linewidth]{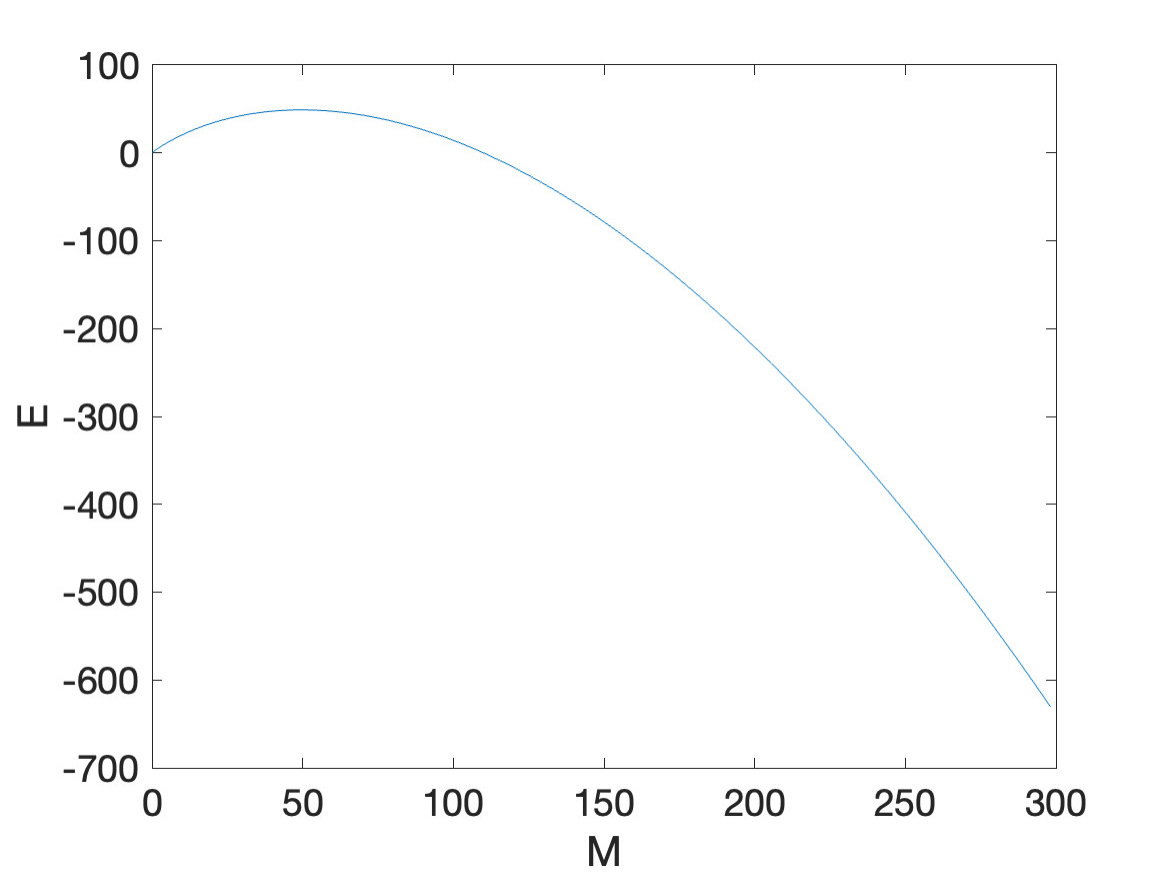}
\subcaption[]{{\footnotesize $E = E(M)$}}
\end{subfigure}\\
\begin{subfigure}{.32\textwidth}
\includegraphics[width=1\linewidth,height=0.85\linewidth]{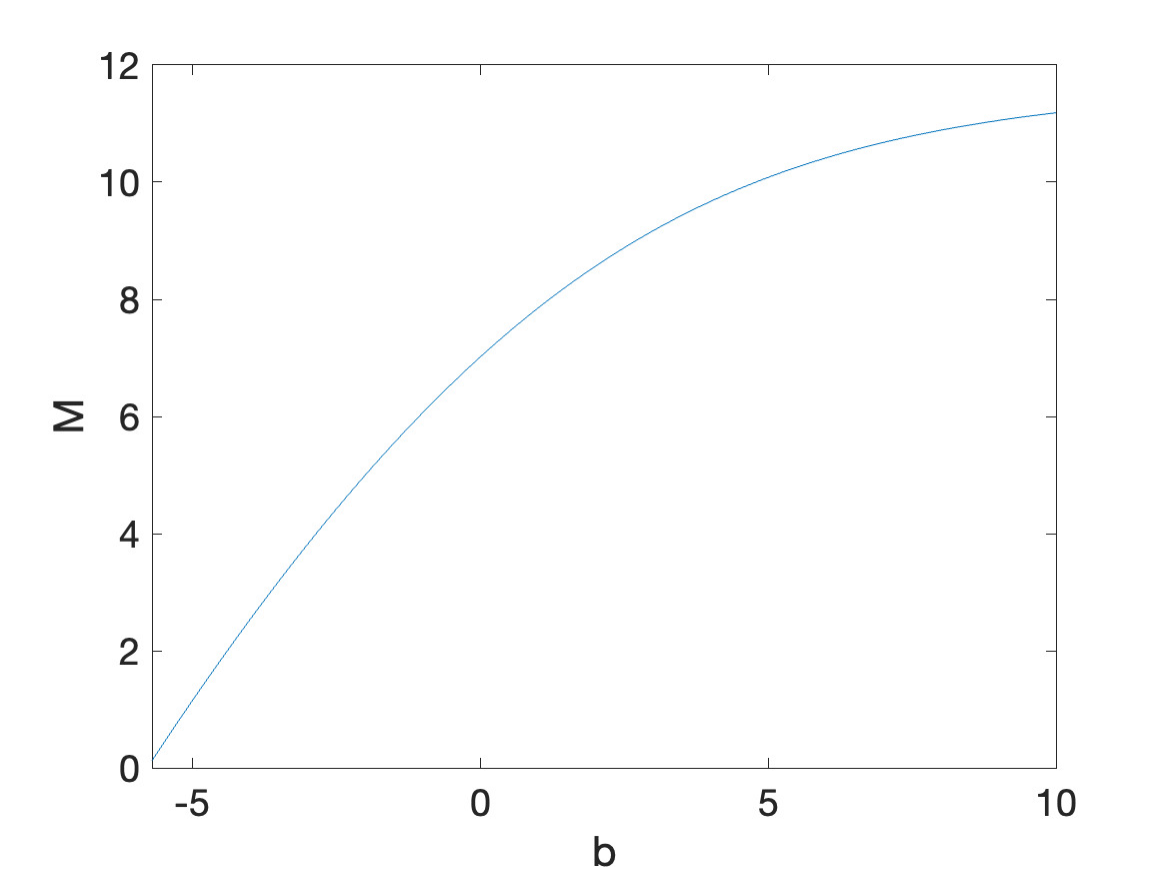}
\subcaption[]{{\footnotesize $\alpha=2$, $M(Q_b) \to 11.70$}}
\end{subfigure}
\begin{subfigure}{.32\textwidth}
\includegraphics[width=1\linewidth,height=0.85\linewidth]{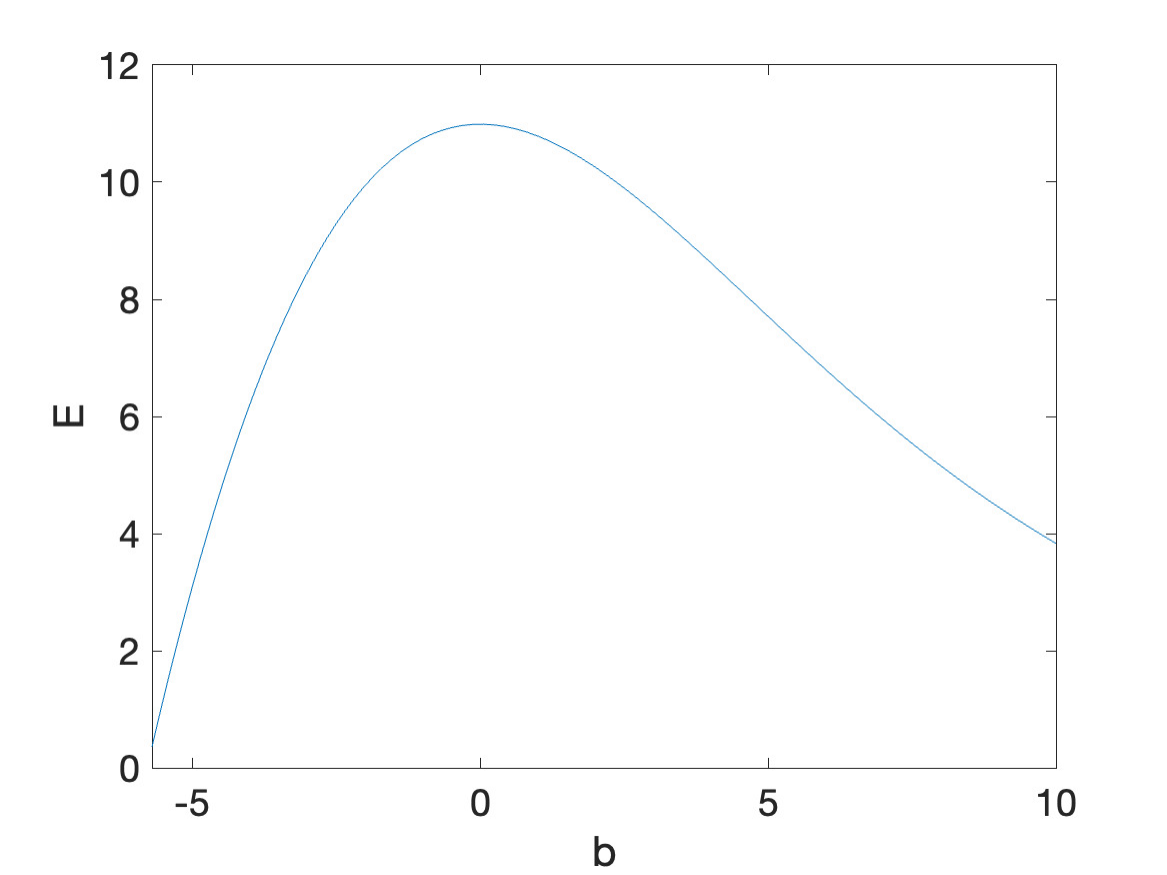}
\subcaption[]{{\footnotesize $E(Q_b)$ as function of $b$}}
\end{subfigure}
\begin{subfigure}{.32\textwidth}
\includegraphics[width=1\linewidth,height=0.85\linewidth]{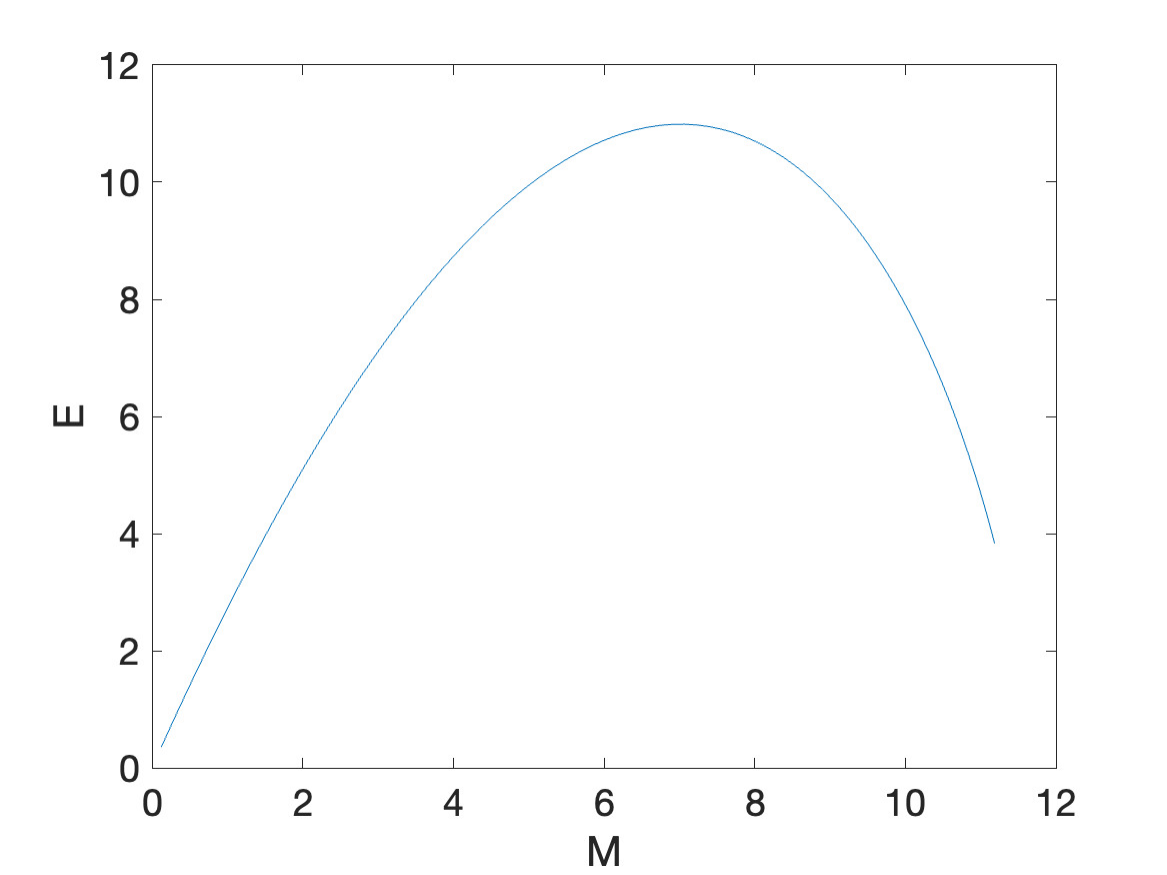} 
\subcaption[]{{\footnotesize $E=E(M)$}}
\end{subfigure}\\
\begin{subfigure}{.32\textwidth}
\includegraphics[width=1\linewidth,height=0.85\linewidth]{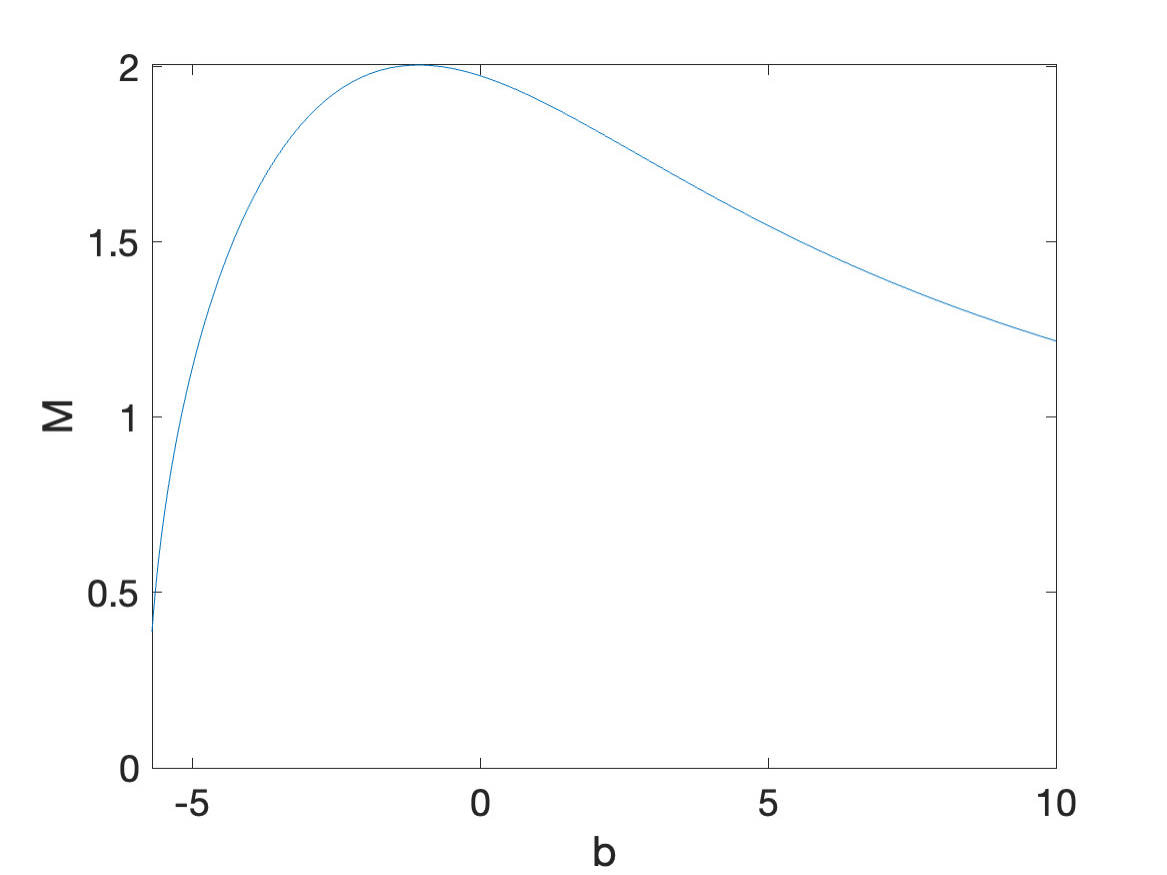}
\subcaption[]{{\footnotesize $\alpha=4$, $M(Q_b) \searrow 0$}}
\end{subfigure}
\begin{subfigure}{.32\textwidth}
\includegraphics[width=1\linewidth,height=0.85\linewidth]{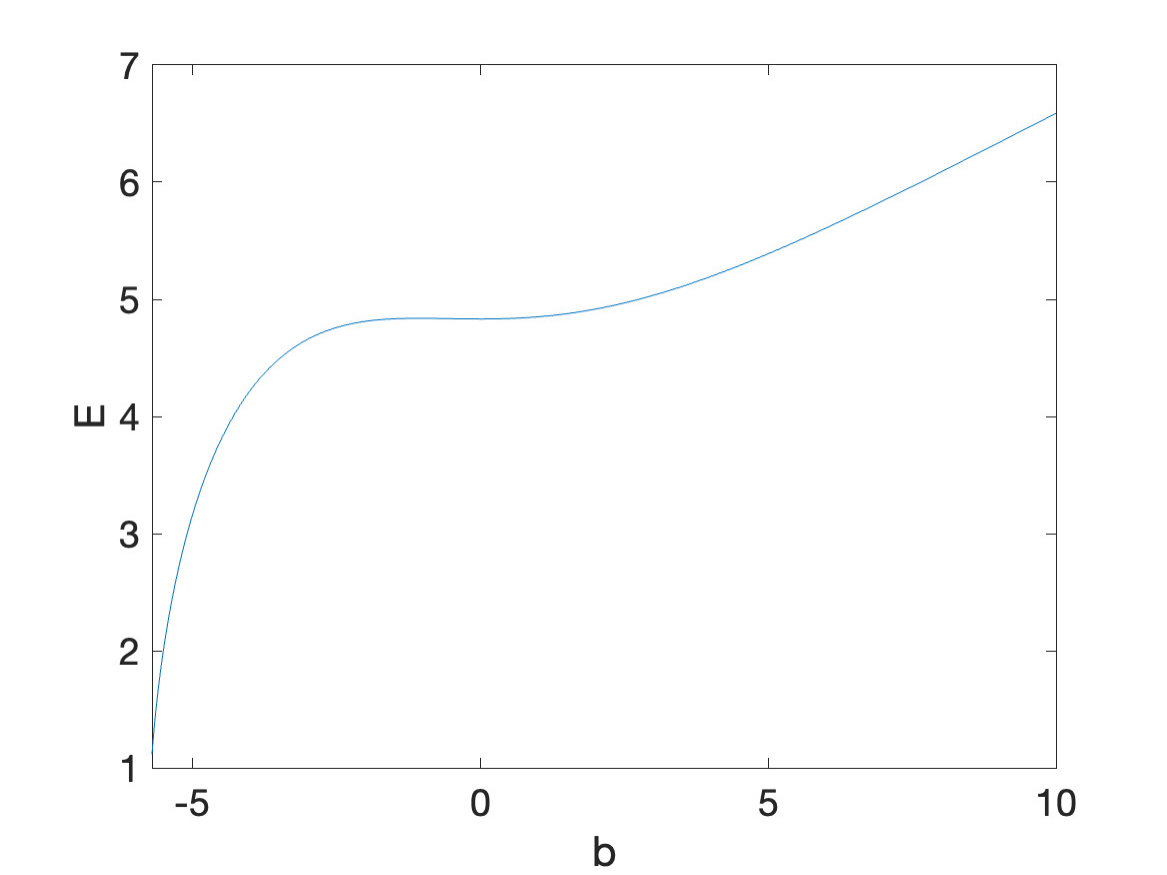}
\subcaption[]{{\footnotesize $E(Q_b)$ as function of $b$.}}
\end{subfigure}
\begin{subfigure}{.32\textwidth}
\includegraphics[width=1\linewidth,height=0.85\linewidth]{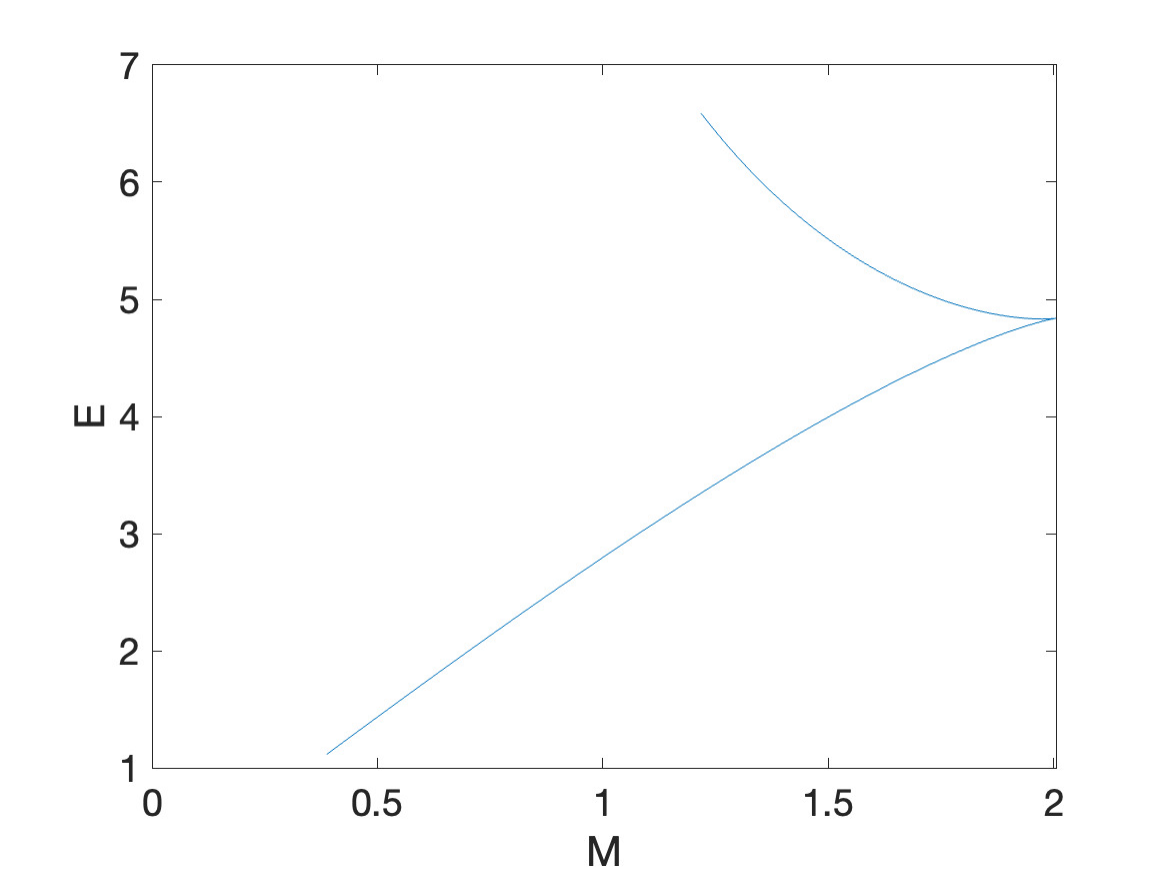}
\subcaption[]{{\footnotesize $E=E(M)$.}}
\end{subfigure}
\caption{\small 2D NLS: $\alpha=1$ (top), $\alpha=2$ (middle) and $\alpha=4$ (bottom). Mass and energy dependence on $b$  (left and middle columns). Dependence of energy as a function of mass, $E = E(M)$ (right column). The maximum value of the mass as a function of $b$ is reached at  $b^\ast \approx -1.0311$ only in the supercritical case ($\alpha=4$). This is the point where the split into two branches occurs in the right bottom column.}  
\label{F:ME-2D}
\end{figure}

We next examine the two dimensional case, and compare the subcritical ($\alpha=1$), critical ($\alpha=2$) and supercritical ($\alpha=4$) cases, see Figure \ref{F:ME-2D}. 
We observe that in all 2D subcritical and critical cases, the mass 
$M(Q_b)$ as a function of $b$ monotonically increases, see (A) and (D) plots in Figure \ref{F:ME-2D}, furthermore, in the critical case it asymptotically approaches the value of the mass of the whole space 
ground state $M(\mathcal{R}) \approx 11.70$. In the supercritical case, similarly to the 1D case, the mass reaches the maximum value at $b=b^\ast$ and then decreases to zero, thus, creating branching in the energy vs. mass graph, see the right bottom plot in Figure \ref{F:ME-2D}. We investigate the stability of different branches in Section \ref{S:branching}.

The critical 2D cubic NLS case was initially studied in \cite{FM2001}. We confirm, in particular, that our computations in the middle row ($\alpha=2$) of Figure \ref{F:ME-2D} for mass and energy are in agreement with those computed in \cite[Fig. 1 and Fig. 2]{FM2001} up to the factor $2\pi$ (which we include in our computations) and corresponding scaling in the energy by a factor $\frac12$.

We make a remark about turning points in 2D. 

\begin{remark}[Turning points for radial ground states in 2D]\label{R:2D-turning}
The radial equation $-Q_b'' - \frac1{r} Q_b' + bQ_b - Q_b^{\alpha+1}=0$
with $Q_b'(0)=0$ and $Q_b(1)=0$ gives
$$
\frac{d}{dr}\left[- \frac12 (Q_b'(r))^2 +\frac{b}2 (Q_b(r))^2
-  \frac{1}{\alpha+2}(Q_b(r))^{\alpha+2}\right]
=\frac1r(Q_b'(r))^2.
$$
Comparing with the 1D equation \eqref{E:1D-integralQ}, we notice that the left-hand sides match (with $x$ replaced by $r$) and the corresponding square bracket is constant in 1D, as the right-hand side in 1D is zero, while in 2D we have an extra radial term, which is not straightforward to deal with when integrating: it gives $\int_0^r(Q'_b(\rho))^2/\rho \, d\rho$, which depends on the unknown profile, and thus, a similar integral quadrature used in Proposition \ref{P:1D-b*} is not available in 2D. 
Nevertheless, we are still able to obtain an approximate value $b^*$ numerically. For instance, in 2D quintic ($\alpha=4$) NLS, we obtained $b^* \approx -1.0311$, see the maximum of mass in Figure \ref{F:ME-2D} (G).
\end{remark}

\subsubsection{Branching}\label{S:branching}

In the supercritical cases, the value $b^*$ as described above can be obtained either from Prop. \ref{P:1D-b*} in 1D or from the maxima of the computed curve $b \mapsto M(Q_b)$. 
We use these turning points as a threshold to select values of $b$ on both sides of the transition and then evolve perturbations of ground state $Q_b$ of the form 
$$
u_0 = A \, Q_b \quad \mbox{with} \quad A \sim 1.
$$ 
In the simulations shown below we take $A=0.99$ or $A=1.01$ (we also ran simulations for other amplitudes, e.g., $A = \{0.95, 0.98, 1.02, 1.05\}$, and observed the same qualitative types of dynamics; for brevity we omit those figures).  
\smallskip

\begin{figure}[htb!]
\begin{subfigure}{.24\textwidth}
\includegraphics[width=1\linewidth,height=0.85\linewidth]{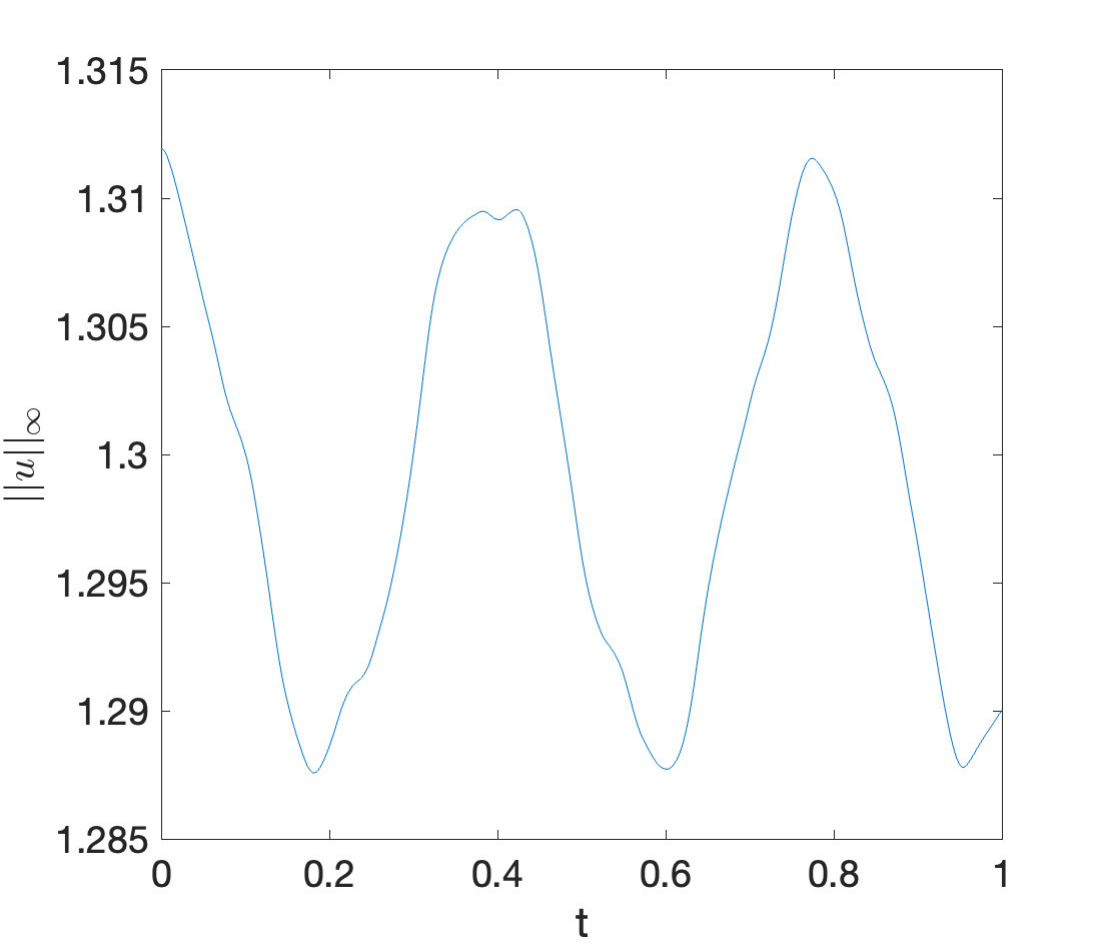} 
\subcaption[]{{\footnotesize $b=0$, $A=0.99$.}}
\end{subfigure}
\begin{subfigure}{.24\textwidth}
\includegraphics[width=1\linewidth,height=0.85\linewidth]{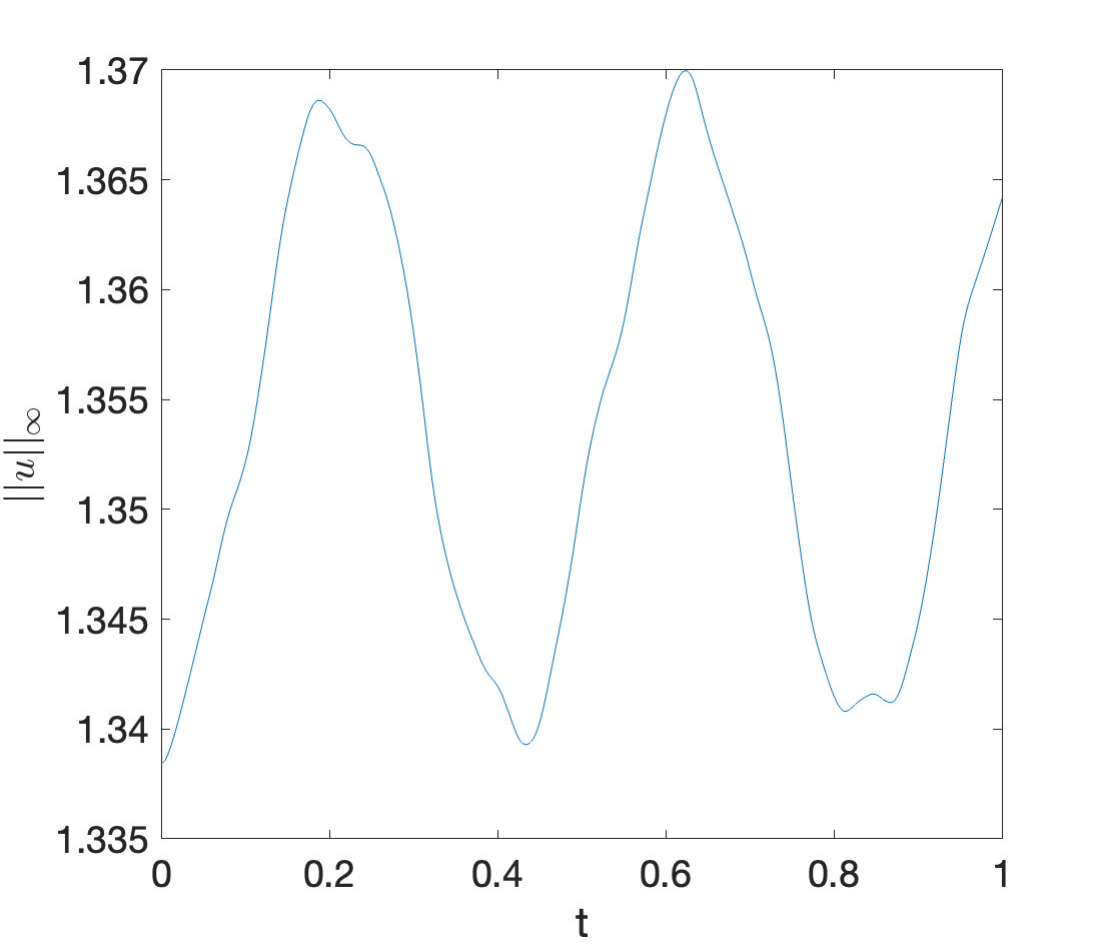}
\subcaption[]{{\footnotesize $b=0$, $A=1.01$.}}
\end{subfigure}
\begin{subfigure}{.24\textwidth}
\includegraphics[width=1\linewidth,height=0.85\linewidth]{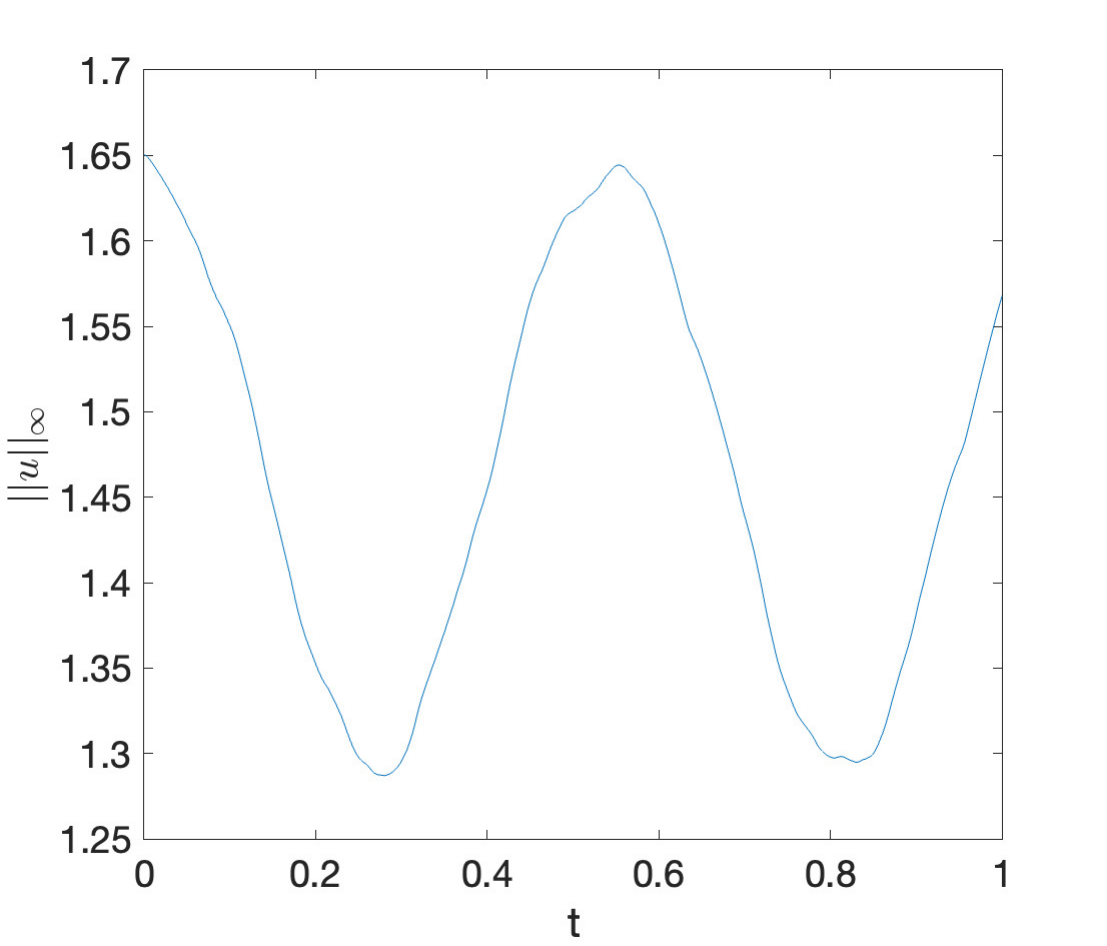} 
\subcaption[]{{\footnotesize $b=5$, $A=0.99$.}}
\end{subfigure}
\begin{subfigure}{.24\textwidth}
\includegraphics[width=1\linewidth,height=0.85\linewidth]{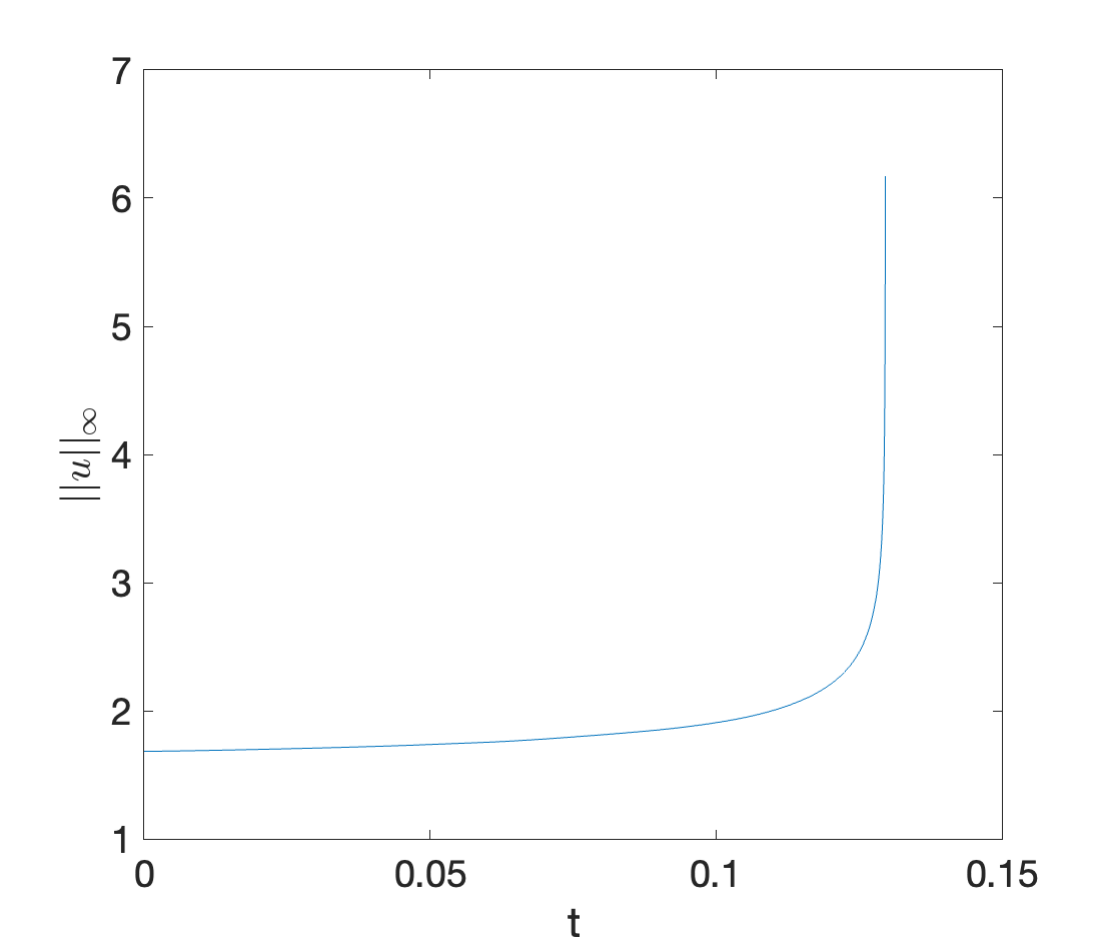}
\subcaption[]{{\footnotesize $b=5$, $A=1.01$.}}
\end{subfigure}
\caption{\footnotesize 
1D septic NLS, $\alpha=6$. (A) \& (B): stable branch, $b=0 < b^\ast \approx 3.2945$, $u_0=AQ_{0}$. 
(C) \& (D): unstable branch, $b = 5 > b^\ast \approx 3.2945$, $u_0 = AQ_{5}$. 
The values of $A$ as indicated.}
\label{F:1D-septic}
\end{figure}

$\bullet$ \underline{1D, $\alpha=6$.} ~~ In Figure \ref{F:1D-septic} in plots (A) and (B) we show the perturbations of the stable branch in 1D septic NLS, $\alpha=6$, with $b=0 < b^\ast \approx 3.3$. Note that both perturbations have tiny oscillations (in the second or third digit) around the $Q_b$ itself. On the other hand, in plots (C) and (D) one can see perturbations of different kind (there we take $b=5 > b^\ast$): solutions with $A=1.01$ blow up in finite time, while the initial data $u_0 = 0.99\, Q_b$ oscillates more drastically (e.g., in (C) between 1.25 and 1.65), indicating that the oscillations are not around one state but rather between two different states. (See Proposition \ref{P:blowup} and Remark \ref{R:1D-?} for further explanations on blow-up.)

\begin{figure}[h!]
\begin{subfigure}{.24\textwidth}
\includegraphics[width=1\linewidth,height=0.85\linewidth]{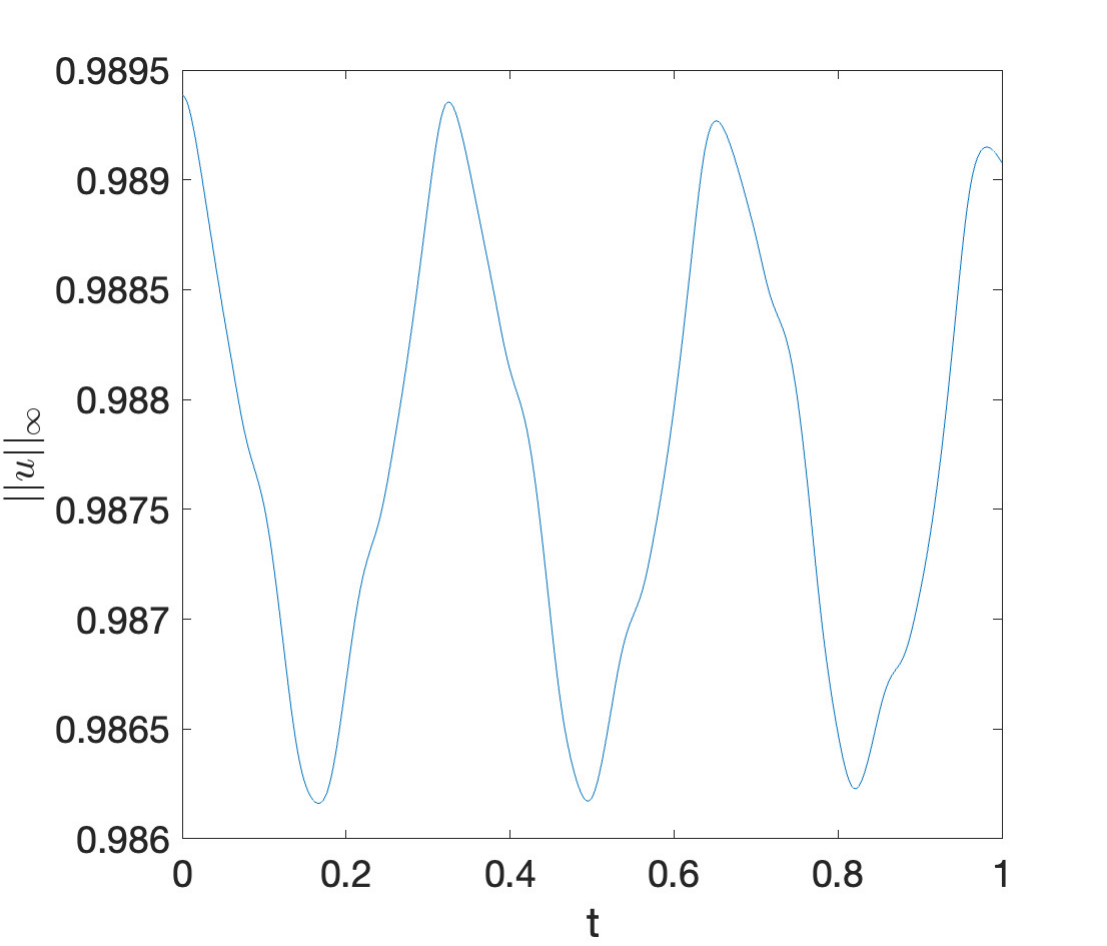} 
\subcaption[]{{\footnotesize $b=-2$, $A=0.99$.}}
\end{subfigure}
\begin{subfigure}{.24\textwidth}
\includegraphics[width=1\linewidth,height=0.85\linewidth]{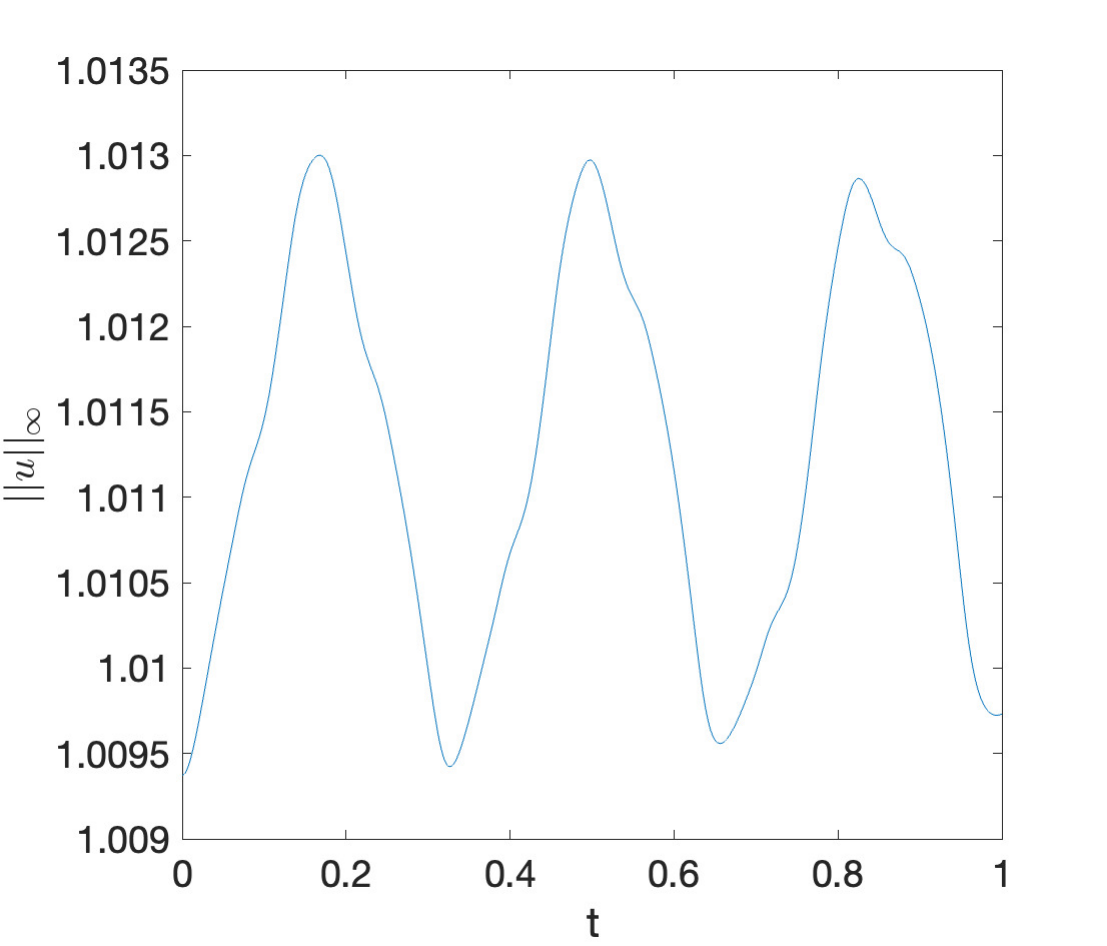}
\subcaption[]{{\footnotesize $b=-2$, $A=1.01$.}}
\end{subfigure}
\begin{subfigure}{.24\textwidth}
\includegraphics[width=1\linewidth,height=0.85\linewidth]{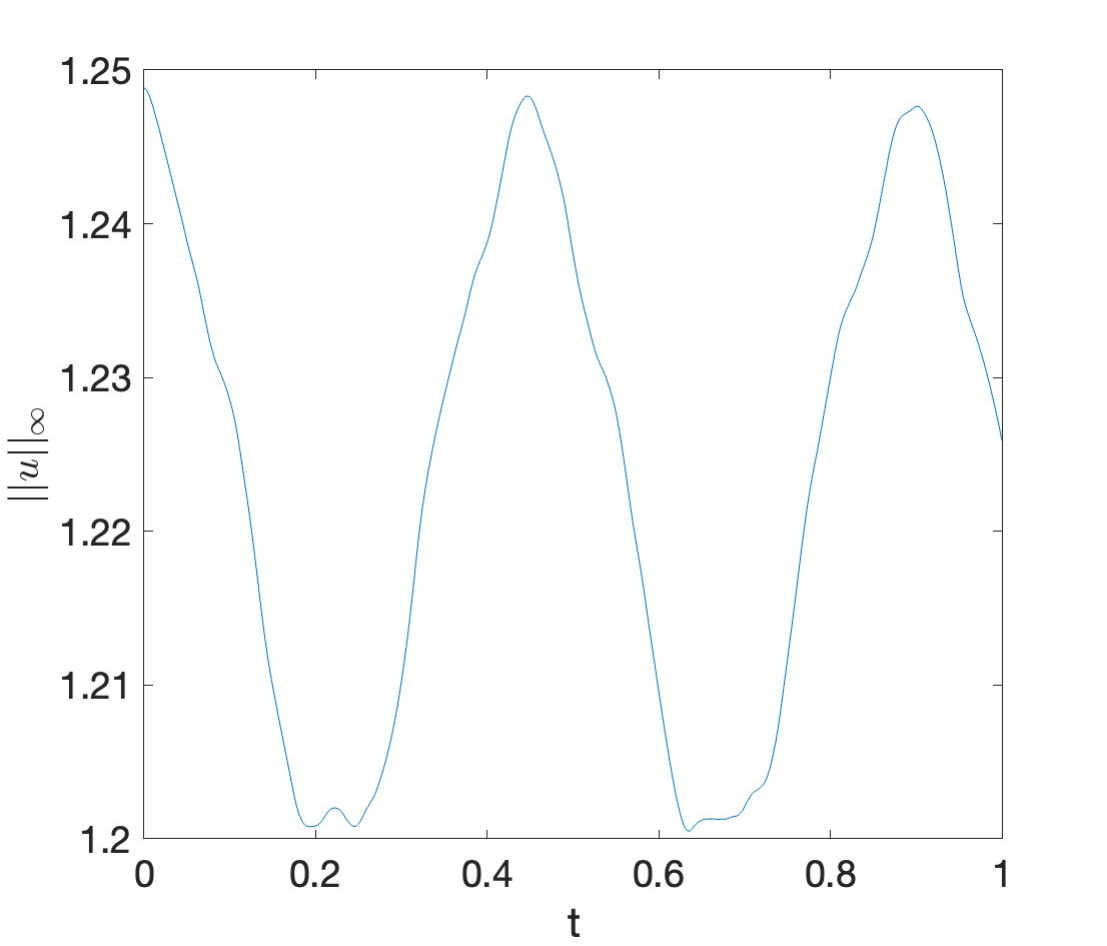} 
\subcaption[]{{\footnotesize $b=0$, $A=0.99$.}}
  \end{subfigure}
\begin{subfigure}{.24\textwidth}
\includegraphics[width=1\linewidth,height=0.85\linewidth]{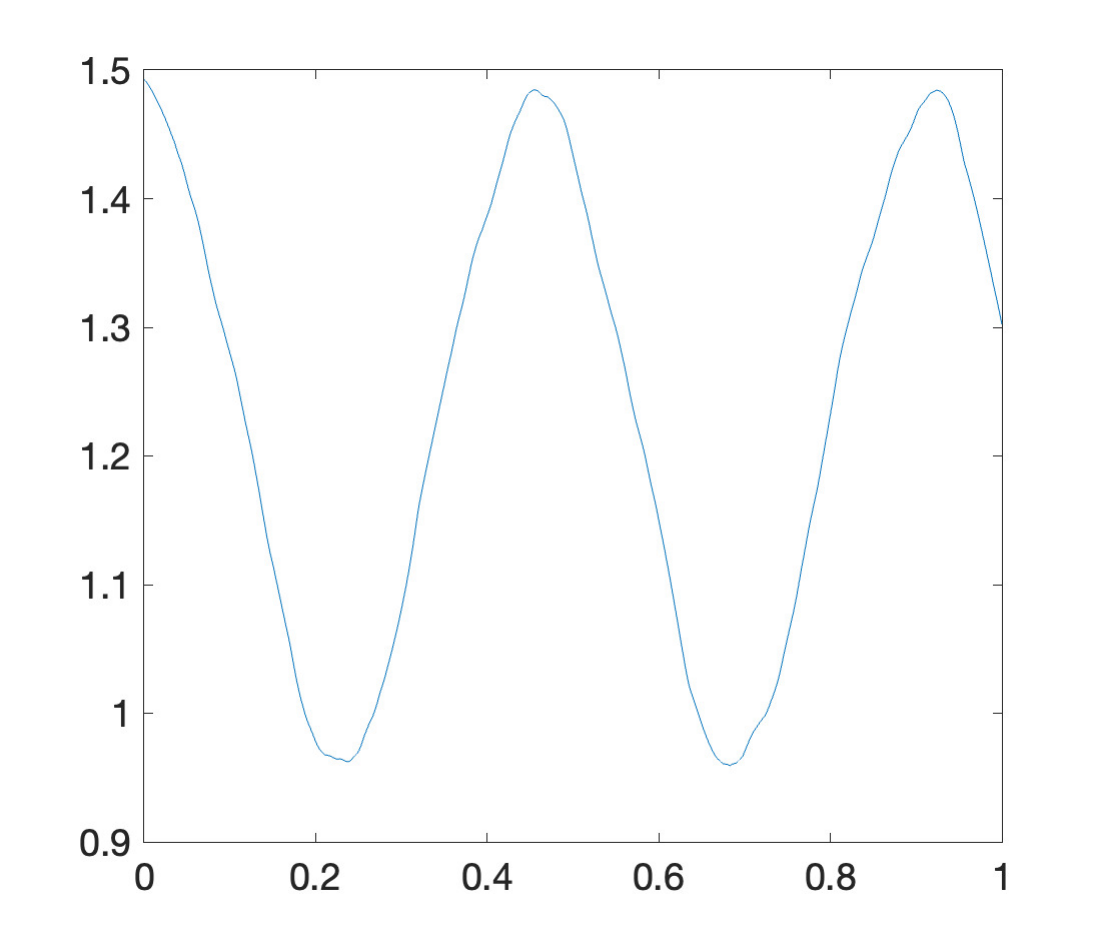}
\subcaption[]{{\footnotesize $b=5$, $A=0.99$.}}
\end{subfigure}
\caption{\footnotesize 
1D nonic NLS, $\alpha=8$. (A) \& (B): stable branch, $b=-2 < b^\ast \approx 1.06$, $u_0=AQ_{-2}$, $A = 0.99$ or $1.01$. 
(C): stable branch, $b = 0 < b^\ast \approx 1.06$, $u_0 = 0.99 Q_{0}$; (D): unstable branch, $b = 5 > b^\ast \approx 1.06$, $u_0 = 0.99 Q_{5}$.  }
\label{F:1D-nonic}
\end{figure}

$\bullet$ \underline{1D, $\alpha=8$.} ~~Solution dynamics for another supercritical case with the nonic nonlinearity $\alpha = 8$ is examined in Figure \ref{F:1D-nonic}, where the behavior seems to be analogous to the septic case.
There, we compare the 1D nonic NLS perturbations at $b=-2, 0, 5$, on either side of the turning point (for brevity, while we observe that for $b>b^*$ and $A>1$ solutions blow up, we omitted the plots). 

\begin{figure}[htb!]
\centering 
\begin{subfigure}[b]{0.53\textwidth}
\includegraphics[width=\textwidth,height=.85\textwidth]{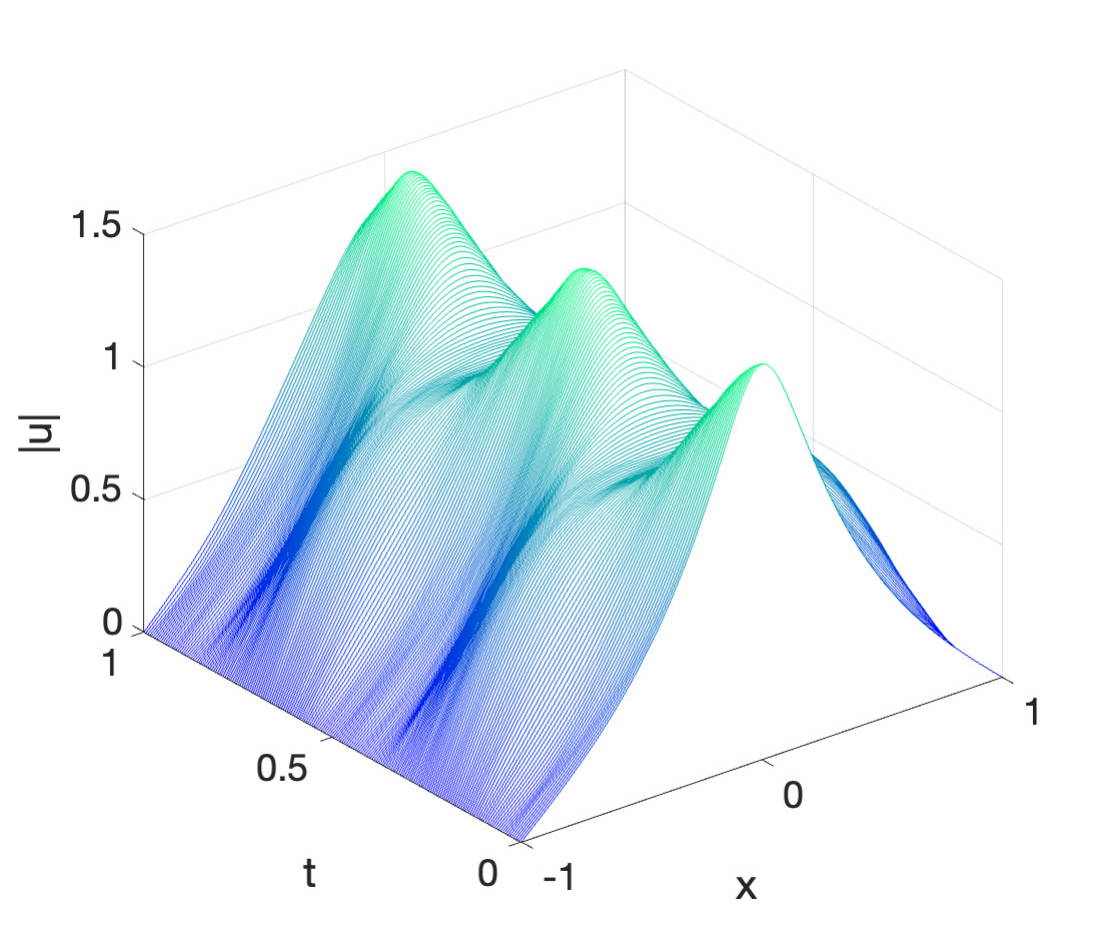} 
\caption{Oscillations between two coherent profiles.}
\end{subfigure}
\begin{subfigure}[b]{0.46\textwidth}
\centering  
\begin{subfigure}[b]{\textwidth}
\centering
\includegraphics[width=.8\textwidth,height=.4\textwidth]{NLS1Dp9_Qb5_099inf-eps-converted-to.pdf} 
\caption{\footnotesize Time dependence of the $L^\infty$ norm.}
\end{subfigure}
\vspace{1mm} 
\begin{subfigure}[b]{0.49\textwidth}
\centering
\includegraphics[width=\textwidth,height=.85\textwidth]{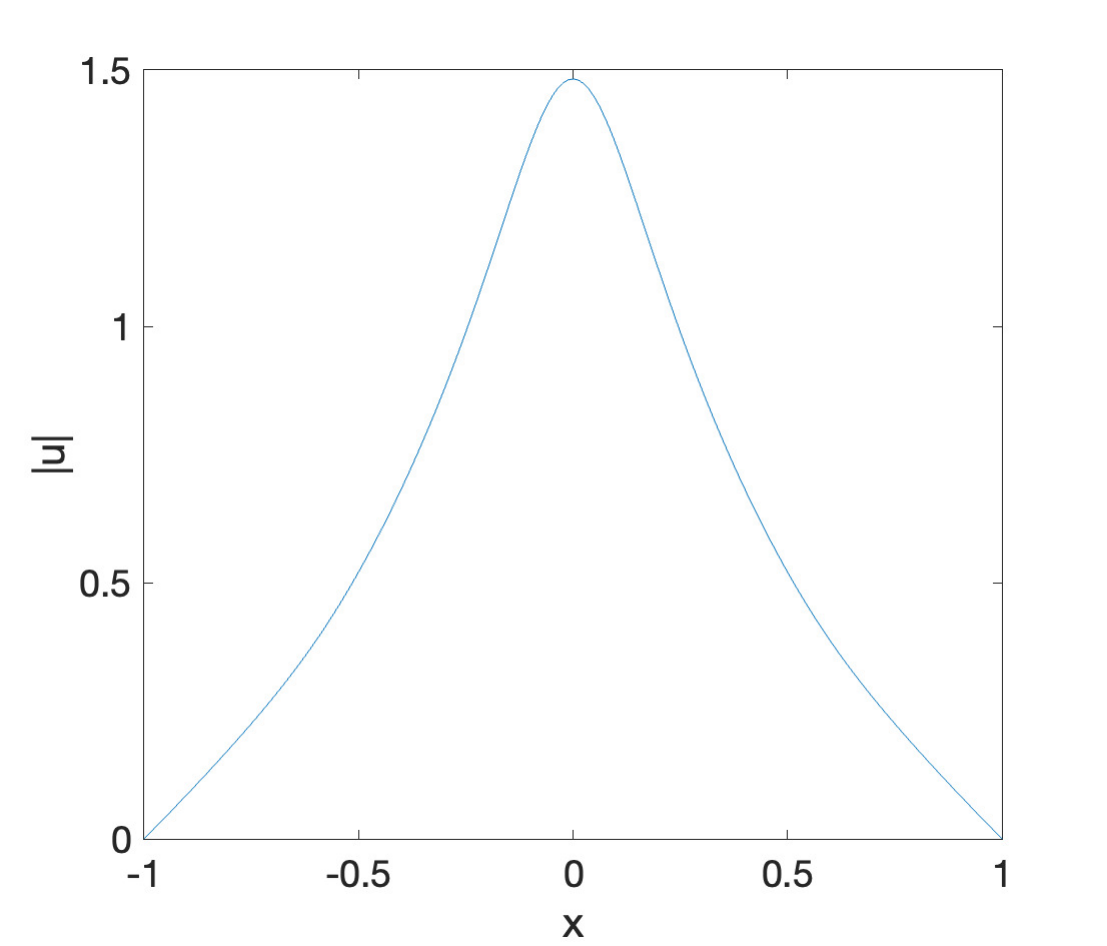} 
\caption{\footnotesize Profile 1 at $t=0.45$.}
\end{subfigure}%
\hfill 
\begin{subfigure}[b]{0.49\textwidth}
\centering
\includegraphics[width=\textwidth,height=.85\textwidth]{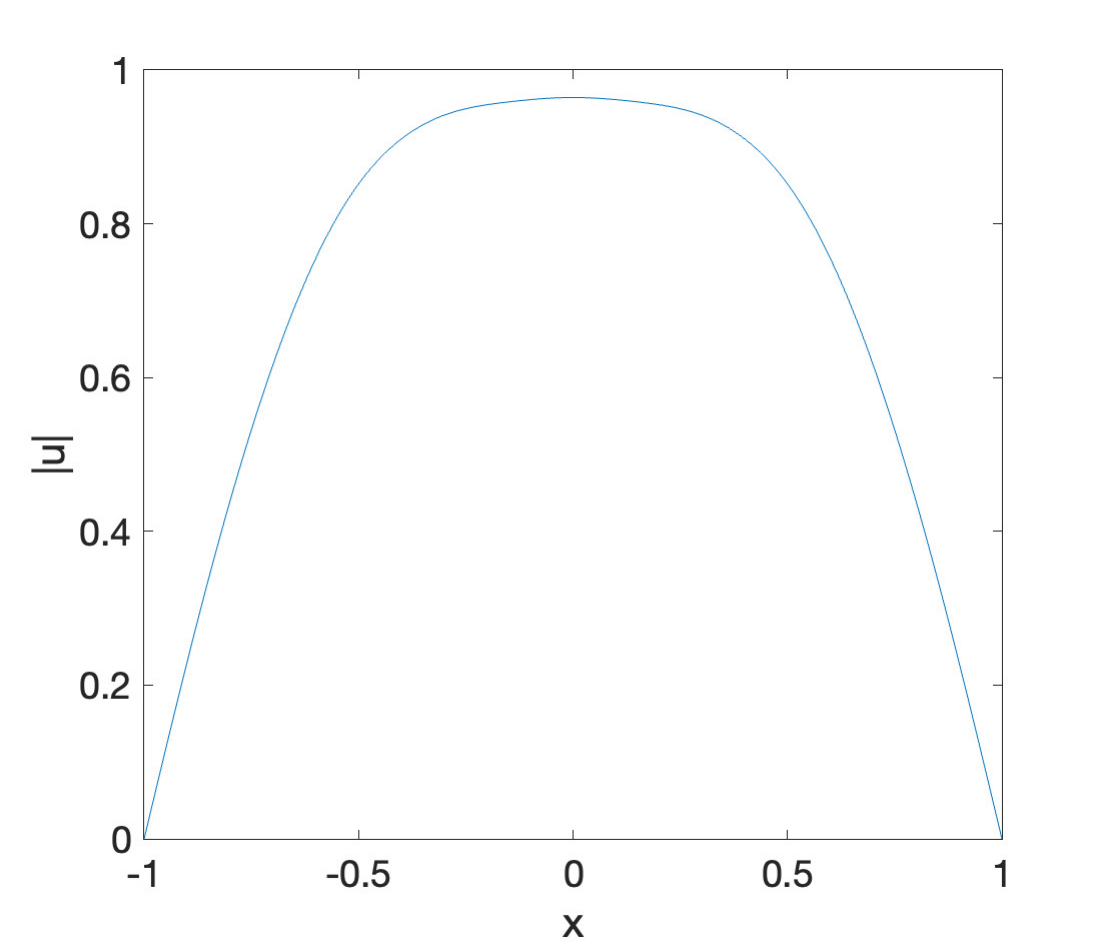} 
\caption{\footnotesize Profile 2 at $t=0.23$.}
\end{subfigure}
\end{subfigure}
\caption{\small 1D nonic NLS, $\alpha=8$: recurrent oscillations for a perturbation of the unstable branch of $Q_b$ with $b = 5 > b^\ast$, $u_0=0.99Q_{5}$. The profiles shown are two coherent profiles visited during the oscillation.}
\label{F:1D-nonic-water}
\end{figure}


Further examining small-amplitude perturbations, we find that the oscillation at least between two states (or profiles) is, indeed, the case, see Figure \ref{F:1D-nonic-water} with three-dimensional visualization of recurrent oscillations between two coherent profiles. 
In that figure 
we keep studying the 1D nonic NLS, $\alpha=8$, and look at the unstable branch with $b=5 > b^\ast \approx 1.06$. Plot (A) shows the solution with the initial condition $u_0=0.99Q_{5}$, which exhibits pronounced recurrent oscillation between two states. The $L^\infty$ norm is plotted in (B) and profiles of the two states are given in (C) and (D), where the profiles are taken at the maximum and minimum locations (and seem to continue to oscillate between these states). 
 Note that since the initial condition has mass less than the ground state $Q_b$, the solution starts to disperse, however, the Dirichlet boundary condition prevents escape from the ball and reflects the solution so that it concentrates into another profile (D).  The height of this state (D) is below that of the stable branch ground state, so its height starts increasing but passes it and reaches back the profile (C), continuing this recurrent behavior. 

$\bullet$ \underline{2D, $\alpha=2$.} ~~
To understand better the 2D case, we start with the critical case first, and 
consider perturbations of the ground state $Q_b$ for several values of $b$, together 
with solutions that have the initial data above the mass of the whole-space ground state $M(\mathcal{R})$, which we schematically indicated in the left plot of Figure \ref{F:cloud}. 

\begin{figure}[htb!]
\includegraphics[width=0.32\textwidth,height=0.21\textwidth]{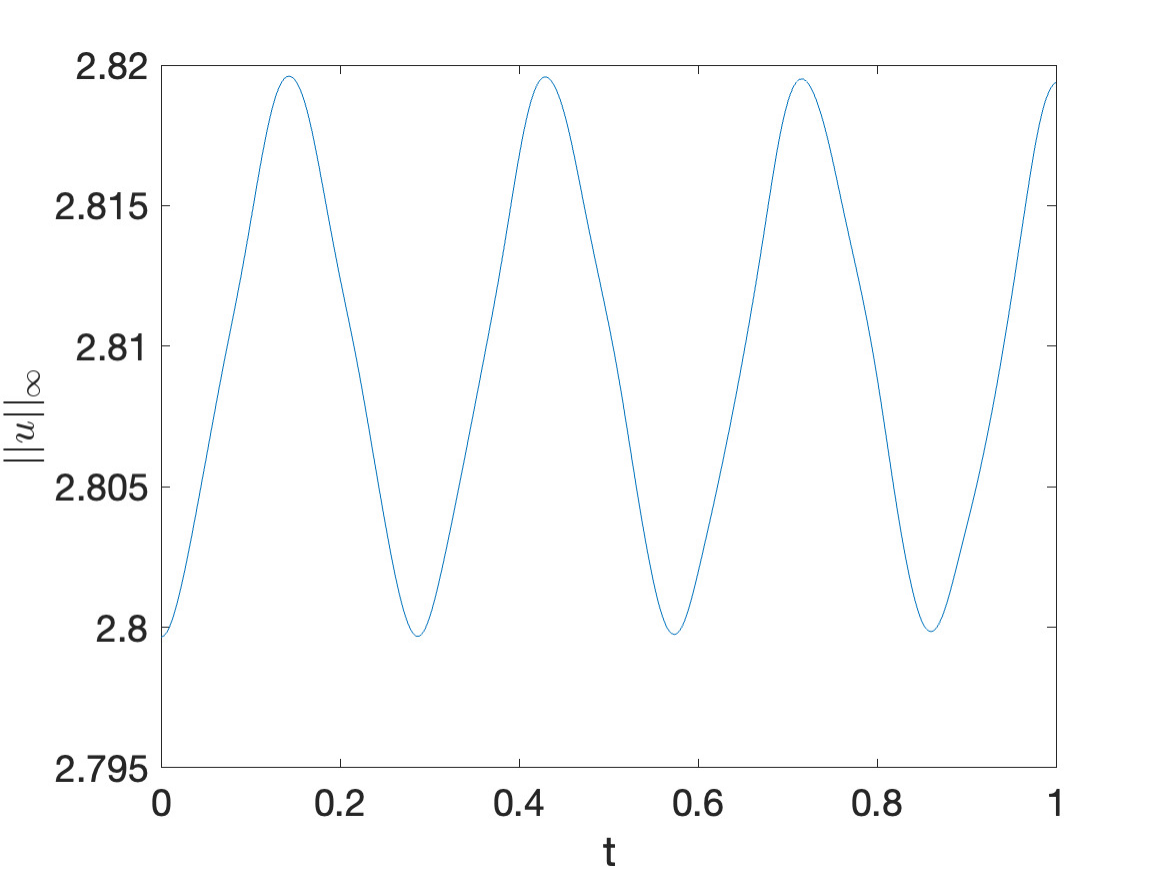} 
\includegraphics[width=0.32\textwidth,height=0.21\textwidth]{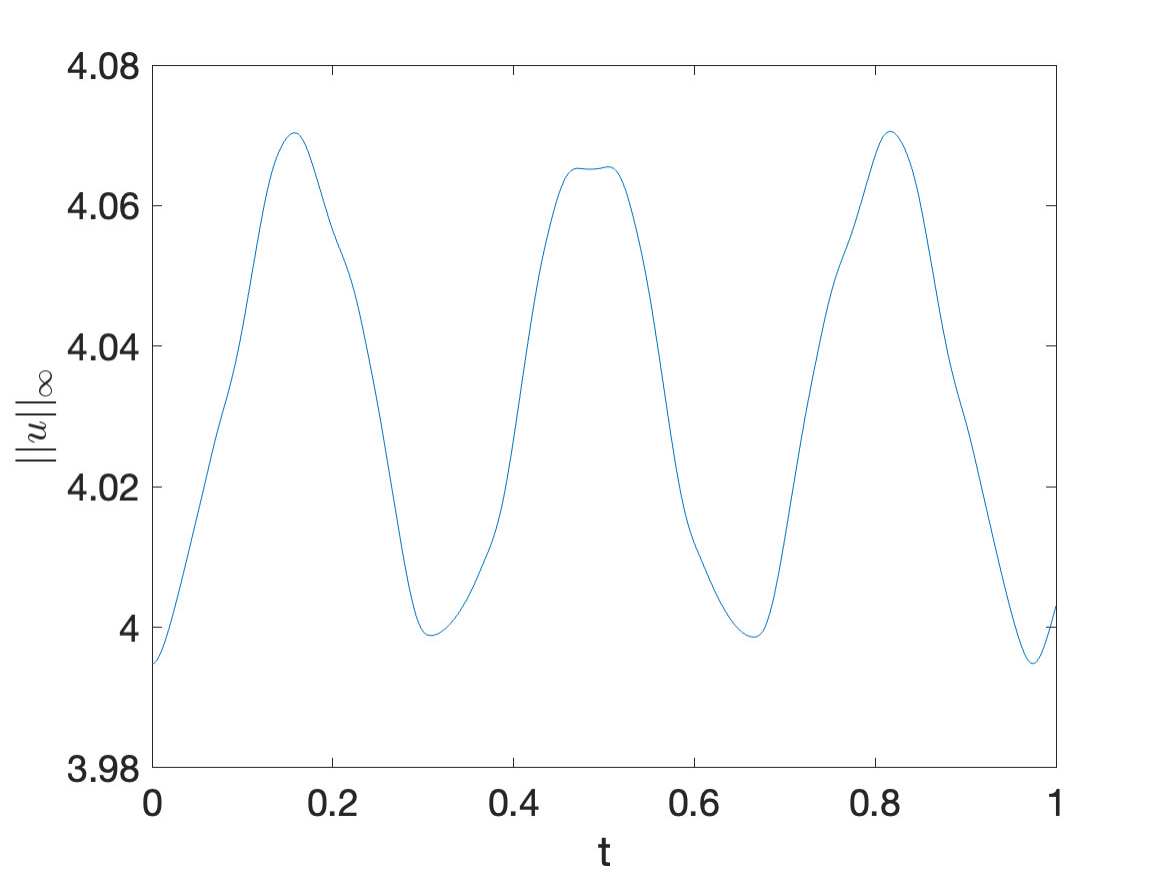}
\includegraphics[width=0.32\textwidth,height=0.21\textwidth]{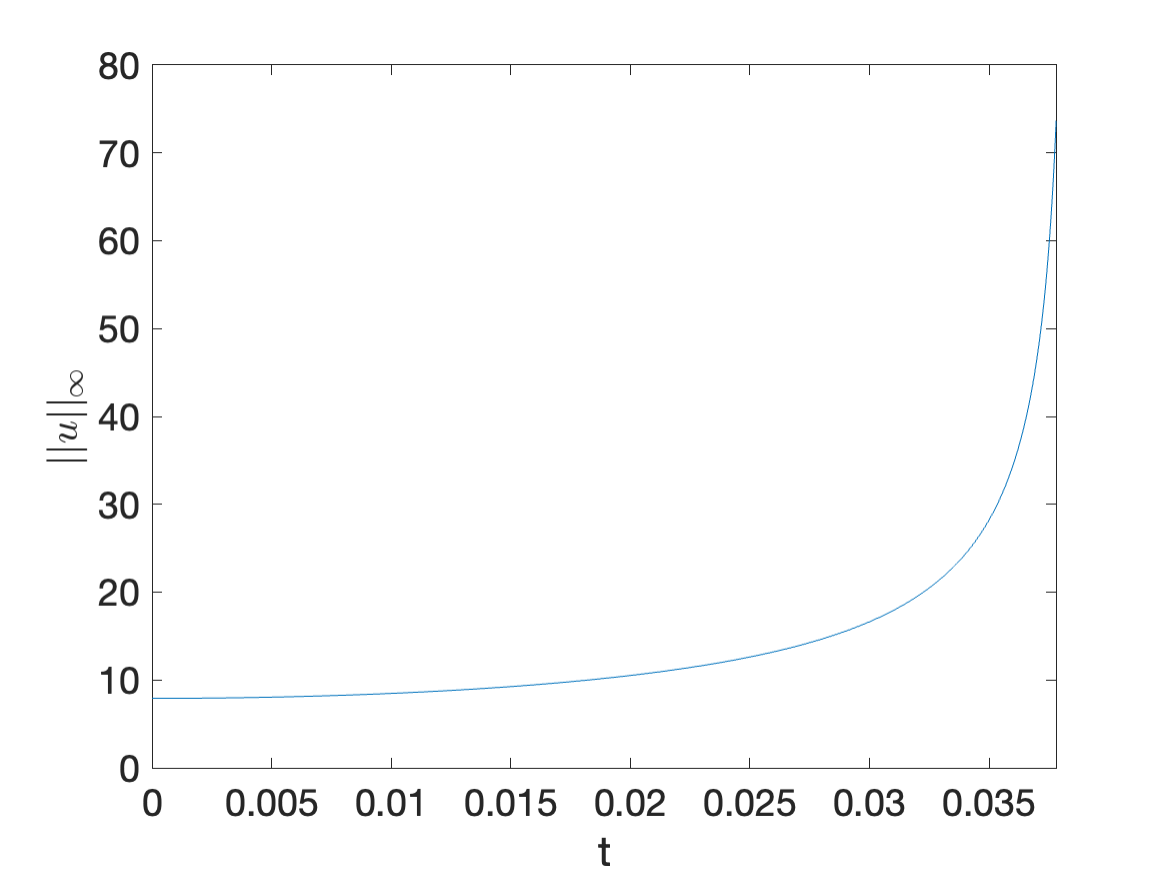}
\caption{\small 2D cubic NLS: time evolution of the $L^\infty$ norm for $u_0 = 1.01 Q_b$ with $b=-2.5$ (left), $b=1$ (middle), and for $u_0 = 2 Q_1$ (right), where $M(2 Q_1) > M(\mathcal{R})$.} 
\label{F:2D-cubic-threshold}
\end{figure}

In Figure \ref{F:2D-cubic-threshold} we show that in the 2D cubic NLS small perturbations of $Q_b$ are stable: for $A=1.01$, the observed oscillations in the $L^\infty$ norm occur only in the second or third digit. We also notice that, for smaller $b$, the stable perturbation interval in the amplitude $A$ is visibly wider, while for larger $b$, as $M(Q_b)$ approaches $M(\mathcal R)$, the admissible upward perturbation becomes smaller. 
Above $M(\mathcal{R})$ the tested amplified data blow up, for example, the datum $u_0=2Q_1$ ($b=1$) has mass above $M(\mathcal R)$, its blow-up dynamics is shown in the right plot of Figure \ref{F:2D-cubic-threshold} (see Remark \ref{R:2Q1} for the proof).  

\begin{figure}[htb!]
\begin{subfigure}{.24\textwidth}
\includegraphics[width=1\linewidth,height=0.81\linewidth]{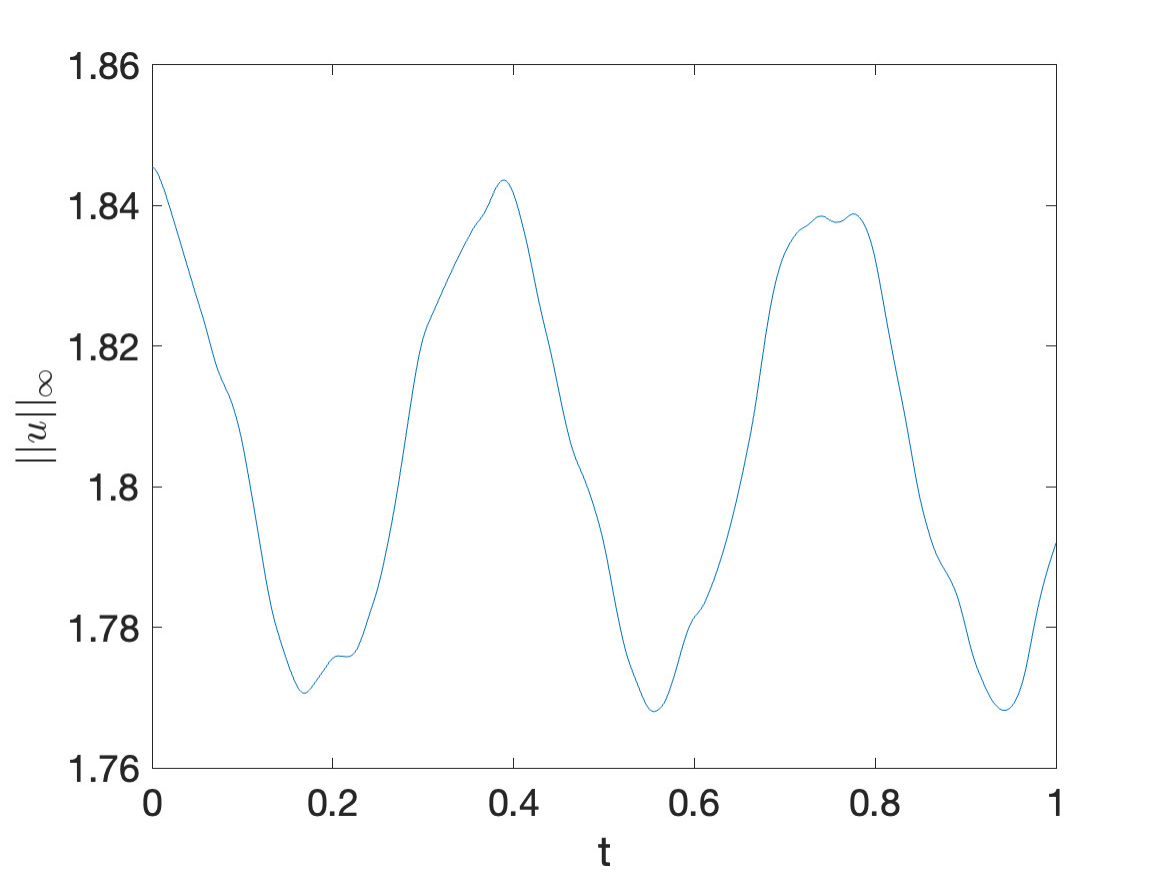}
\subcaption[]{{\footnotesize $b=-2.5$, $A=0.99$.}}
\end{subfigure}
\begin{subfigure}{.24\textwidth}
\includegraphics[width=1\linewidth,height=0.81\linewidth]{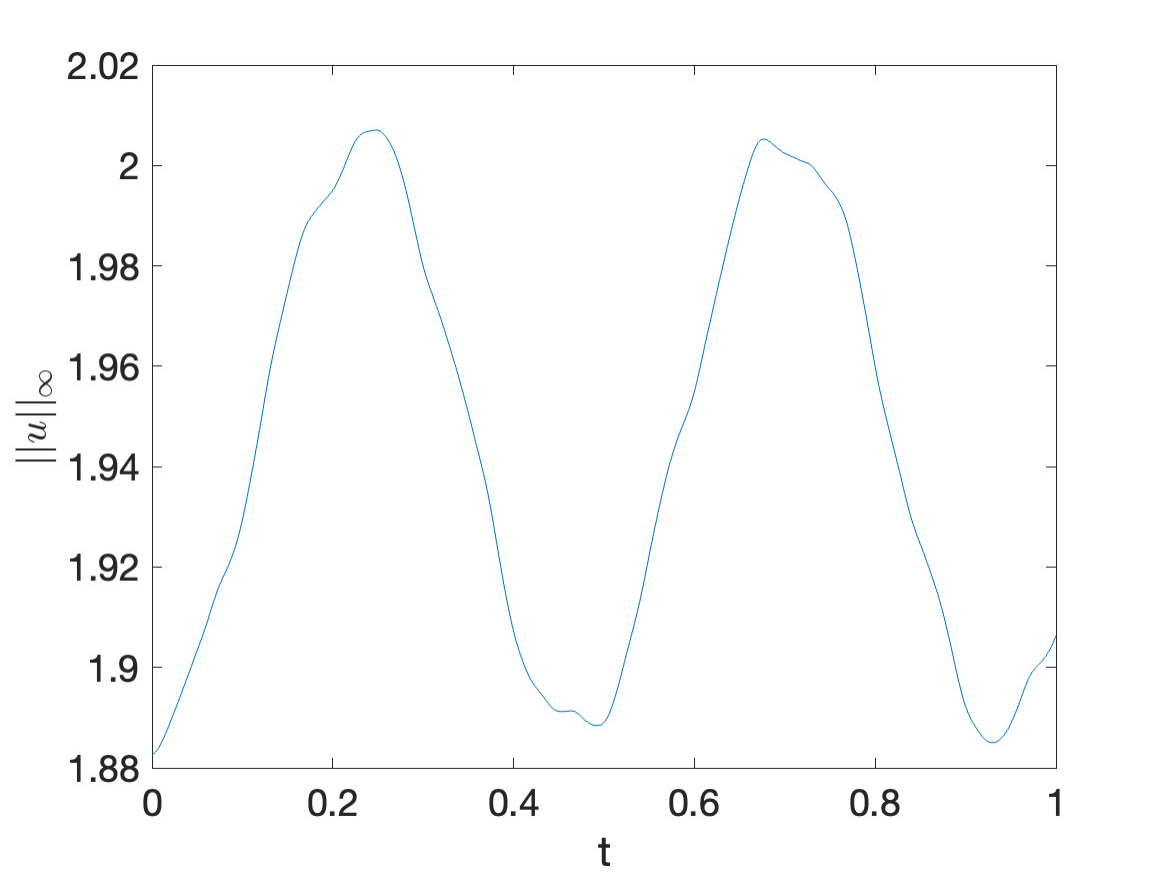}
\subcaption[]{{\footnotesize $b=-2.5$, $A=1.01$.}}
\end{subfigure}
\begin{subfigure}{.24\textwidth}
\includegraphics[width=1\linewidth,height=0.81\linewidth]{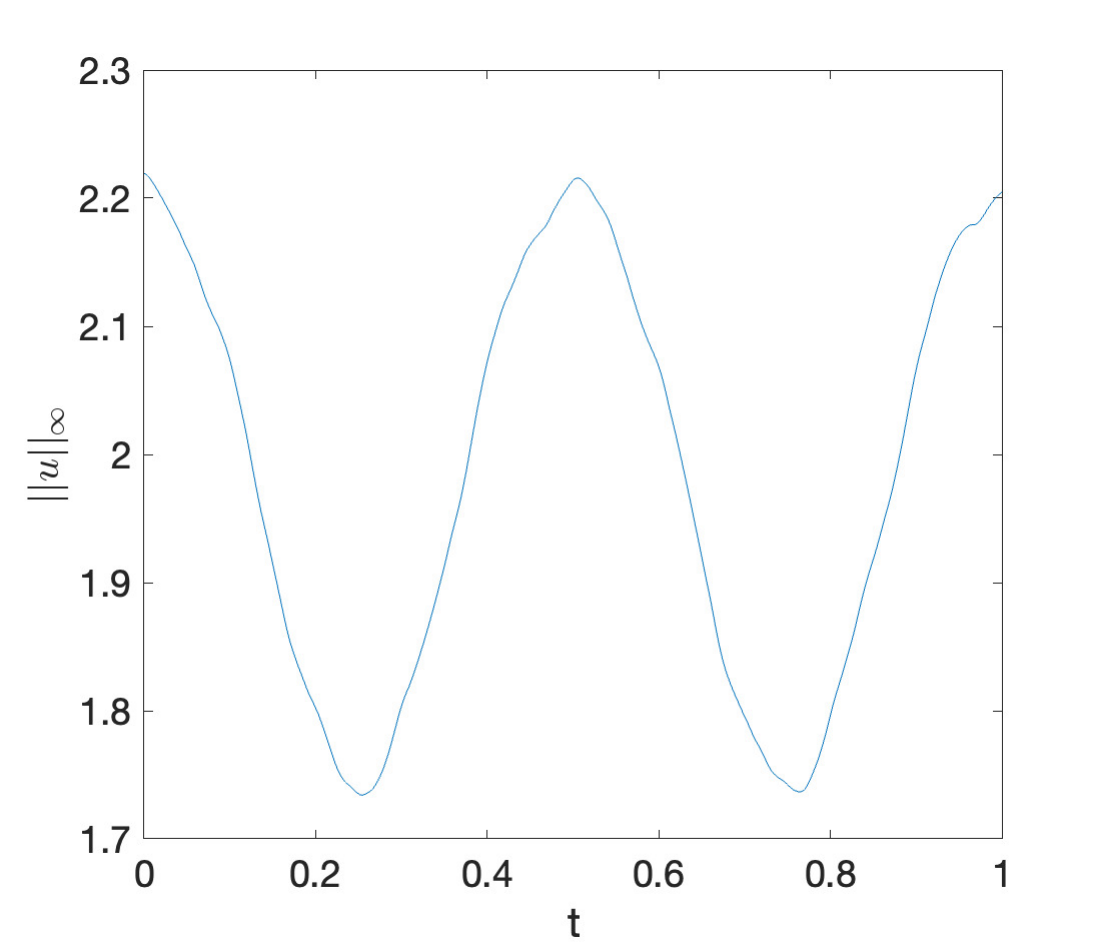} 
\subcaption[]{{\footnotesize $b=-0.5$, $A=0.99$.}}
\end{subfigure}
\begin{subfigure}{.24\textwidth}
\includegraphics[width=1\linewidth,height=0.81\linewidth]{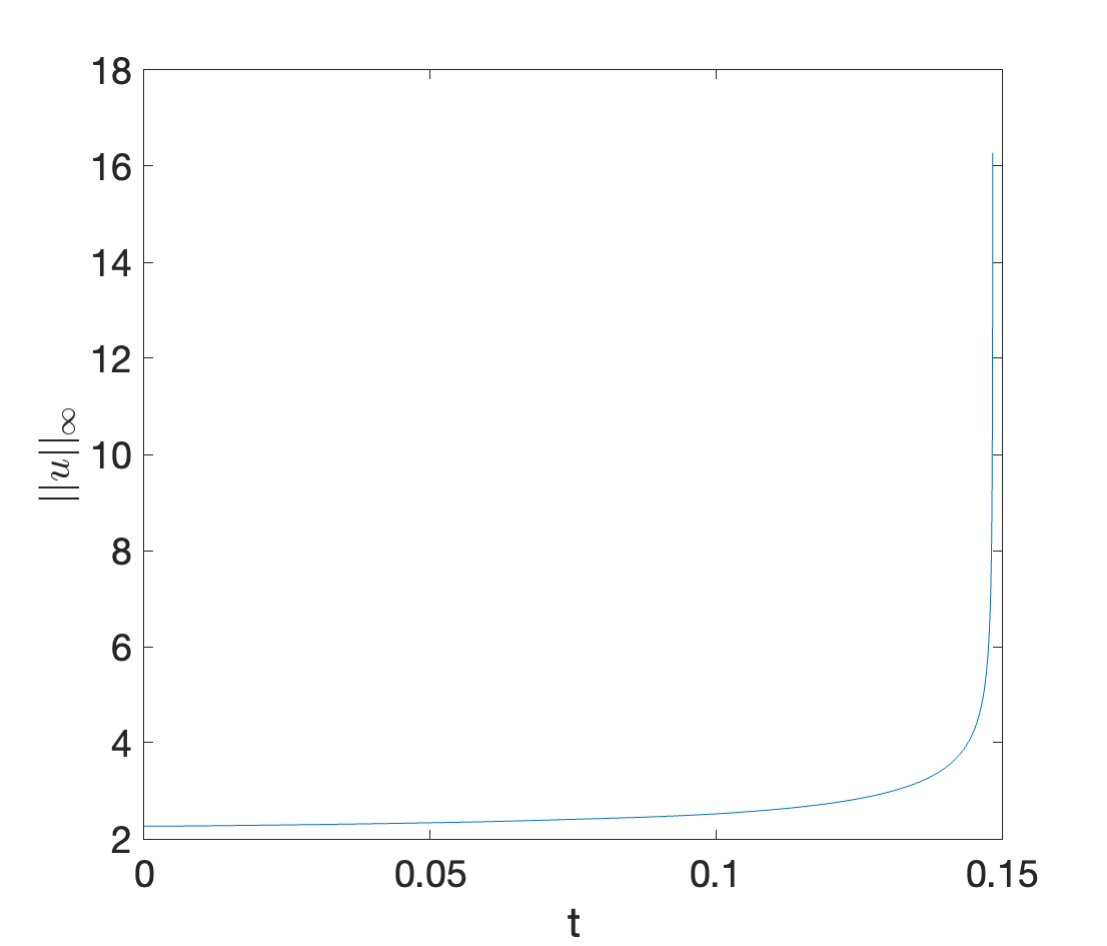}
\subcaption[]{{\footnotesize $b=-0.5$, $A=1.01$.}}
\end{subfigure}
\caption{\footnotesize 
2D quintic NLS, $\alpha=4$. (A) \& (B): stable branch, $b=-2.5 < b^\ast \approx -1.03$, $u_0=AQ_{-2.5}$. (C) \& (D): unstable branch, $b = -0.5 > b^\ast$, $u_0 = AQ_{-0.5}$. }
\label{F:2Dquintic}
\end{figure}

$\bullet$ \underline{2D, $\alpha=4$.} ~~ Considering the two dimensional supercritical cases, we show in Figure \ref{F:2Dquintic} perturbations of the ground states $Q_b$ in the 2D quintic NLS, $\alpha=4$. 
Perturbations of the stable branch in plots (A) and (B) remain close to the ground state itself, with changes in the second digit for both $A=0.99$ and $A=1.01$ (here, $b=-2.5 < b^\ast \approx -1.03$). In contrast, for the unstable branch at $b=-0.5>b^\ast$, the perturbation with $A=0.99$ shows recurrent oscillations between two different states, while the perturbation with $A=1.01$ blows up in finite time, see plots (C) and (D) (we examine blow-up below in Proposition \ref{P:blowup}).

\begin{figure}[htb!]
\centering     
\begin{subfigure}[b]{0.53\textwidth}
\centering
\includegraphics[width=\textwidth,height=.8\textwidth]{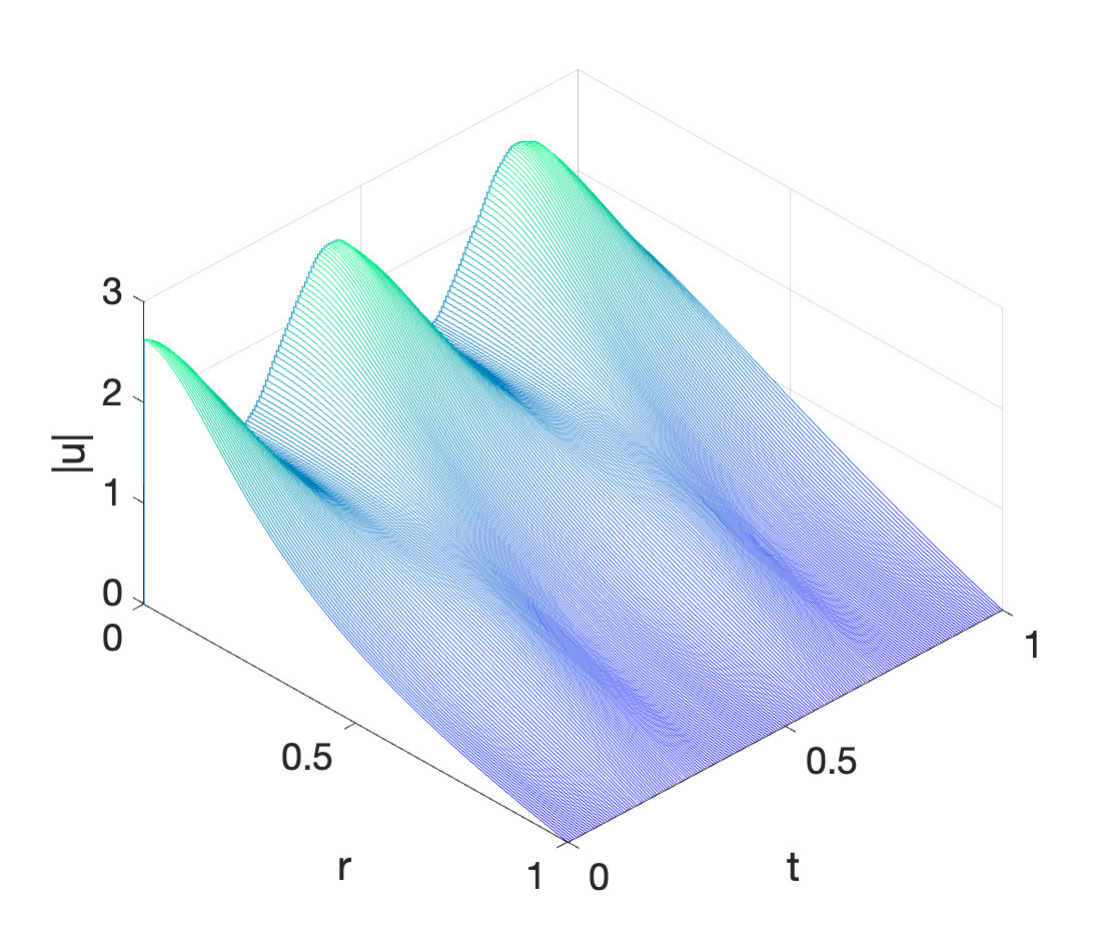} 
\caption{Oscillations between two coherent profiles.}
\end{subfigure}
\begin{subfigure}[b]{0.46\textwidth}
\centering  
\begin{subfigure}[b]{\textwidth}
\centering
\includegraphics[width=.8\textwidth,height=.4\textwidth]{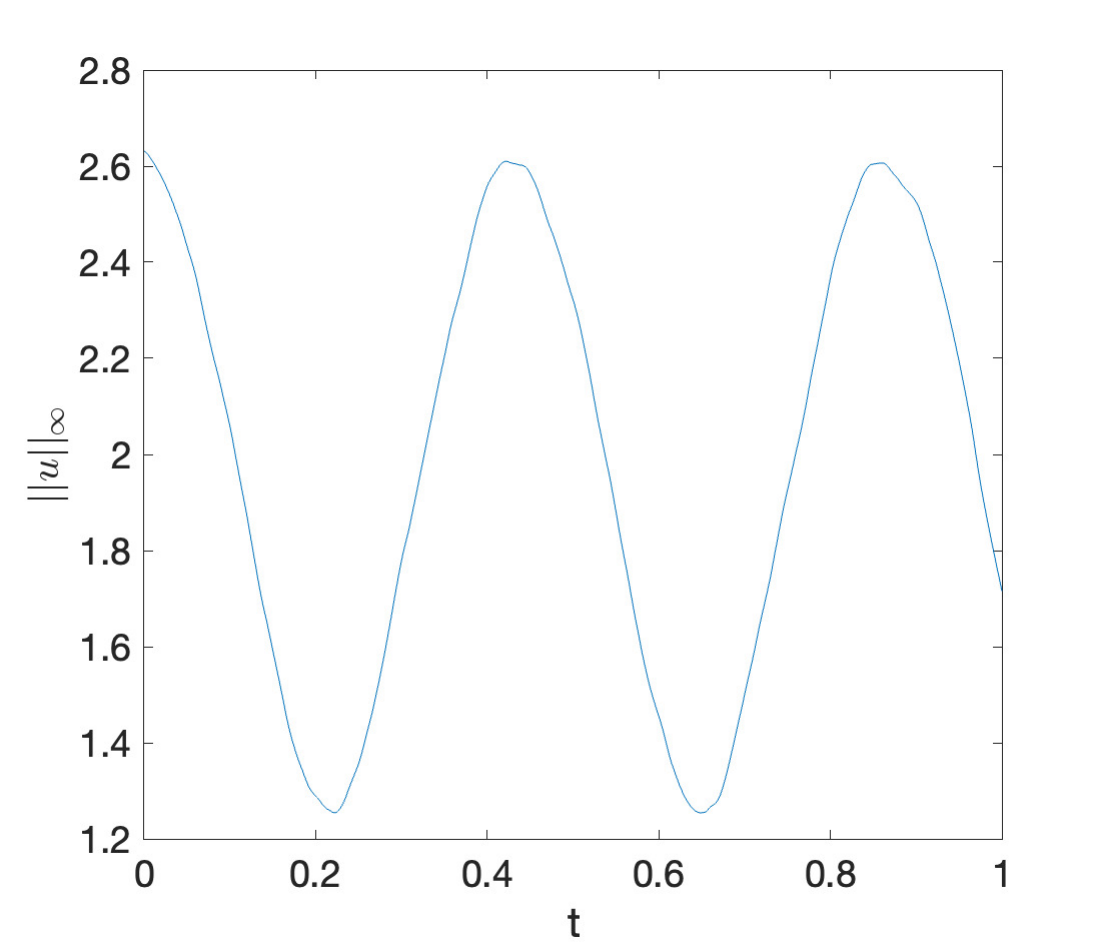} 
\caption{\footnotesize Time dependence of the $L^\infty$ norm.}
\end{subfigure}
\vspace{1mm} 
\begin{subfigure}[b]{0.49\textwidth}
\centering
\includegraphics[width=\textwidth,height=.8\textwidth]{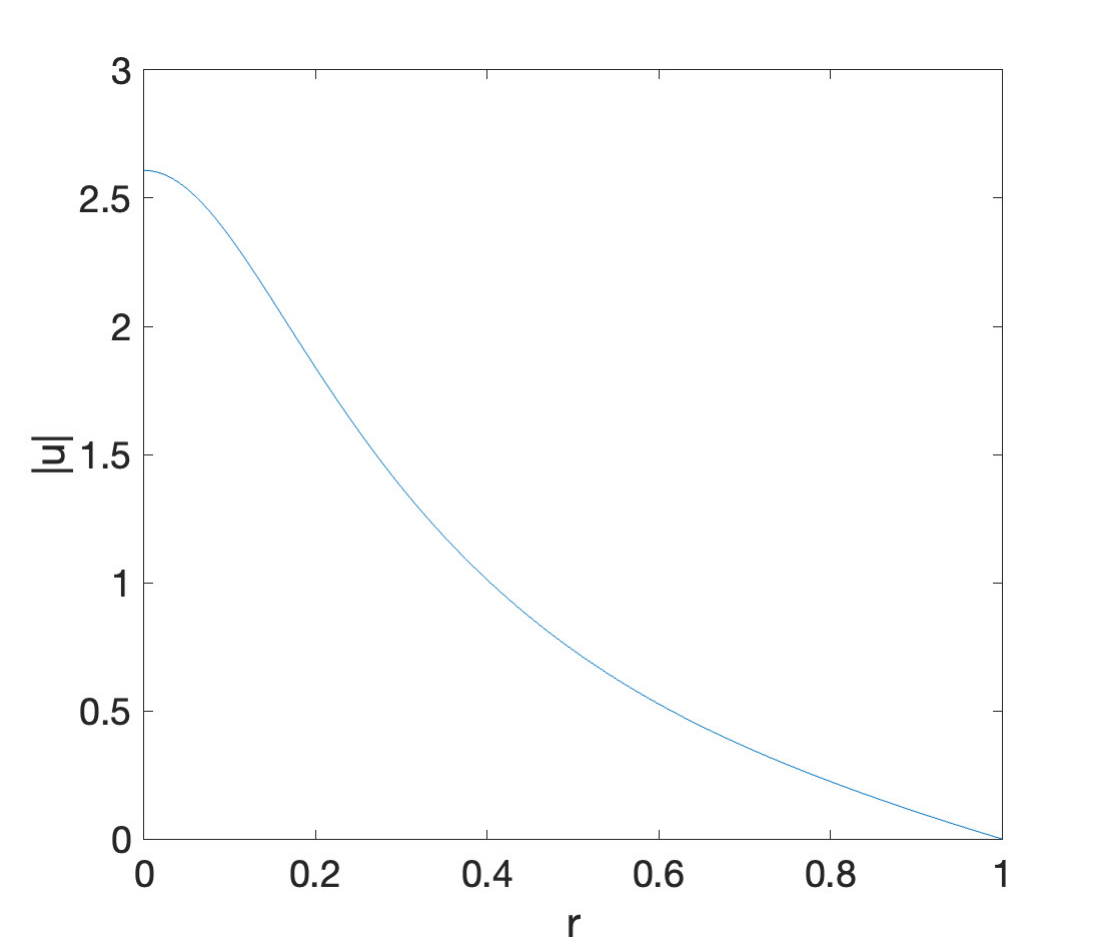} 
\caption{\footnotesize Profile 1 at $t=0.42$.}
\end{subfigure}
\hfill
\begin{subfigure}[b]{0.49\textwidth}
\centering
\includegraphics[width=\textwidth,height=.8\textwidth]{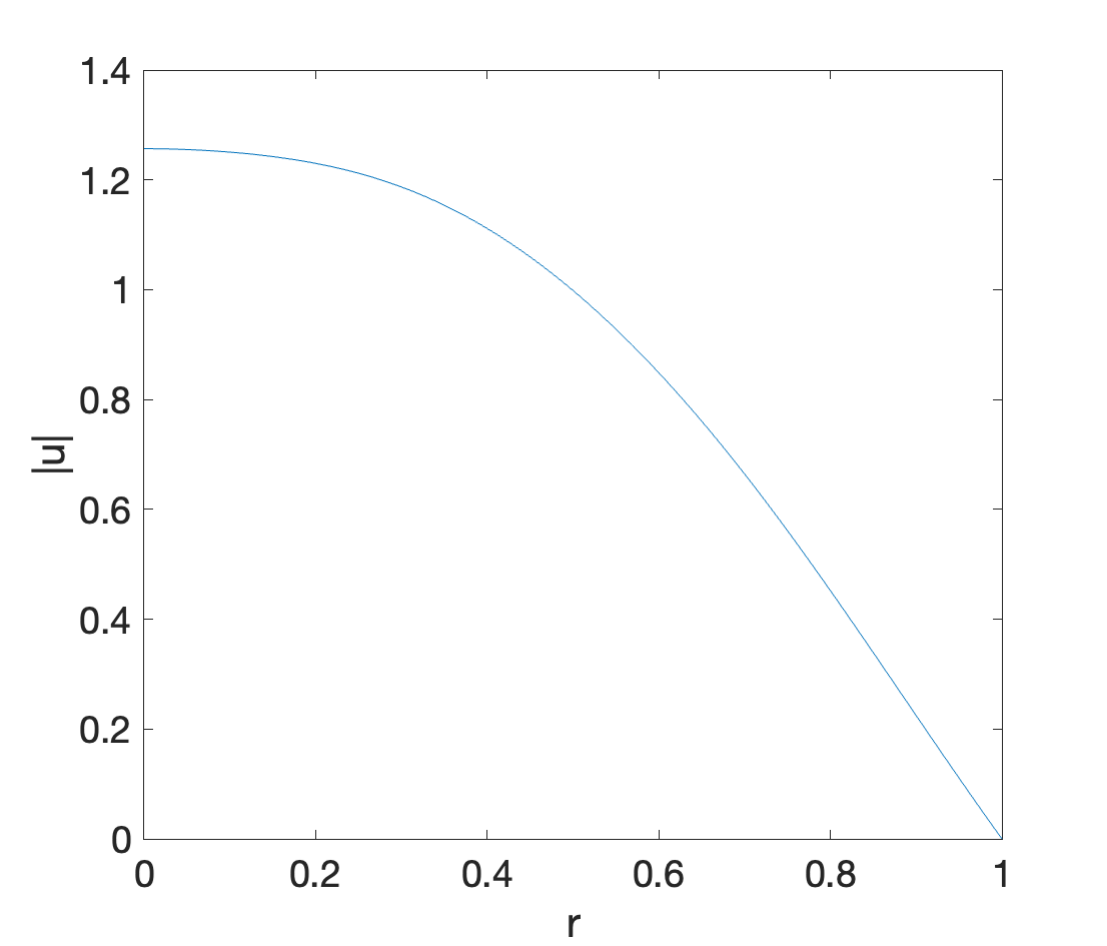} 
\caption{\footnotesize Profile 2 at $t=0.21$.}
\end{subfigure}
\end{subfigure}
\caption{ {\small 2D quintic NLS: recurrent oscillations for a perturbation of the unstable branch of $Q_b$ with $b = 2 > b^\ast \approx -1.03$, $u_0=0.99Q_{2}$. The profiles shown are two coherent radial states visited during the oscillation. } }
\label{F:2D-quintic-water}
\end{figure}
\begin{figure}[htb!]
\includegraphics[width=0.5\textwidth,height=0.28\textwidth]{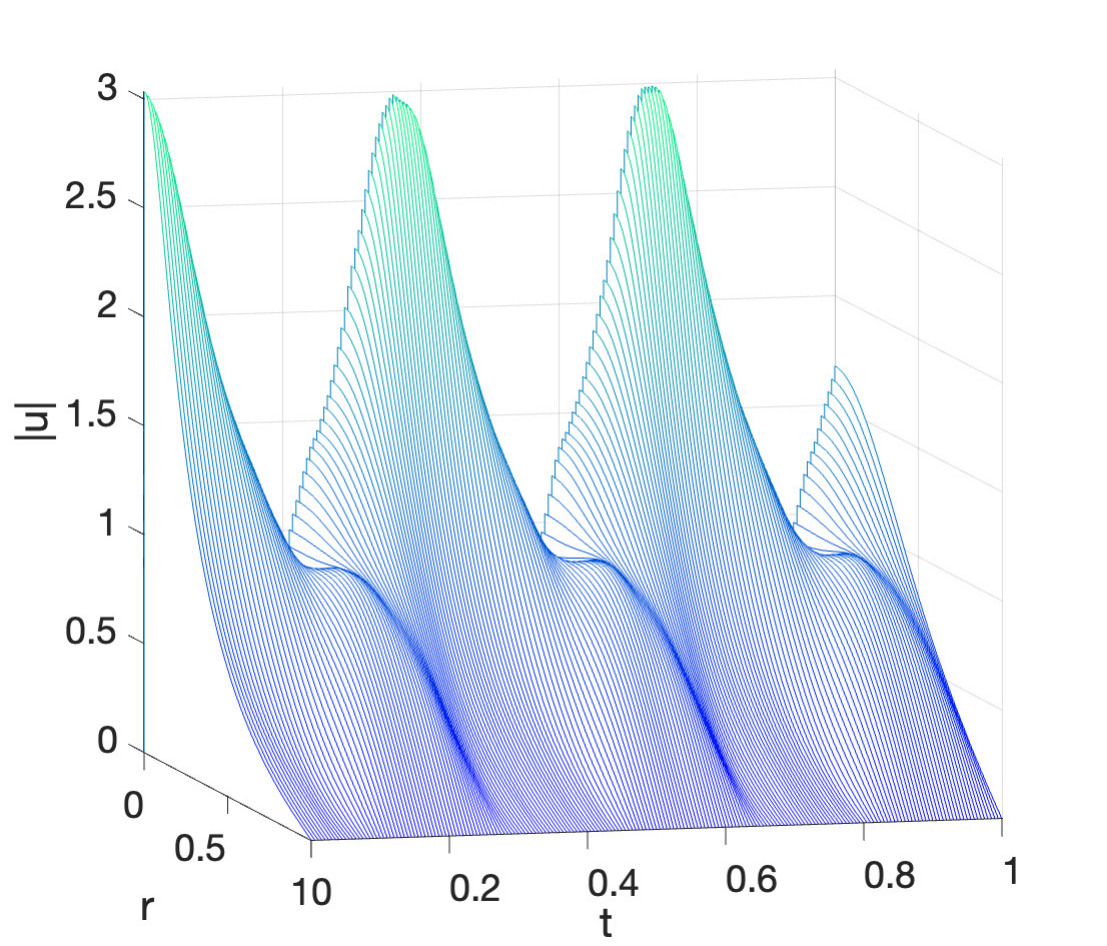}
\includegraphics[width=0.49\textwidth,height=0.28\textwidth]{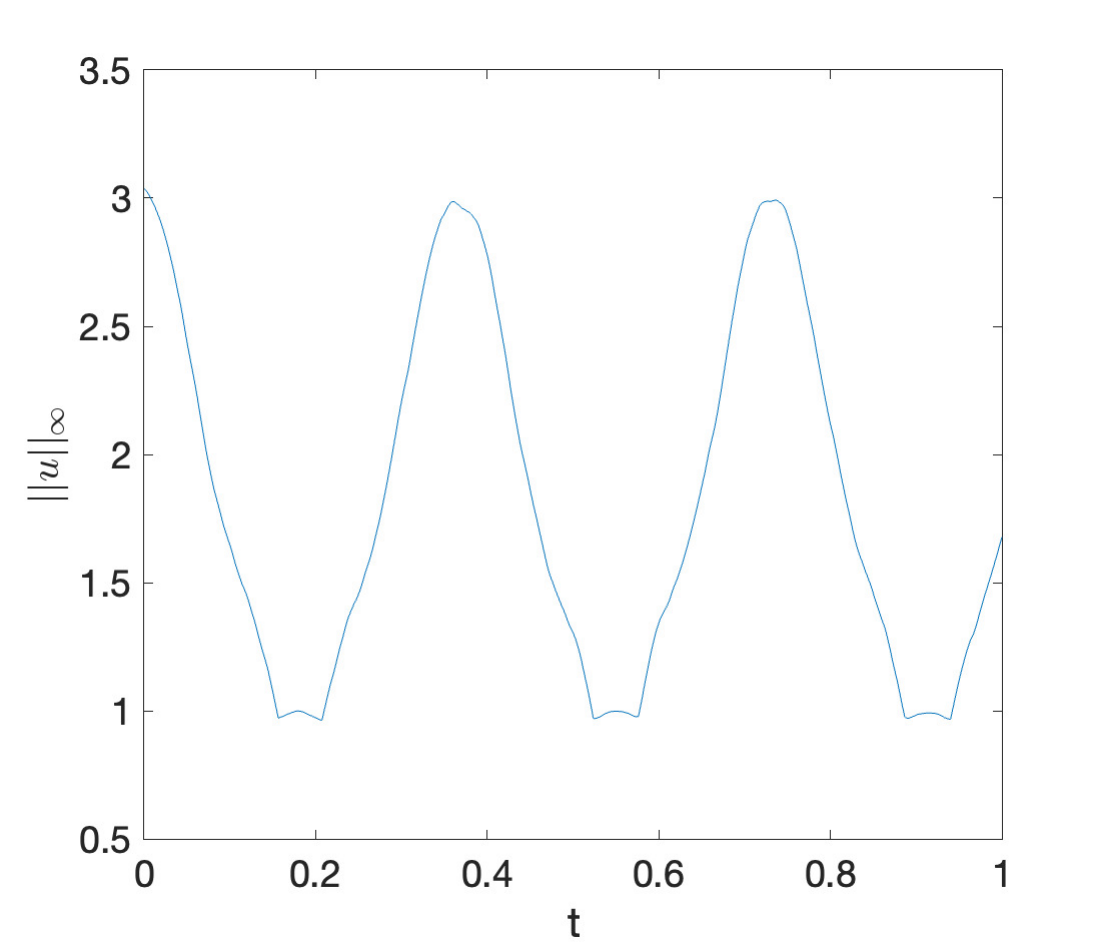}
\caption{ {\small 2D quintic NLS. Unstable branch, $b = 5 > b^\ast$, $u_0=0.99Q_{5}$. Compare oscillations with the case $b=2$.} }
\label{F:2D-quintic-b5}
\end{figure}

The oscillations following perturbation of the unstable branch are shown for $b=2$ in Figure \ref{F:2D-quintic-water}  and for $b=5$ in Figure \ref{F:2D-quintic-b5}, allowing a comparison of the oscillation amplitudes for two values of $b$ on the unstable branch. 
We show the radial profiles of the two states in the case of $b=2$ in plots (C) and (D) of Figure \ref{F:2D-quintic-water}. One may notice that the larger $b$ is, the bigger oscillations are between the two states, and it might happen that there are more states between which solutions may oscillate.

\subsubsection{Blow-up}
We now investigate in more detail the blow-up conditions. If we take a profile $f$ with Dirichlet boundary condition, then we can increase the magnitude of this profile by some positive constant $A$, and choosing $A$ large enough, make its energy negative. Then by a virial argument for negative energy or by scale-invariant mass-energy argument of the second author with Holmer from \cite{HR2007, HR2008}, it would follow that such initial data blow up in finite time. In particular, magnification of $Q_b$ by $A>1$ would blow up in finite time, starting from sufficiently large $b$.  We prove this in the following statement (we assume local theory with mass and energy conservation, the $H^1_0$ continuation and approximation by regular solutions as needed, and by the finite time blow-up we mean that for some positive time $T<\infty$ we have $\ds \lim_{t \nearrow \, T} \|\nabla u(t) \|_{L^2(B_1)}=\infty$).

\begin{proposition}[Blow-up for amplified profiles]
\label{P:blowup}
Let $\alpha \geq 4/d$, $d \geq 1$ (critical or supercritical cases) and let $f \in H_0^1(B_1)$ (and in $L^{\alpha+2}(B_1)$ if $\alpha> 4/(d-2)$, $d > 2$) be nonzero and real-valued. 
Set $u_0=A f$ with $A>0$. 
\begin{itemize}
\item[\underline{Part 1.}](Negative energy case)
The solution $u(t)$ from the initial datum $u_0$ blows up in finite positive time if $A$ satisfies
\begin{equation}\label{E:A-blowup}
A^\alpha > \frac{(\alpha+2)}{2} \frac{\|\nabla f\|_{L^2(B_1)}^2}
{\|f\|_{L^{\alpha+2}(B_1)}^{\alpha+2}}.
\end{equation}
In particular, taking any positive single-peaked profile with the above regularity and Dirichlet boundary conditions, produces a blow-up after multiplying by a sufficiently large amplitude.

\item[\underline{Part 2.}](Non-negative energy case)
Let $\alpha > 4/d$ such that $s_c = \frac{d}2-\frac2{\alpha} \in (0,1)$ (intercritical case). Then the solution with the initial condition $u_0$ and $E(u_0) \geq 0$ blows up in finite time if 
\begin{equation}\label{E:ME<1}
E(u_0) \,M(u_0)^{\frac{1-s_c}{s_c}} < E(\mathcal R) \, M(\mathcal R)^{\frac{1-s_c}{s_c}}
\end{equation}
and
\begin{equation}\label{E:grad>1}
\|\nabla u_0\|_{L^2(B_1)}^{s_c} \|u_0\|_{L^2(B_1)}^{1-s_c}
> \|\nabla\mathcal R\|_{L^2(\mathbb R^d)}^{s_c}
\|\mathcal R\|_{L^2(\mathbb R^d)}^{1-s_c},
\end{equation}
where $\mathcal R = \mathcal R_{\alpha}$.

\item[\underline{Part 3.}](Amplification of a ground state)
Assume $\alpha < 4/(d-2)$ when $d \geq 3$ (and $\alpha \geq 4/d$ as above). In the case of the ground state profile (i.e., $f = Q_b$), for any fixed $A > 1$, there exists $b_0=b_0(A)>0$ such that the initial condition $u_0= A\, Q_b$ evolves into a finite-time blow-up for any $b > b_0(A)$.
\end{itemize}
\end{proposition}

\begin{proof}
We write the energy of the amplified profile $Af$ and use the condition \eqref{E:A-blowup} to deduce
\begin{align} \notag
E(A f) &= \frac{A^2}{2}\|\nabla f\|_{L^2(B_1)}^2
-\frac{A^{\alpha+2}}{\alpha+2} \|f\|_{L^{\alpha+2}(B_1)}^{\alpha+2}\\
& =  \frac{A^2}{\alpha+2} \|f\|^{\alpha+2}_{L^{\alpha+2}(B_1)} \Big( \frac{\alpha+2}{2}\frac{\|\nabla f\|^2_{L^2(B_1)}}{\|f\|^{\alpha+2}_{L^{\alpha+2}(B_1)}} - A^\alpha \Big)<0. \label{E:negE}
\end{align}
Recalling the variance 
$$
V(t) = \int_{B_1} |x|^2 |u(t,x)|^2 \, dx,
$$
differentiating it twice and using the boundary condition, we get the virial identity with Dirichlet boundary (e.g., see \cite[p.\,3 and Sec.\,4.1]{L2022} adjusted to our case)
$$
V''(t) = 16\, E(u_0)- \frac{4(d\alpha-4)}{\alpha+2} \|u(t)\|^{\alpha+2}_{L^{\alpha+2}(B_1)} - 4 \int_{\partial B_1} (x\!\cdot\!n) |\partial_n u|^2 dS \leq 16 E(u_0),
$$
where $n$ is the outward unit normal and $x\!\cdot\!n=1$ on $\partial B_1$ (and in 1D the last term is the sum of the derivatives squared at the two endpoints).
 
Integrating and using that $u_0$ is real-valued, and hence, $V'(0) = 0$, we obtain
$$
V(t) \leq V(0) +8 E(u_0)\, t^2.
$$
Since $V(t)$ is always non-negative, the negative energy condition $E(u_0)<0$, obtained in \eqref{E:negE}, yields a contradiction at some finite positive time $T > \Big(\frac{V(0)}{8|E(u_0)|} \Big)^{1/2}$, which means that blow-up occurs before or at that time $T$, thus, proving Part 1. 

For the second part, assume $\alpha > 4/d$ (such that $0<s_c<1$). 
The argument is essentially the same as in \cite[proof of Thm. 2.1]{HR2007}, applying the sharp Gagliardo-Nirenberg inequality to the zero extension of $u(t)$ and using mass and energy conservation, omitting the repetition, we obtain 
\begin{equation}\label{E:grad-t}
\|\nabla u(t)\|_{L^2(B_1)}^{s_c} \|u_0\|_{L^2(B_1)}^{1-s_c}
> \|\nabla\mathcal R\|_{L^2(\mathbb R^d)}^{s_c}
\|\mathcal R\|_{L^2(\mathbb R^d)}^{1-s_c}
\end{equation}
for any $t$ in the lifespan. Then the virial identity with Dirichlet condition gives (discarding the non-positive boundary term) becomes
$$
V''(t) \leq 4d\alpha E(u_0) - 2(d\alpha-4) \|\nabla u(t)\|^2_{L^2(B_1)} \leq 
4d\alpha \left( E(u_0) - \Big(\frac{M(\mathcal R)}{M(u_0)} \Big)^{\frac{1-s_c}{s_c}} E(\mathcal R) \right)< 0,
$$ 
where we applied the mass-energy threshold \eqref{E:ME<1} as well as the gradient bound from \eqref{E:grad-t}. Then a similar argument as in Part 1 yields the conclusion.

For the third part, for the ground state solutions, we recall \eqref{E:Q-rescale} or Lemma \ref{L:conv}\,(II), stating that $\widetilde Q_b(X) \to \mathcal R$ in $H^1(\mathbb R^d)$ as $b \to \infty$ (with $\widetilde Q_b$ extended by zero outside $B_{\sqrt{b}}$). 
Note that 
$$
E_{B_{\sqrt b}}(A\, \widetilde Q_b) = b^{\frac{d}2-1-\frac2{\alpha}} E_{B_1}(A\, Q_b),
$$
which in the mass-critical case $\alpha=4/d$ becomes $E_{B_{\sqrt b}}(A\, \widetilde Q_b) = b^{-1} E_{B_1}(A\, Q_b)$. 

For the whole-space ground state we apply $\|\mathcal R\|_{L^{\alpha+2}}^{\alpha+2}= \frac{2(\alpha+2)}{d\alpha}\|\nabla\mathcal R\|_{L^2}^2,$ coming from the Pokhozhaev identities (e.g., see \cite[Sec. 2]{L2022}) to get
$$
E(\mathcal R) = \frac{d\alpha-4}{2d\alpha}\| \nabla \mathcal R\|^2_{L^2(\mathbb R^d)} \quad \mbox{and} \quad 
E(A \, \mathcal R) = \frac{A^2}2 \big( 1-\frac4{d\alpha}A^\alpha \big) \|\nabla \mathcal R\|^2_{L^2(\mathbb R^d)}.
$$ 
Thus, in the critical case, $\alpha=4/d$, we obtain 
$$
E(\mathcal R)=0 \qquad \mbox{while} \quad E(A \, \mathcal R)<0, ~~\mbox{provided}~~A>1.
$$ 
Using that $E_{\mathbb R^d}(A \, \tilde Q_b) \to E (A \,\mathcal R)$, we deduce  
$$
b^{-1}E(A Q_b) = E_{\mathbb R^d}(A\, \widetilde Q_b) 
\xrightarrow[b \to \infty]{}  E (A \mathcal R) <0.
$$
Since the limit is negative, we can find large enough $b_0$ such that $E(A \, Q_b) <0$ for all $b > b_0$. 

For the intercritical case, for $A>1$, rescaling mass, gradient and energy, we obtain 
$$
\frac{M(AQ_b)^{\frac{1-s_c}{s_c}}
\|\nabla(AQ_b)\|_{L^2(B_1)}^2}{M(\mathcal R)^{\frac{1-s_c}{s_c}}
\|\nabla\mathcal R\|_{L^2(\mathbb R^d)}^2}
~\to ~ A^{\frac2{s_c}} > 1
$$
and
$$
\quad 
\frac{M(AQ_b)^{\frac{1-s_c}{s_c}}E(AQ_b)}{M(\mathcal R)^{\frac{1-s_c}{s_c}}E(\mathcal R)} ~ \to ~
\frac{d\alpha A^{\frac2{s_c}}-4A^{\alpha+\frac2{s_c}}}{d\alpha-4} < 1.
$$
Thus, \eqref{E:ME<1} and \eqref{E:grad>1} hold for sufficiently large $b$, and a similar argument as in the previous parts concludes the proof. 
\end{proof} 

\begin{remark}
Note that in the critical case because of the boundary term \eqref{E:EMP2} yields $E(Q_b) >0$ (and thus, the exact standing wave $e^{ibt} Q_b$ is global). However, multiplication of the ground state by an amplitude $A>1$ forces negative energy $E(AQ_b) < 0$ for all sufficiently large $b$ as we show in the above proof, and thus, such amplified ground state solution blows up in finite time.  
\end{remark}

\begin{remark}\label{R:2Q1}
Returning back to Figure \ref{F:2D-cubic-threshold}, right plot, we see that the initial condition $u_0 = 2Q_1$ has to blow up in finite time by Part 1 of Proposition \ref{P:blowup}: for the 2D cubic equation, by \eqref{E:P1} we have 
$$
E(2Q_1) =2 \|\nabla Q_1\|^2_{L^2(B_1)}-4 \|Q_1\|^4_{L^4(B_1)} = -2\|\nabla Q_1\|^2_{L^2(B_1)} - 4M(Q_1) < 0,
$$  
where we used the Pokhozhaev identity $\|Q_1\|^4_{L^4(B_1)} = \|\nabla Q_1\|^2_{L^2(B_1)} + \|Q_1\|^2_{L^2(B_1)}$.
\end{remark}

\begin{remark}\label{R:1D-?}
The blow-up observed for $A=1.01$ in Figures \ref{F:1D-septic} and \ref{F:2Dquintic} is consistent with
the large-$b$ amplitude increase result in Part 3 of
Proposition \ref{P:blowup}. However, for the $b$ values used in these figures
($b=5$ and $b=-0.5$, respectively), the initial data have
positive energy as well as the mass-energy condition \eqref{E:ME<1} does not hold, so Part 2  
does not apply. Therefore, these simulations exhibit blow-up outside the explicit
sufficient criteria in Parts 1 and 2 of Proposition \ref{P:blowup}. While this is consistent with a general picture, it remains a question for further investigation. 
\end{remark}

\begin{remark}[About sign-changing profiles]\label{R:sign}
We note that the condition $M(u_0) > M(\mathcal R)$ alone on the bounded domain with Dirichlet boundary condition does not imply blow-up.
For example, in 1D with $\alpha=4$, we can take $Q_0$ to be the positive, even solution of
$$
-Q_0''-Q_0^5=0,\qquad Q_0(-1)=Q_0(1)=0,
$$
and define
$$
W(x)=
\begin{cases}
-\sqrt2\,Q_0(2x+1),& -1 \leq x \leq 0,\\
\sqrt2\,Q_0(2x-1),& 0 \leq x \leq 1,
\end{cases}
$$
which is continuous at $x=0$ with its first and second derivatives. 
Then $-W''-|W|^4W=0$ and $W(-1)=W(1)=0$, so $u(t,x) = W(x)$ is a time-independent global solution. 

To compute mass, we multiply by $Q_0'$ the above equation for $Q_0$ and integrate, which gives $(Q'_0)^2 = \frac13 (Q_0(0)^6 - Q_0^6)$. Similar argument as in Prop. \ref{P:1D-b*} is applied here: changing variables to $t = Q_0(x)$ on $(0,1)$ and then to $s = t/Q_0(0)$, we obtain
$$
M(Q_0) = 2\sqrt{3} \int_0^1 \frac{s^2ds}{\sqrt{1-s^6}} = \frac{\pi}{\sqrt 3}.
$$ 
Then $M(W) = 2M(Q_0) = \frac{2\pi}{\sqrt 3}$ and recalling the mass of $\mathcal R$ from \eqref{E:MR} (with $\alpha=4$), we get $M(\mathcal R) = \frac{\sqrt{3} \, \pi}{2}$. Thus, $M(W) > M(\mathcal R)$, but $u = W$ as mentioned above does not blow up. (Note that it has positive energy $E(W) = \frac13 \|W'\|^2_{L^2(-1,1)}$, since $\|W'\|^2_{L^2} = \|W\|^6_{L^6}$.)
\end{remark}

In summary, Proposition \ref{P:blowup} proves Part II\,(ii) under the negative energy assumption and it also verifies the blow-up statement for $A>1$ and $b>b_0(A)$ in Part III\,(iii) of Conjecture \ref{C:1} (note we were not able to show that $b_0(A) = b^*$), Remark \ref{R:sign} explains why the mass assumption alone is insufficient.



\subsubsection{Small-amplitude}
\begin{figure}[htb!]
\includegraphics[width=0.5\textwidth,height=0.27\textwidth]{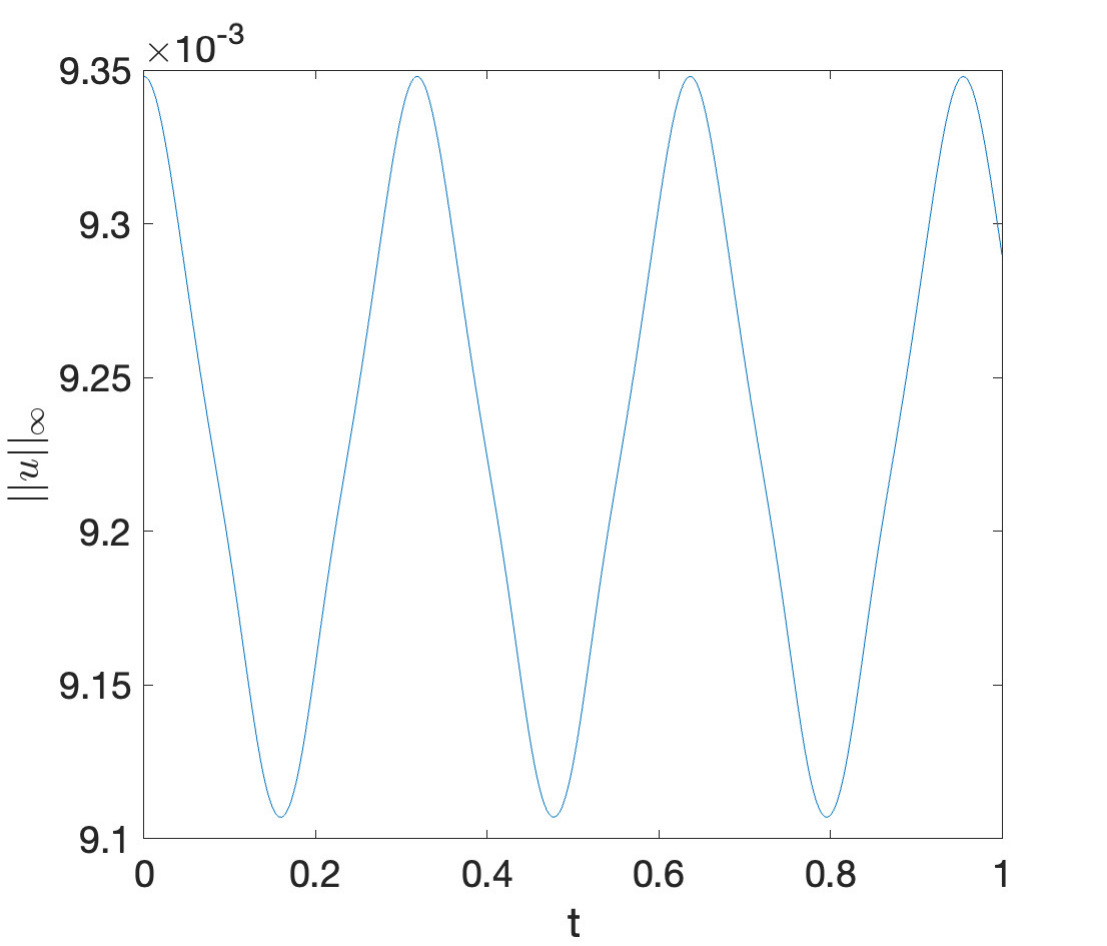}
\caption{ {\small 1D quintic NLS, $\alpha=4$. Small amplitude initial data  $u_0=0.01Q_{-2}$, $b=-2$, resolves into an oscillation around a small soliton, consistent with soliton resolution.} }
\label{F:small}
\end{figure}

We conclude with investigations of radiation, or more precisely, of small-amplitude dynamics, for which our simulations did not show any dispersion or radiation in such solutions and instead exhibited coherent oscillations. An example in Figure \ref{F:small} shows that the initial condition $u_0 = 0.01Q_{-2}$ does not disperse and evolves into oscillations around some small soliton or coherent state (recall that the masses of ground states $Q_b$ cover the entire interval $(0,M(\mathcal{R}))$), note the scale of the $y$-axis (on the order of $10^{-3}$). This, therefore, provides numerical evidence for the proposed soliton resolution conjecture on bounded domains: any such solution resolves asymptotically into a sum of persistent coherent oscillations.

\section{Conclusions}\label{S:conclusions} 
In this paper, we presented an analytical and  numerical study of least-energy ground states for the nonlinear Schr\"odinger equation on the unit 
ball in one and two dimensions with Dirichlet boundary condition. Our computations illustrate the stabilizing role of confinement and give high-accuracy numerical profiles for ground states that are not available in closed form. In the subcritical and critical cases considered here, the mass curves have positive slope throughout the computed parameter range, in agreement with the known orbital-stability theory and with the time evolutions of ground state perturbations.  

In the critical case, the mass of the ground state 
$M(Q_{b})$ increases towards the mass of the whole-space ground state 
$\mathcal{R}$. While the ground states $Q_{b}$ are stable under small 
perturbations in the critical case, our computations indicate stable dynamics for finite perturbations as long
as the mass of the initial data remains below the mass of $\mathcal{R}$. Above this threshold, the negative energy solutions blow up in finite time as we prove in Proposition \ref{P:blowup}.

In the supercritical case, the new feature is the appearance of a turning point $b^*$ in the curve $b \mapsto M(Q_b)$ and the corresponding splitting of the ground state family into two branches. 
In 1D, we prove in Proposition \ref{P:1D-b*} existence of a mass maximum at a turning point and characterize the location of the turning point via an integral equation. 
The phenomenon of {\it branching} splits the ground state family into stable and unstable branches, depending on the size of the parameter $b$. Perturbations of ground states on a stable branch leave the 
solutions oscillating around that ground state, while perturbations 
of the unstable branch force solutions either to blow  up in finite 
time (if perturbations have an amplitude larger than the solitary 
state) or to oscillate between two coherent states (if 
perturbations have amplitude smaller than the original ground 
state). In particular, we establish blow-up conditions in Proposition \ref{P:blowup} for $A\,Q_b$ with $A>1$ and sufficiently large $b$ in the mass-critical and intercritical cases (including conditions with positive energy) and its generalization for an arbitrary single-peaked real-valued profile $f$. Furthermore, we observe that this problem does not produce any radiation or scattering (due to reflection at the Dirichlet boundary), and thus, our simulations support a bounded-domain soliton-resolution conjecture: no scattering or radiation holds for all data, splitting any solution into coherent structures, such as solitons, even for very small initial data.  

\bigskip

{\bf Conflict of Interest:} The authors declare that they have no
conflicts of interest.

\bibliographystyle{abbrv}
	
\bibliography{references}

\end{document}